\documentclass[12pt,leqno]{article}
\usepackage[latin1]{inputenc}
\usepackage{amsxtra}
\usepackage{amsmath}
\usepackage{amssymb}
\usepackage{amsfonts}
\usepackage[all]{xy}
\usepackage{mathrsfs}
\usepackage{amsthm}
\usepackage{esint}
\usepackage{enumerate}

\newtheorem{thm}[equation]{Theorem}
\newtheorem{cor}[equation]{Corollary}
\newtheorem{lem}[equation]{Lemma}
\newtheorem{prop}[equation]{Proposition}

\newtheoremstyle{example}{\topsep}{\topsep}%
     {}
     {}
     {\bfseries}
     {.}
     {2pt}
     {\thmname{#1}\thmnumber{ #2}\thmnote{ #3}}

   \theoremstyle{example}
   
   \newtheorem{Defi}[equation]{Definition}
   \newtheorem{rem}[equation]{Remark}
   \newtheorem{rems}[equation]{Remarks}

   \newtheorem{ex}[equation]{Example}
    \newtheorem{ques}[equation]{Question}
   
  \numberwithin{equation}{subsection}

\newtheorem{exa}[equation]{Example}

\def\la{{\langle}}
\def\ra{{\rangle}}

\def\MM{\mathbb{M}}

\def\RR{\mathbb{R}}
\def\ZZ{\mathbb{Z}}
\def\QQ{\mathbb{Q}}
\def\TT{{\mathbb{T}}}

\def\fen{\mathfrak{f}}
\def\gen{\mathfrak{g}}
\def\hen{\mathfrak{h}}

\def\ken{\mathfrak{k}}

\def\Gg{\mathfrak{g}}

\def\Ac{\mathcal{A}}

\def\Cc{\mathcal{C}}

\def\Mc{\mathcal{M}}

\def\Pc{\mathcal{P}}

\def\Sc{\mathcal{S}}

\def\Gf{\mathfrak{f}}

\def\ab{{\operatorname{ab}\nolimits}}

\def\Ad{\operatorname{Ad}\nolimits}
\def\Aut{\operatorname{Aut}\nolimits}

\def\Bar{\operatorname{Bar}}

\def\Cat{\operatorname{Cat}\nolimits}
\def\CcM{{\Cc\Mc}}
\def\cl{\operatorname{cl}\nolimits}
\def\CM{{\operatorname{CM}\nolimits}}
\def\Coker{\operatorname{Coker}\nolimits}

\def\deg{{\operatorname{deg}\nolimits}}
\def\der{\operatorname{\it \mathscr{D}\kern-.25em er}}
\def\Der{\operatorname{Der}\nolimits}
\def\dim{{\rm{dim}}}
\def\DR{{\operatorname{DR}\nolimits}}

\def\End{\operatorname{End}\nolimits}

\def\fd{{\operatorname{fd}\nolimits}}
 \def\FL{\operatorname{FL}\nolimits}

 \def\Glb{\widetilde{\Pi}}
 
 \def\Glob{\on{Glob}}

 \def\Hom{\operatorname{Hom}\nolimits}
 
 \def\id{\operatorname{id}\nolimits}
\def\Id{\operatorname{Id}\nolimits}
 \def\Im{\operatorname{Im}\nolimits}
 
 \def\Ker{\operatorname{Ker}\nolimits}

 \def\Lie{\operatorname{Lie}\nolimits}

 \def\Mor{{\operatorname{Mor}\nolimits}}
 
 \def\Ob{\operatorname{Ob}\nolimits}
 \def\on{\operatorname}

  \def\pt{\operatorname{pt}\nolimits}

 \def\sab{{\operatorname{sab}\nolimits}}
 \def\Sh{\operatorname{Sh}}

 \def\Vect{{\operatorname{Vect}\nolimits}}

\def\1{{\bf 1}}
\def\lra{\longrightarrow}

\def\(({(\hskip -1mm (}
\def\)){)\hskip -1mm )}
\def\be{\begin{equation}}
\def\ee{\end{equation}}
\def\ed{\end{document}}
\def\gdot{{ \gamma{^{^{\hskip -.17cm \bullet}}}}}
\def\lan{\langle}
\def\llan{{\lan\hskip -.1cm \lan}}
\def\rran{\rangle\hskip -.1cm\rangle}
\def\loint{\varointctrclockwise}
\def\roint{\varointclockwise}

\title{Membranes and higher groupoids}
\author{M. Kapranov
  }

\begin{document}


\maketitle

  \tableofcontents

  \addtocounter{section}{-1}

\section {Introduction.}

The goal of this paper is threefold.

\vskip .3cm

\noindent {\bf (0.1)} First, we generalize the classical result of C. Reutenauer 
(\cite{reutenauer}, \S 8.6.12) on the abelianization of the commutant of the free Lie algebra.
Denote by $\FL(V)$ the free Lie algebra on a finite-dimensional vector space $V$
over a field $k$ of characteristic $0$.
Reutenauer's theorem  (in its dual formulation)
says that $H^1([\FL(V), \FL(V)], k)$ is identified with the
space of formal germs of closed 2-forms on $V$, see formula \eqref{eq:reutenauer2}
below.
 Our generalization 
(Theorems \ref{thm:semiabelianization} and  \ref{thm:cartesian}), uses a certain free dg-Lie algebra
$\fen^\bullet(V)$, containing $\FL(V)$ as the degree $0$ part, and
(if formulated in the dual version similar to the above)
produces the full de Rham complex of formal germs of forms,
 truncated in degrees $\geq 2$. 
Instead of abelanization of the commutant, we use the so-called
semiabelianization of $\fen^\bullet(V)$.  We call a dg-Lie algebra
$\gen^\bullet$, situated in degrees $\leq 0$,
{\em semiabelian}, if the commutators $\gen^i\otimes\gen^j\to\gen^{i+j}$
vanish for $i,j\leq -1$. 

\vskip .2cm

By Quillen's theory, dg-Lie $\QQ$-algebras situated in
degrees $\leq 0$, correspond to all rational homotopy types.
The class of semiabelian 
dg-Lie algebras corresponds to those homotopy types that are represented
by strict $\infty$-groupoids of Brown-Higgins
\cite{brown-higgins}.  Indeed, for such homotopy types the Whitehead
products $\pi_{i+1}\otimes\pi_{j+1}\to\pi_{i+j+1}$, vanish for  $i,j\geq 1$. 
On the other hand,  by Grothendieck's philosophy, all homotopy types should
correspond to appropriately defined weak $\infty$-groupoids. The functor
of taking the maximal semiabelian quotient
of a dg-Lie algebra $\gen^\bullet\mapsto \gen^\bullet_{\sab}$
is therefore the analog of  that of ``strictification",  from the
category of ``true" $\infty$-groupoids to that of Brown-Higgins ones. 
Our result can be thus seen as calculating the (Lie algebra
version of the) nonabelian derived
functor of this strictification procedure.

\vskip .3cm

\noindent {\bf (0.2)} Our second goal is to use the theory of 2-dimensional
holonomy of connections with values in crossed modules of Lie
groups and algebras, as developed
by Baez-Schreiber \cite{baez-schreiber}, in order to construct a representation
of 2-dimensional membranes in $\RR^n$ by certain formal series-type data.
This representation extends the classical construction of K.-T. Chen 
which associates to a path $\gamma$ in $\RR^n$ the noncommutative 
formal power series
$$
E_\gamma(Z_1, ..., Z_n)\,\,\, = \,\,\,\sum_{d=0}^\infty \sum_{i_1, ..., i_d=0}^n
Z_{i_1} ... Z_{i_d} \int_\gamma (dx_{i_1}, ..., dx_{i_d}),
$$
the generating function for all the iterated integrals of the coordinate 1-forms $dx_i$.
This series, being group-like, lies in $\widehat G^0_n$, the free prounipotent proalgebraic group on $n$ generators.
We extend $\widehat G^0_n$ to a crossed module of prounipotent groups
$G^{\geq -1}_n= \{G^{-1}_n\buildrel\partial\over \to G^0_n\}$, which carries a natural connection over $\RR^n$
with vanishing fake curvature in the sense of Breen-Messing \cite{breen-messing}.  
The Lie algebra crossed module of $G^{\geq -1}_n$ is a natural completion of the
truncation
  of $\fen^\bullet(V)$ in degrees $\geq -1$. To each 2-brane $\sigma$ we then associate 
  \eqref{eq:2-chen} an
  element $\widehat M(\sigma)\in G^{-1}_n$, with $\partial \widehat M(\sigma)\in G^0_n$ being the Chen series
  corresponding to the boundary path of $\sigma$.

 \vskip .3cm
 
 \noindent {\bf (0.3)} Third, we generalize the construction of Baez-Schreiber
 \cite{baez-schreiber} to that of
 $p$-dimensional holonomy for connections with values in {\em crossed complexes} of Lie groups. 
Crossed complexes correspond to Brown-Higgins $\infty$-groupoids with one object
 and consist of a non-abelian ``head" in degrees $0$ and $-1$ attached to an  abelian
 ``tail" in lower degrees.   In this paper we consider only what should be understood as connections in
 trivial $n$-gerbes with a given structure $n$-groupoid. 
  Our construction is based on the same techniques
 as in \cite{baez-schreiber}, namely a covariant generalization of Chen's iterated
 integrals in the presence of a background connection, see \cite{hofman}. 
 It gives a strict functor from an appropriately defined $n$-groupoid of 
 unparametrized membranes in  a manifold $X$ to the $n$-groupoid corresponding
 to the crossed complex, see Theorem \ref{thm:p-holonomy} .
 
  Although Brown-Higgins $n$-groupoids are rather
 restrictive from the general homotopy-theoretic point of view, this seems to be the maximal
 generality in which we have a holonomy corresponding to a
 $p$-brane as a {\em geometric object}, without any extra data such as a choice of
 a parametrization  or  of a slicing into lower-dimensional membranes. 
 A more general construction  given by E. Getzler \cite{getzler}, depends on
 such choices but comes with a system of higher homotopies accounting for
 making a different choice. The usual homotopy-theoretic source of
 Brown-Higgins $n$-groupoids is provided not by single spaces but by filtered spaces
 \cite{brown:new}. Considering thin homotopies (reparametrizations, cancellations etc.)
 to pass from parametrized to unparametrized (``geometric") 
  membranes can be seen as introducing a differential-geometric analog
 of the skeletal filtration.

 The place  of Brown-Higgins $n$-groupoids 
 among more general  homotopy types (strict 2-dimensional associativity, 
  more solid geometric nature of higher holonomy with values in  them) 
 is somewhat similar to the place of commutative rings among all associative rings. 
 It may be therefore interesting to study infinitesimal non-Brown-Higgins deformations
 of Brown-Higgins $n$-groupoids in geometric terms, similarly to noncommutative
 deformations of commutative algebras.

 \vskip .2cm
 
 We then apply the theory of $p$-dimensional holonomy to extend the
 representation of 2-branes above to the case of $p$-branes in $\RR^n$ with
 arbitrary $p$. For $p$-branes $\Sigma$ whose boundary consists of a point, this
 representation associates to $\Sigma$ the corresponding de Rham current
 (integration functional on polynomial $p$-forms). In general, this is
 a certain nonabelian twist of such de Rham representation.
 This construction can be seen as a nonlinear and nonabelian
 analog of Grassmann's geometric calculus:
we represent $p$-branes in $\RR^n$ by data constructed out of exterior powers $\Lambda^i(\RR^n)$,
$i\leq p$. 
 
 \vskip .3cm
 
 \noindent {\bf (0.4)}  A different generalization of Reutenauer's theorem
 was found  by B. Feigin and
B. Schoikhet  \cite{feigin-shoikhet}. Their result
provides  a relation of the free associative, not Lie, algebra with $\Omega^{2\bullet, \cl}$,
the space of closed differential forms of arbitrary even degrees.  Our Theorem \ref{thm: abelianization}
provides a generalization in a different direction, involving forms of all, not just even,  degrees.

It seems that the dg-Lie algebra $\fen^\bullet(V)$, associated to a vector space $V$,
admits  a ``curvilinear" analog 
which is a natural  dg-Lie agebroid $\Pc_X^\bullet$ associated to any smooth (or complex
or algebraic) manifold $X$ and containing the free Lie algebroid $\Pc_X$
from \cite{kapranov-Lie-algebroids} as the  degree 0 part. 
 We discuss
this possibility in more detail in Remark \ref{rem:nonlin}. Note that nonlinear analogs of
the Feigin-Shoikhet result have been subject of several recent works
\cite{dobr1, dobr2, jordan-orem}.

An approach to 2-dimensional holonomy  somewhat different from that of J. Baez and U. Schreiber
(and thus from the approach of the present paper) was developed
by A. Yekutieli \cite{yekutieli}. That approach is based on direct approximation by
``Riemann products", not on iterated integrals.  It can probably be applied to the
$p$-dimensional holonomy (situation of Chapter \ref{sec:crossed-complex-n-branes})
as well.  An approach  involving the holonomy of connections on the
space of paths (and thus using iterated integrals, albeit implicitly) was developed by
J. F.  Martins and R. Picken \cite{martins-picken}. 

\vskip .3cm

 \noindent {\bf (0.5)} The three chapters of the paper correspond to the three main subjects
 outlined above. In the Appendix we provide a self-contained account
 of the main points of the theory of
 covariant iterated integrals. 
 
 \vskip .2cm
 
   I would like to thank P. Bressler, R. Mikhailov, B. Shoikhet and A. Yekutieli
 for useful discussions and/or correspondence. This work was
 supported by   the World Premier International Research Center Initiative (WPI Initiative), MEXT, Japan
 and by the Max-Planck-Institut f\"ur Mathematik, Bonn, Germany.

\vfill\eject

\section { A resolution of an abelian Lie algebra.}

\subsection { The resolution.}\label{subsec:resolution}
 Let $k$ be a field of characteristic
0. We denote by $\Vect_k^{\ZZ}$ the symmetric monoidal category of 
$\ZZ$-graded $k$-vector spaces $V^\bullet =\bigoplus_{i\in\ZZ} V^i$, with the usual
graded tensor product and the symmetry $V^\bullet \otimes W^\bullet\to W^\bullet\otimes V^\bullet$ 
given by the Koszul sign rule:
\[
v\otimes w \longmapsto (-1)^{\deg(v)\deg(w)} w\otimes v. 
\]
By a {\em graded Lie k-algebra} we mean a Lie algebra object $\gen^\bullet$ in $\Vect_k^{\ZZ}$,
so the bracket on $\gen^\bullet$ satisfies
 the antisymmetry and the Jacobi identity twisted, in the usual way, by the Koszul sign factors. 
In particular, for $V^\bullet\in \Vect_k^{\ZZ}$ we have the {\em free graded Lie algebra}
$\FL(V^\bullet)$. By a {\em dg-Lie $k$-algebra} we mean a graded Lie $k$-algebra
equipped with a differential $d$ of degree $+1$, satisfying the graded Leibniz rule and such that
$d^2=0$.  

\vskip .2cm

Let $V$ 
 be a finite-dimensional $k$-vector space (concentrated in degree $0$). Let $Z_1, ..., Z_n$ be
a basis of $V$. We will write $\FL(Z_1, ..., Z_n)$ for
the free Lie $k$-algebra $\FL(V)$. 
For each nonempty subset $I=\{ i_1 < ... < i_p\}$ of $\{1, ..., n\}$ we
introduce a generator $Z_I = Z_{i_1, ..., i_p}$ of degree $(-p+1)$.
In  the case $p=1$ these generators will be identified with the $Z_i$ above.
Let  $\fen^\bullet = \fen^\bullet (Z_1, ..., Z_n)$ be
 the free graded Lie $k$-algebra on all
the generators $Z_I$. The degree 0 part of this algebra
is $\FL(Z_1, ..., Z_n)$. 
One can express this more invariantly by saying that 
$\fen^\bullet = \fen^\bullet(V)$
is constructed  from $V$ as follows:
 \be 
 \fen^\bullet(V) = \FL\bigl (\Lambda^{\geq 1}(V)\bigr), 
\quad \deg(\Lambda^p(V)) = -p+1.
\ee
Indeed, the $Z_I$ for $|I|=p$, form a basis in $\Lambda^p(V)$. 

\vskip .2cm

We introduce a derivation  $d$ of $\fen^\bullet$ of degree 1 by defining
it on the generators (and then extending uniquely using the Leibniz rule)
as follows:
\be\label{eq:differential-f}
dZ_I ={1\over 2}  \sum_{I=J\sqcup K} \sigma(J,K) [X_J, X_K]. 
\ee
Here the summation is over all partitions of $I$ into the union of two
disjoint subsets $J,K$. The number $\sigma(J,K)$ is the sign of the
shuffle  permutation induced by this partition. For example,
\be dZ_{ij}= [Z_i, Z_j], \quad dZ_{ijp} = [Z_i, Z_{jp}] - [Z_j, Z_{ip}] + 
[Z_p, Z_{ij}].  
\ee

\begin{prop} \label{prop:f}  (a) The differential $d$ satisfies $d^2=0$
and thus makes $\fen^\bullet$ into a dg-Lie algebra.\hfill\break
(b) The dg-Lie algebra $\fen^\bullet(V)$ is a resolution of the abelian
Lie algebra $V$, i.e., 
$$H^0(\fen^\bullet(V)) = V, \quad H^i(\fen^\bullet(V))=0, \,\, i\neq 0.$$

\end{prop}

\noindent {\sl Proof:} Part (a) can be verified directly on the generators. 
An alternative way to see the validity of (a) and also to prove (b) is
to consider the graded commutative algebra $\Lambda$ without unit defined
as $\Lambda= \Lambda^{\geq 1}(V^*)$, with $\deg(\Lambda^p(V^*)) = p$. 
Then
$$\fen^\bullet(V) = \FL(\Lambda^*[-1]) =
 {\operatorname{Harr}}^\bullet(\Lambda)$$
is nothing but the Harrison cochain complex of $\Lambda$. See \cite{barr}
for background on Harrison (co)homology. To be precise, the
standard definition of  the Harrison chain
complex of a commutative algebra $A$,  is as the space of indecomposable elements of
the Hochschild complex $\bigoplus_r A^{\otimes r[r]}$ with respect to the shuffle
multiplication \cite{barr}. This space
  is the same as $\FL(A[1])$. The Harrison cochain complex
is dual: it is $\FL(A^*[-1])$. 
 The differential
\eqref{eq:differential-f} is just the map dual to the multiplication in $\Lambda$,
so it is indeed the Harrison differential, and the fact
that $d^2=0$ follows on general grounds. 
 Further, $\Lambda$ is free as a graded
commutative algebra without unit, and $V$ is its space of generators.
 So part (b)
follows from vanishing of the higher
 Harrison cohomology of a free (graded) commutative
 algebra.    \qed

\vskip .1cm

Taking the universal enveloping algebra of $\fen^\bullet(V)$, i.e., the
free associative dg-algebra on the $Z_I$, we get a free associative dg-resolution
of $S^\bullet (V) = k[Z_1, ..., Z_n]$. Compare with \cite{feigin-shoikhet}.

\subsection { The connection.} \label{subsec: connection} 

Let $X$ be a $C^\infty$-manifold and $\gen^\bullet$ a dg-Lie $\RR$-algebra.
By a {\em graded connection} on $X$ with values in $\gen^\bullet$ we will
simply mean a differential form $A\in (\Omega^\bullet_X\otimes\gen^\bullet)^1$
of total degree 1. We associate to $A$ the operator $\nabla_A=d-A$ in
the de Rham complex $\Omega^\bullet_X\otimes\gen^\bullet$, 
 where $d = d_{\operatorname {DR}}
+ d_{\gen}$ is the sum of the de Rham differential and the differential
induced by the one in $\gen^\bullet$. In particular, $\nabla_A^2$ is given
by multiplication with the {\em curvature} of $A$ which is the form
\be\label{eq:graded-curvature}
F_A\,\,=\,\, dA - {1\over 2} [A, A] \,\,\in\,\, (\Omega^\bullet_X\otimes\gen^\bullet)^2.
\ee
We say that $A$ is {\em formally flat}, if $F_A=0$.

We apply this to the case when $X=V = \RR^n$ with basis $Z_1, ..., Z_n$
and the corresponding coordinate system
 $t_1, ..., t_n$, and $\gen^\bullet = \fen^\bullet=\fen^\bullet(V)$.
  Introduce a differential form $A$ on $V$
with values in $\fen^\bullet$ as follows:
\be\label{eq:the-connection} 
A = \sum_{p=1}^n \sum_{1\leq i_1< ...< i_p\leq n} Z_{i_1, ..., i_p} 
dt_{i_1} ... dt_{i_p}. 
\ee
Note that $A$ has total degree 1, so it is a graded connection with values in $\fen^\bullet$.

 \begin{prop}The  graded connection  $A$
 is formally flat. 
\end{prop}

\noindent {\sl Proof:} Let us write $dt_I = dt_{i_1} ... dt_{i_p}$. Our form
$A = \sum Z_I dt_I$ has constant coefficients. 
So $d_{\operatorname{DR}} (A)=0$. On the other hand,
$$\begin{gathered}
d_{\fen}(A) \,\,= \,\,\sum_I d(Z_I) dt_I = {1\over 2}
 \sum_{J, K: \,\, J\cap K =\emptyset} \sigma(J,K) [Z_J, Z_K] dt_{J\sqcup K} \,\,= \cr
=\,\,{1\over 2} \sum [Z_J, Z_K] dZ_J dZ_K\,\, =\,\, {1\over 2} [A,A]. \qed
\end{gathered} $$

It is clear moreover, that  requiring $A$ to be formally flat, is equivalent to imposing the differential
\eqref{eq:differential-f}.  Therefore
$A$ is the {\em universal translation equivariant flat graded connection} on $V$,
and so $\fen^\bullet$ has the intepretation as
 the structure dg-Lie algebra of this universal graded connection.

 \begin{rem}\label{rem:nonlin}
 It would be very interesting to construct a ``curvilinear" version of the dg-Lie algebra $\fen^\bullet(V)$.
Recall that in \cite{kapranov-Lie-algebroids}, to each smooth manifold $X$
(in $C^\infty$, analytic or algebraic category) we associated a Lie algebroid  $\Pc_X\buildrel\alpha\over\to T_X$
(in the corresponding category)
which locally ``looks like" $\FL(T_X)$ but with the Lie algebra structure mixing the Lie bracket of vector
fields and the formal bracket in the free Lie algebra. The free Lie algebra $\FL(\RR^n)$  can be recovered
as the algebra
of global translation invariant sections of $\Pc_{\RR^n}$. 
One can view $\Pc_X$ as the Lie algebroid corresponding to
the (infinite-dimensional Lie) groupoid of formal unparametrized paths. In particular, bundles with (not necessarily flat) connections
can be seen as modules over $\Pc_X$. 

In this direction, it would be natural to extend the construction of  $\Pc_X$ by associating with
each manifold $X$ as above,  a dg-Lie algebroid $\Pc_X^\bullet$
situated in degrees $\leq 0$ with $\Pc_X^0=\Pc_X$, so that $\fen^\bullet(\RR^n)$ is recovered
as the dg-Lie algebra of global translation invariant sections of $\Pc^\bullet_{\RR^n}$. One can possibly
define $\Pc^\bullet_{\RR^n}$ by a version of the Tannakian formalism as in
 \cite{kapranov-Lie-algebroids}, i.e., by making it to
classify (flat) graded connections. In this way the role of $\fen^\bullet(\RR^n)$ in describing
translation invariant graded connections would be extended to the curvilinear case. 
We leave this study for future work.  

\end{rem}

 \subsection { Semiabelianization.}
  Let $k$ be any field
of characteristic 0. A dg-Lie $k$-algebra $\gen^\bullet$
is called abelian, if its bracket vanishes. Thus an abelian dg-Lie algebra is the
same as a cochain complex. For any dg-Lie $k$-algebra $\gen^\bullet$ we denote by
$
\gen^\bullet_{\ab} = \gen^\bullet /[\gen^\bullet , \gen^\bullet] 
$
its maximal abelian quotient.

Further, let  $\gen^\bullet$ be a
dg-Lie $k$-algebra situated in degrees $\leq 0$. We will say that
$\gen$ is {\em semiabelian}, if $[\gen^{\leq -1}, \gen^{\leq -1}] = 0$. Thus
the data needed to define a semiabelian dg-Lie algebra are:

\begin{enumerate}
\item[(SA1)] A Lie algebra structure on $\gen^0$.

\item[(SA2)]  A $\gen^0$-module structure on each $\gen^i$, taken to be the adjoint action for $i=0$
and denoted by $[z,x]$, $z\in \gen^0$, $x\in\gen^i$. 

\item[(SA3)]  A differential $d$ of degree $+1$ such that $d^2=0$.
\end{enumerate}

\begin{prop}\label{prop:semiabelian-conditions}
 In order that the above data define a semiabelian dg-Lie algebra, it is necessary
and sufficient that the following hold:

(a) $d$ is a morphism of $\gen^0$-modules.

(b) For $x,y\in \gen^{-1}$ we have $[dx,y]=[x,dy]$. 

\end{prop}

\noindent {\sl Proof:} Indeed, the data (SA1) and (SA2) completely define the structure of a graded Lie algebra
such that $[\gen^{\leq -1}, \gen^{\leq -1}] = 0$, and all we need is to account for the Leibniz rule for
$d[x,y]$. The condition (a) corresponds to the case $x\in \gen^0$, $y\in \gen^i$, $i<0$,
while (b) corresponds to the case $x,y\in\gen^{-1}$, as $d[x,y]$ in this case
must be equal to 0. Other cases do not present any restrictions. \qed

\vskip .2cm

\begin{ex}[Crossed modules]\label{ex:Lie-crossed} Suppose that
 $\gen^\bullet$ is situated
in degrees $0, -1$ only:
\be\label{eq:crossed-Lie}
\gen^\bullet \quad = \quad \bigl\{ \gen^{-1} \buildrel d\over
\longrightarrow \gen^0 \bigr\}. 
\ee
Then $\gen^\bullet$ is semiabelian, as for any $x, y\in \gen^{-1}$
we have $[x, y]\in\gen^{-2} = 0$. 
Recall the concept of a {\em  crossed module
of Lie algebras}  
 \cite{kassel-loday}. By definition, such a crossed module is 
 a diagram as in \eqref{eq:crossed-Lie}, except that  both $\gen^0$ and $\gen^{-1}$
are assumed to be Lie algebras, $d$ is a homomorphism of Lie algebras and,
in addition, $\gen^0$ is acting on $\gen^{-1}$ by derivations, via a homomorphism
$\alpha: \gen^0\to\Der(\gen^{-1})$ which is required to satisfy
\be
\begin{split}
 [x, y]  = \alpha(d(x))(y), \quad [z, d(x)] = 
d(\alpha(z)(x)), \cr
 z\in \gen^0, x, y\in\gen^{-1}. 
 \end{split}
 \ee
Given a dg-Lie algebra $\gen^\bullet$ situated in degrees $0, -1$, and any
two elements $x, y\in  \gen^{-1}$, we  define their {\em derived bracket} by
\be
[x, y]_{-1}: = [dx, y] = [x,dy]\,\,\in\,\,\gen^{-1},  
\ee
where the second equality follows from $[x,y]=0$. 
We leave to the reader the verification of the following fact. 

 \begin{prop}\label{prop:crossed-Lie-dg}
  The derived bracket  $[-,-]_{-1}$ makes
$\gen^{-1}$ into a Lie algebra (in the usual, ungraded sense),
 and $d$ becomes a homomorphism of Lie algebras. 
Further, the rule 
$$\alpha(z) (x) = [z, x], \quad z\in\gen^0, \, x\in\gen^{-1},$$
defines an action of $\gen^0$ on $(\gen^{-1}, [-,-]_{-1})$ by derivations,
and this action satisfies the axioms of a  crossed module of Lie algebras.
This correspondence establishes an equivalence of categories between 
crossed modules  of Lie algebras, on one hand, and dg-Lie algebras situated
 in degrees
$0, -1$,  on the other hand. \qed

\end{prop}

\end{ex}

Similarly, semiabelian dg-Lie algebras of any length can be seen as Lie algebraic
analogs of crossed complexes of Brown and Higgins \cite{brown-higgins}, see \S 
\ref{subsec:crossed-complexes}
below.

\vskip .2cm

 \begin{ex} \label{ex:semiabelian-filtered}
 Let $\gen^\bullet$ be a dg-Lie algebra
equipped with an increasing  filtration
$$\gen^\bullet_0\subset\gen^\bullet_1 \subset \cdots $$
by dg-Lie subalgebras. Set
$\hen^{-i} = H^{-i}(\gen^\bullet_i/\gen^\bullet_{i-1}).$
These spaces are parts of the term $E_1$ of the spectral sequence of the
filtered complex $\gen^\bullet$. Let us define $d: \hen^{-i}\to\hen^{-i+1}$ to be the 
given by the differential $d_1$ in this spectral sequence.
Then $(\hen^\bullet, d)$ is naturally a semiabelian dg-Lie algebra,
with respect to the induced Lie algebra structure on $\hen^0=H^0(\gen^\bullet_0)$
and its natural action on each $\hen^i$.
Indeed, both conditions of Proposition \ref{prop:semiabelian-conditions}
follow easily from the Leibniz rule. For instance, let us show (b). 

Let $x,y\in \hen^{-1} = H^{-1}(\gen^\bullet_1/\gen^\bullet_0)$, and
choose their representatives $\widetilde x, \widetilde y\in\gen_1^{-1}$.
Then $d\widetilde x, d\widetilde y\in\gen_0^0$, and their classes
in $H^0(\gen^\bullet_0)=\hen^0$ are $dx$ and $dy$. Now, the element
$[dx,y]-[x,dy]\in \hen^{-1}$ is represented by the coboundary
 $[d\widetilde x,\widetilde y]-[\widetilde x, d\widetilde y]= d[\widetilde x,\widetilde y] \in\gen_1^{-1}$
 and therefore is equal to 0.

\end{ex}

Semiabelian dg-Lie algebras share some of the properties of abelian dg-Lie algebras,
i.e., cochain complexes. In particular, if $\gen^\bullet$ is semiabelian,
then the naive truncation
\be
\gen^{\geq -m} \quad =\quad \biggl\{ \gen^{-m} \to \gen^{-m+1} \to \cdots \to \gen^0\biggr\},
 \ee
is again a (semiabelian) dg-Lie algebra. 

\vskip .1cm

Let now $\gen^\bullet$ be any dg-Lie algebra situated in degrees $\leq 0$. We
call the  {\em semiabelianization} of $\gen^\bullet$ and denote by
 $\Gg^\bullet_{\sab}$ the maximal semiabelian quotient of $\Gg^\bullet$:
\be\label{eq:sab}
\Gg^\bullet_{\sab} ={ \Gg^\bullet \over  [\Gg^{\leq -1}, \Gg^{\leq -1}] + 
d([\Gg^{\leq -1}, \Gg^{\leq -1}])}. 
\ee
This construction can be reduced to that of Example \ref{ex:semiabelian-filtered}.
Indeed, let us define the $i$-{\em skeleton} $\gen^\bullet_i\subset\gen^\bullet$
to be the graded Lie subalgebra generated by the subspace $\gen^{\geq -i}$.
This subalgebra is closed under $d$, so we get a filtration of $\gen ^\bullet$
by dg-Lie subalgebras, which we call the {\em skeleton filtration}.

\begin{prop} The dg-Lie algebra $\gen^\bullet_\sab$ is isomorphic to
the semiabelian Lie dg-algebra $\hen^\bullet = (H^{-i}(\gen^\bullet_i/\gen^\bullet_{i-1}))$
associated as in Example \ref{ex:semiabelian-filtered} to the skeleton filtration.
\end{prop}

\noindent {\sl Proof:} By definition, $\gen^{-i}_\sab$ is the quotient of $\gen^{-i}$ by the subspace
$$\sum_{p=1}^{i-1}\, [\gen^{-p}, \gen^{-i+p}]\,\,\,+ \,\,\, \sum_{q=1}^i\,d \, [\gen^{-q}, \gen^{-i-1-q}]. 
$$
On the other hand, $(\gen_i^\bullet/\gen_{i-1}^\bullet)^{-i+1}=0$, while
$$\begin{gathered}
(\gen_i^\bullet/\gen_{i-1}^\bullet)^{-i}\,\,\,=\,\,\, \gen^{-i}\biggl/  \sum_{p=1}^{i-1}\,\,  [\gen^{-p}, \gen^{-i+p}], \cr
(\gen_i^\bullet/\gen_{i-1}^\bullet)^{-i-1}\,\,\, =\,\,\, \sum_{q=1}^i \,\,  [\gen^{-q}, \gen^{-i-1-q}],
\end{gathered}
$$
so taking the cohomology in degree $(-i)$ gives the same answer. We leave the remaining details
to the reader. \qed

\vskip .2cm

We also denote by $\Gg^\bullet_{\CM}$ the
 maximal crossed module quotient of $\Gg^\bullet$, i.e.,
\be\label{eq:cm-quotient}
\Gg^\bullet_{\CM} \quad = \quad \biggl\{{ \Gg^{-1}\over d[\Gg^{-1}, \Gg^{-1}]} \longrightarrow 
\Gg^0\biggr\} \quad = \quad (\Gg^\bullet_{\sab})^{\geq -1}. 
\ee

\subsection{Reminder on differential forms and Schur functors.}
Let   $V$ be a finite-dimensional
$k$-vector space. Consider $V$ as an affine algebraic variety over $k$. 
 The space  $\Omega^p(V)$  of regular $p$-forms on $V$ has the decomposition
\be
\Omega^p(V) = \bigoplus_{d=0}^\infty \Lambda^p(V^*)\otimes S^d(V^*).
\ee
The de Rham differential $d: \Omega^p(V)\to\Omega^{p+1}(V)$ makes $\Omega^\bullet (V)$
into a complex exact everywhere except the $0$th term. 
Let us denote by
\be
\Gamma_p(V) = \bigoplus_{d=0}^\infty \Lambda^p(V)\otimes S^d(V)
\ee
the restricted dual of $\Omega^p(V)$. Geometrically, $\Gamma_p(V)$ can be
seen as the space of de Rham currents on $V$ supported at 0. 
The de Rham differential $d$ on $\Omega^\bullet(V)$
induces differentials $\partial: \Gamma_p(V)\to\Gamma_{p-1}(V)$ by dualization. The dual of
$$\Omega^{p, \cl}(V) = \Ker\bigl\{ d: \Omega^p(V)\to\Omega^{p+1}(V)\bigr\}$$
is then
\be
\Coker\bigl\{ \partial: \Gamma_{p+1}(V)\to\Gamma_p(V)\bigr\} \,\, = \,\,
\Ker \bigl\{ \partial: \Gamma_{p-1}(V)\to\Gamma_{p-2}(V)\bigr\},
\ee
which we denote by $\Gamma_{p-1}^{\cl}(V)$.

\vskip .2cm

Let  $\Vect_k^\fd$ be the category of finite-dimensional
$k$-vector spaces. For any sequence of integers $\alpha = (\alpha_1 \geq ...\geq \alpha_n\geq 0)$
(with arbitrary $n$) we have the Schur functor $\Sigma^\alpha: \Vect^\fd_k\to\Vect^\fd_k$, see
\cite{macdonald}. For $V=k^n$ the space $\Sigma^\alpha(V)$ is the space of irreducible representation of 
the algebraic group $GL_n$ (over $k$)  with highest weight $\alpha$. For example,
$$S^d(V) = \Sigma^{d,0,...,0}(V), \quad \Lambda^d(V) = \Sigma^{1, ..., 1, 0,...,0}(V),$$
where in the last equality there are $d$ occurrences of 1.

The space of closed $p$-forms has the following well known decomposition as a $GL(V)$-module:
\be
\Omega^{p, \cl}(V) \,\,=\,\, d(\Omega^{p-1}(V))\,\, =\,\, \bigoplus_{d=1}^\infty \Sigma^{d, 1, ..., 1}(V^*)
\ee
(with $(p-1)$ occurrences of 1 in the RHS). See, e.g., \cite{GKZ}, Proposition 14.2.2.

\subsection {The semiabelianization of $\Gf^\bullet(V)$.}\label{subsec:semiab-f}
As before, let   $V=k^n$ with basis
$Z_1, ..., Z_n$. We denote by $\Gf^\bullet_{\sab}(V)$ the
semiabelianization of $\fen^\bullet(V)$. 
 Explicitly, it is generated by the elements $Z_I$ as in \S \ref{subsec:resolution}
with differential given by \eqref{eq:differential-f},
but the $Z_I$ are subject to the relations expressing \eqref{eq:sab}. The following gives a first sketch
of the structure of $\Gf^\bullet_{\sab}(V)$.

 \begin{prop}\label{prop:in-f-sab}
 In $\Gf^\bullet_{\sab}(V)$ we have:\hfill\break
(a) Any Lie monomial containing at least two of the $Z_I$, $|I|\geq 2$, vanishes.
\hfill\break
(b) The Lie monomials $[Z_{i_1}, [Z_{i_2}, ..., [Z_{i_p}, Z_J]...]$, $p\geq 0$,
$|J|=m+1$, span $\Gf^{-m}_{\sab}(V)$.\hfill\break
(c) If $m\geq 2$, then the Lie monomials in (b) are symmetric in $i_1, ..., i_p$. 
\end{prop}

\noindent {\sl Proof:} Part (a) is obvious, while (b) follows from (a) and the 
Jacobi identity. Finally, to see (c), let $M= [Z_{i_{a+1}}, [Z_{i_{a+2}}, ..., Z_{i_p}, Z_J]...]$
for some $a$. Then the relations (1.3.9) give
$$0= d([Z_{i_{a-1}, i_a}, M]) = [[Z_{i_{a-1}}, Z_{i_a}], M].$$
Indeed, the term $[Z_{i_{a-1}, i_a}, dM]$ vanishes since both arguments of the bracket
lie in $\Gf^{\leq -1}_{\sab}(V)$. This equality implies
$$[Z_{i_{a-1}}, [Z_{i_a}, M]] = [Z_{i_a}, [Z_{i_{a-1}}, M]],$$
whence the symmetry. \qed

Our main result about the structure of $\Gf^\bullet_{\sab}$ is as follows.

\begin{thm}\label{thm:semiabelianization}
 (a) The dg-Lie algebra $\Gf^\bullet_{\sab}(V)$
is a resolution of $V$, i.e., the projection $\Gf^\bullet(V)\to \Gf^\bullet_{\sab}(V)$
is a quasiisomorphism.\hfill\break
(b) For any $m\geq 1$ the space
$$H^{-m}(\Gf^{\geq -m}_{\sab}(V)) \quad = \quad \Ker\,\,\bigl\{ d: \Gf^{-m}_{\sab}(V)
\longrightarrow \Gf^{-m+1}_{\sab}(V)\bigr\}$$
is identified, as a $GL(V)$-module, with $\Gamma^{\cl}_{m+1}(V)$.
\end{thm}

This theorem will be proved by giving an explicit description of $\Gf^\bullet_{\sab}(V)$.
Let $\widetilde{\Gf}^\bullet(V) \subset \Gf^\bullet(V)$ be the dg-Lie subalgebra given by
\be
\widetilde{\Gf}^0(V) = [\FL(V), \FL(V)] \subset \FL(V) = \Gf^0(V); \quad
\widetilde{\Gf}^i(V) = \Gf^i(V), \, i\neq 0.
\ee
Thus $\widetilde{\Gf}^\bullet(V)$ is acyclic by Proposition \ref {prop:f} 
  (b). We start with describing the
abelianization of this subalgebra. 

 \begin{thm} \label{thm: abelianization}
 The abelianization $\widetilde{\Gf}^\bullet_{\ab}(V)$
is identified, as a cochain complex of $GL(V)$-modules, with the
``co-de Rham complex''
$$\widetilde{\Gamma}^\bullet(V) = \biggl\{ \Gamma_n(V)\buildrel \partial\over\longrightarrow
 \cdots \buildrel\partial\over\longrightarrow \Gamma_2(V) \buildrel \partial\over
\longrightarrow \Gamma_1^{\cl}(V)\biggr\}, $$
with indexing given by
$$\widetilde{\Gamma}^0(V) = \Gamma_1^{\cl}(V), \quad \widetilde{\Gamma}^{-i}(V) = \Gamma_{i+1}
(V), \, i>0.$$
\end{thm}

\begin{rem} This generalizes the result of Reutenauer 
(\cite{reutenauer}, \S 8.6.12),
which says that
\be\label{eq;reutenauer}
[\FL(V), \FL(V)]_{\ab} \,\, = \,\, \bigoplus_{d=1}^\infty \Sigma^{d,1}(V) \,\,= \,\,
\Gamma^{\cl}_1(V) 
\ee
as a $GL(V)$-module.  One can also formulate \eqref{eq;reutenauer}  in terms of the dual spaces
as follows:
\be\label{eq:reutenauer2}
H^1_{\Lie}([\FL(V), \FL(V)], \RR)\,\,=\,\, \widehat\Omega^{2,\cl}_{V,0},
\ee
where the right hand side stands for the space of formal germs of closed 2-forms on $V$ at $0$.
 \end{rem}

The proof of Theorem \ref{thm: abelianization}, given in the next subsection, provides an explicit homomorphism of
dg-Lie algebras (with abelian target)
\be\label{eq:rho}
\rho_\bullet: \widetilde{\Gf}^\bullet(V)\to \widetilde{\Gamma}^\bullet(V),
\ee
inducing an isomorphism on $\widetilde{\Gf}^\bullet_{\ab}(V)$. 

\vskip .2cm

For simplicity, let us drop $V$ from the notation, writing $\Gamma_m$ for $\Gamma_m(V)$,
$\Gf^\bullet$ for $\Gf^\bullet(V)$ etc. Since abelianization factors through semiabelianization,
the homomorphism $\rho_\bullet$ gives rise to a commutative diagram:
\be\label{eq:matrix}
\begin{matrix} \cdots&\buildrel d\over \longrightarrow & \Gf^{-2}_{\sab} &
\buildrel d \over  \longrightarrow & \Gf^{-1}_{\sab}&\buildrel d\over\longrightarrow & [\Gf^0, \Gf^0]\cr
&&&&&&\cr
&&\big\downarrow \overline{\rho}_{-2}&&\big\downarrow \overline{\rho}_{-1}&&\big\downarrow \rho_0\cr
&&&&&&\cr
\cdots& \buildrel\partial\over\longrightarrow& \Gamma_3&\buildrel\partial\over\longrightarrow&\Gamma_2&
\buildrel\partial\over\longrightarrow&\Gamma^{\cl}_1
\end{matrix},
\ee
where $\overline{\rho}_{-m}$, $m\geq 1$, is induced by $\rho_{-m}$. Theorem 
\ref{thm:semiabelianization} will be a consequence
of the following fact.

\begin{thm}\label{thm:cartesian} The rightmost square in \eqref{eq:matrix}
is Cartesian. For $m\geq 2$ the map $\overline{\rho}_{-m}$ is an isomorphism, so
$\Gf^{-m}_{sab}= \Gf^{-m}_{\ab}$. 
\end{thm}

\subsection {Proofs of Theorems \ref{thm: abelianization} and \ref{thm:cartesian}.} We start with constructing the
homomorphism $\rho_\bullet$ from \eqref{eq:rho}.
The map $\rho_0$ is constructed as follows, cf \cite{reutenauer}, \S 5.3.
  Recall that $\Gamma^{\cl}_1(V)$
is the restricted dual of $\Omega^{2, \cl}(V)$. Let $t_1, ..., t_n$ be the linear coordinates
in $V$ corresponding to the basis $Z_1, ..., Z_n$. Then a (closed) 2-form $\omega\in\Omega^2(V)$ is written in these
coordinates as
$$\omega \,\,= \,\,\sum_{i,j} \, \omega_{ij} dt_i dt_j, \quad \omega_{ij}\in k[t_1, ..., t_n].$$
Writing 
$\partial_i = \partial/\partial t_i$,
we define, for any $p\geq 2$:
\be\label{eq:rho-0}
\rho_0\biggl( [Z_{i_1}, [Z_{i_2}, ..., [Z_{i_{p-1}}, Z_{i_p}]...]\biggr)(\omega) \quad=\quad
\bigl(\partial_{i_1} ... \partial_{i_{p-2}} \omega_{i_{p-1}, i_p}\bigr)(0).
\ee
It follows from Reutenauer's theorem (cf. also Corollary 4.4.5 of 
\cite{kapranov-Lie-algebroids} for the case of connections on line bundles)
that $\rho_0$ induces an identification as claimed in Theorem \ref{thm: abelianization}
for the $0$th level.

\vskip .2cm

Further, let us define $\rho_{-m}, m\geq 1$. We view $\Gamma_{m+1}(V)$ as the restricted dual of
$\Omega^{m+1}(V)$, and write any form $\omega\in\Omega^{m+1}(V)$ as $\sum_{ |I|=m+1}  \omega_I dt_I$
similarly to the above. 
We then define:

\vskip .2cm

\be\label{eq-rho-m-a}
 \begin{split}
\rho_{-m}(M)=0,\text{ if } M \text{ is a Lie monomial containing} \cr \text{at least two generators }
 Z_J \text{ with } |J|>1;
 \end{split}
\ee
\be\label{eq-rho-m-b}
\rho_{-m}\biggl([Z_{i_1}, [Z_{i_2}, ..., [Z_{i_p}, Z_I]...]\biggr)(\omega) = 
\bigl(\partial_{i_1} ... \partial_{i_p} \omega_I\bigr)(0). \quad |I|=m+1.
\ee

 \begin{lem}The map $\rho_\bullet$ of graded vector spaces, defined above,
is a homomorphism of dg-Lie algebras. 

\end{lem}

\noindent {\sl Proof:} To show that $\rho_\bullet$ is a morphism of dg-Lie algebras (with
abelian target), it is enough to show that it vanishes on commutators.  For any commutator
involving at least two of the $Z_J, |J|>1$, it follows from \eqref{eq-rho-m-a}. For a commutator
involving none or one of the $Z_J, |J|>1$, it follows from the Jacobi identity and the
fact that the partial derivatives in \eqref{eq:rho-0} or \eqref{eq-rho-m-b} commute with each other.
For example,
$$\rho_{-m} \bigl([[Z_a, Z_b], Z_I]\bigr)(\omega) \quad  =\quad  \rho_{-m} \bigl([Z_a, [Z_b, Z_I]]\bigr)(\omega) - 
\rho_{-m} \bigl( [Z_b, [Z_a, Z_I]]\bigr)(\omega)=$$
$$= (\partial_a\partial_b \omega_I)(0) - (\partial_b\partial_a \omega_I)(0)=0.$$
Next, let us show that $\rho_\bullet$ commutes with the differentials. We denote by
$d_{\Gf}$ the differential in $\widetilde{\Gf}^\bullet(V)$, and by $d_{\DR}$ the de Rham
differential on forms (dual to $\partial$). Then on the generators, we have, by \eqref{eq:differential-f}:
$$d_{\Gf}(Z_{i_0...i_m}) \quad\equiv \quad \sum_{\nu=0}^m (-1)^\nu [Z_{i_\nu}, Z_{i_0, ..., \widehat{i_\nu}, ..., 
i_m}]$$
modulo terns annihulated by $\rho_{-m+1}$. Therefore
\[
\begin{split}(\rho_{-m+1} d_{\Gf}(Z_{i_0, ..., i_m}))(\omega)   = 
\sum_{\nu=0}^m (-1)^\nu (\partial_{i_\nu} \omega_{i_0, ..., \widehat{i_\nu}, ..., 
i_m})(0) =\cr
= (d_{\DR}\omega)_{i_0 ... i_m}(0) = (\partial \rho_{-m} Z_{i_0...i_m})(0).
\end{split}
\]
The statement for more complicated brackets, involving one generator $Z_{i_0...i_m}$ and
several degree 0 generators $Z_{j_1}, ..., Z_{j_p}$, follows in a similar way,
by applying the derivatives $\partial_{j_1}, ..., \partial_{j_p}$ to the above equality.
\qed

\vskip .2cm

Notice further that $\rho_\bullet$ is surjective. So to finish the proof of Theorem 
\ref{thm: abelianization},
it remains to establish the following fact.
 
 \begin{prop}\label{prop:lie-mononials-span}
 The Lie monomials
$$[Z_{i_1}, [Z_{i_2}, ..., [Z_{i_{p-1}}, Z_{i_p}]...], \quad i_1\geq i_2 \geq ... \geq i_{p-1}< i_p, 
\quad p\geq 2, m=0;$$
$$[Z_{i_1}, [Z_{i_2}, ..., [Z_{i_q}, Z_J]...], \quad i_1 \geq ... \geq i_q, \quad q\geq 0, |J|=m+1, m\geq 1,$$
span $\widetilde{\Gf}^{-m}(V)$ modulo $[\widetilde{\Gf}^\bullet(V), \widetilde{\Gf}^\bullet(V)]^{-m}$. 
\end{prop}

Indeed, such monomials correspond to the standard bases in $\Sigma^{p,1}(V)$, resp. $S^q(V)\otimes \Lambda^{m+1}(V)$
which are the direct summands of $\Gamma_1^{\cl}(V)$, resp. $\Gamma_m(V)$. 
So in the case $m\geq 1$, the proposition would imply that $\widetilde\fen^{-m}_{\sab}$ is a
$GL(V)$-module, having an equivariant surjective map $\rho^{\ab}_{-m}$ to the space $\bigoplus_q S^q(V)\otimes\Lambda^m(V)$
and having a spanning set mapping bijectively onto the basis of that space. 
This of course means that $\rho^{\ab}_{-m}$ is an isomorphism. Similarly for $m=0$. 

\vskip .2cm

\noindent {\sl Proof of Proposition \ref{prop:lie-mononials-span}:}
 The case $m=0$ follows from the cited result of Reutenauer. The case $m\geq 1$
is obvious. Indeed, let us neglect all the iterated brackets containing at least one $[Z_{i_\mu}, Z_{i_\nu}]$
and  one $Z_J$, as such brackets belong to $[\widetilde{\Gf}^\bullet(V), \widetilde{\Gf}^\bullet(V)]$. Then,
modulo such brackets, each Lie monomial $[Z_{i_1}, [Z_{i_2}, ..., [Z_{i_q}, Z_J]...]$ depends
on $i_1, ..., i_q$ in a symmetric way, similarly to the proof of Proposition \ref{prop:in-f-sab}(c). 
\qed

\vskip .3cm

We now prove Theorem  \ref{thm:cartesian}. First, we prove that the rightmost square in
\eqref{eq:matrix} is Cartesian. For this, consider the crossed module of Lie algebras
\[
\fen^{\geq -1}_\sab\,\,=\,\,\bigl\{ \fen^{-1}_\sab\buildrel d\over\lra \fen^0_\sab\bigr\}, \quad
\fen^0_\sab=\fen^0=\FL(V)
\]
(with $\fen^{-1}_\sab$ equipped with the derived bracket). As $\fen^0$ is free, any Lie subalgebra
of it is again free (the Shirshov-Witt theorem \cite{reutenauer}). In particular, 
 $\Im(d) = [\fen^0, \fen^0]$ is free.  So the surjection $d: \fen^{-1}_\sab\to \Im(d)$ has a section,
 which is a Lie algebra homomorphism
 $s: \Im(d)\to \fen^{-1}_\sab$ s.t. $ds=\Id$. As for any crossed module of Lie algebras,
 $\Ker(d)$ lies in the center of $\fen^{-1}_\sab$. This implies that 
 \be\label{eq:Lie-direct product-whitehead}
 \fen^{-1}_\sab \simeq \Ker(d) \oplus [\fen^0, \fen^0]
 \ee
 as a Lie algebra. 
Now, to prove that the square in question is Cartesian, means to prove that the (surjective) map
$\overline\rho_{-1}$ maps $\Ker\{d: \fen^{-1}_\sab\to [\fen^0, \fen^0]\}$
isomorphically onto $\Ker\{\partial: \Gamma_2\to\Gamma_1^{\cl}\}$. 
But we know that $\overline\rho_{-1}$ is the degree $(-1)$ component of the abelianization map
  for $\{\fen^{-1}_\sab\to [\fen^0,\fen^0]\}$
considered as a dg-Lie algebra. This abelianization map 
consists  in quotienting  by the action of $[\fen^0, \fen^0]$ ( the degree $0$ part) on 
$\fen^{-1}_\sab$ (the degree $(-1)$ part).
Using the decomposition  \eqref{eq:Lie-direct product-whitehead},
we can describe this action directly:
it   is trivial on the summand $\Ker(d)$ and is the adjoint action on the summand
$[\fen^0, \fen^0]$. This implies that
 $\overline\rho_{-1}$ gives an isomorphism
\[
\Ker(d)  \,\oplus {[\fen^0, \fen^0]\over [[\fen^0, \fen^0], \fen^0, \fen^0]]}
\buildrel \sim\over\lra \Gamma_2. 
\]
and therefore identifies $\Ker(d)$ with $\Ker(\partial)$. 

\vskip .2cm

We now prove that $\overline\rho_{-m}$ is an isomorphism for $m\geq 2$. Note that
by definition \eqref{eq-rho-m-b} of $\overline\rho_{-m}$ the Lie monomials in part (b) of Proposition 
 \ref{prop:in-f-sab} are mapped to the standard basis vectors of 
$\Gamma_{m+1}$. By part (c) of the same proposition, these monomials
are symmetric in $i_1, ..., i_p$, if $m\geq 1$. This means that the set of
distinct monomials in the source of $\overline \rho_{-m}$
maps bijectively on the standard basis in the target, i.e., that $\overline \rho_{-m}$
is an isomorphism, as claimed. \qed

\vfill\eject

\section  {The crossed module of formal 2-branes.}

\subsection{Crossed modules of  groups.}\label{sec:crossed-modules-groups}

\paragraph{A. Crossed modules.}
 The following concept is classical,
see \cite{brown-spencer, mikhailov-passi, sharko} for more background. 

 \begin{Defi}\label{def:CM} 
 A {\em crossed module} (of groups) is a homomorphism of groups
$$
G^\bullet \quad = \quad \bigl\{ G^{-1}\buildrel \partial\over\longrightarrow G^0\bigr\}
$$
together with a left action of $G^0$ on $G^{-1}$ by automorphisms, via a homomorphism
$\beta: G^0\to \Aut(G^{-1})$ which is required to satisfy the axioms:
$$\begin{gathered} g h g^{-1} \,\,= \,\,\beta(\partial(g))(h), \quad u \partial(g) u^{-1} 
\,\,=\,\, 
 \partial(\beta(u)(g)), 
\cr
 u\in G^0, \, g, h\in G^{-1}.
\end{gathered}$$
Morphisms of crossed modules are defined in the obvious way.
We denote by $\CcM$ the resulting category of crossed modules.
\end{Defi}

It follows from the axioms that $\Im(\partial)$ is a normal subgroup in $G^0$, so 
we have the group
$H^0(G^\bullet) = \Coker(\partial)$. We also have the group $H^{-1}(G^\bullet) = \Ker(\partial)$
and this group is contained in the center of $G^{-1}$ in virtue of the axioms;
in particular, it is abelian.  The group $H^0(G^\bullet)$
acts on the abelian group $H^{-1}(G^\bullet)$. 
In the situation of Definition \ref{def:CM} we will sometimes say that $G^{-1}$ is
a {\em crossed $G^0$-module}. 

\begin{ex} \label{ex:crossed-modules-top}
 Let $Y$ be a CW-complex, $X$ a subcomplex and $x\in X$ a point. Then
\be\label{eq:CM-pi}
G^\bullet \,\,=\,\, G^\bullet(Y,X,x)\,\,=\,\,\bigl\{
\pi_2(Y,X,x)\buildrel \partial\over\lra \pi_1(X,x)\bigr\}
\ee
is a crossed module. 
\end{ex}

\paragraph{B. Crossed modules and 2-categories.}
It is  well known \cite{brown-spencer} (see also
\cite{loday}, Lemma 2.2 and \cite{brown}, p. 127) 
that to every crossed module $G^\bullet$ there corresponds a 
strict monoidal 
category $\Cat_\otimes(G^\bullet)$. By definition
\be
\Ob \bigl( \Cat_\otimes(G^\bullet)\bigr)\,\,=\,\, G^0, \quad \Mor\bigl( \Cat_\otimes(G^\bullet)\bigr) \,\,=\,\,
G^0 \ltimes G^{-1}
\ee
(semidirect product with respect to the action $\beta$), with the source,  target 
and the composition maps
  given by
\be\label{eq:2-group-CM}
\begin{gathered}
s,t: G^0 \ltimes G^{-1}\lra G^0, \quad s(u,g)  = u\partial(g), \,\,\, t(u,g)=u.\cr
\circ: (G^0 \ltimes G^{-1})\times_{G^0}^{s,t} (G^0 \ltimes G^{-1})\lra G^0 \ltimes G^{-1}, \cr
(u,g)\circ (v,h) = (u, gh), \quad \text{if}\quad  s(u,g)=t(v,h). 
\end{gathered}
\ee 
  Thus or any two objects $u,u'\in G^0$ the set $\Hom(u,u')$ is identified with
  the set of elements $g\in G^{-1}$
such that $\partial(g) = u^{-1}u'$. 

Note that all the sets in \eqref{eq:2-group-CM} are groups and the axioms of a crossed
module imply that all the maps are group homomorphisms. Therefore 
 $\Cat_\otimes(G^\bullet)$ is a categorical object in the category of groups (a $\Cat^1$-group,
in the terminology of \cite{loday}). In particular, it is a monoidal category, with
the monoidal operation $\otimes$ on objects, resp. morphisms,  given by the group
operation in $G^0$, resp.  $G^0 \ltimes G^{-1}$.

\vskip .2cm

We further recall that a (strict, globular, small) {\em 2-category} $C$ consists of
sets $C_0, C_1, C_2$, whose elements are referred to as {\em objects, 1-morphisms}
and {\em 2-morphisms} of $C$, equipped with maps
(called 1- and 0-dimensional {\em source,  target} and {\em unit maps})
 \be
 \xymatrix{C_2\ar@<.5ex>[r]^{ s_1}\ar@<-.5ex>[r]_{t_1}&
C_1\ar@/_1.2pc/[l]_{\1_1} \ar@<.5ex>[r]^{ s_0\hskip 3mm}\ar@<-.5ex>[r]_{t_0\hskip 3mm}&C_0
\ar@/_1.2pc/[l]_{\1_0}
}, 
\quad s_0s_1=s_0t_1, \,\,t_0s_1=t_0t_1,\,\,s_i\1_i = t_i\1_i
 \ee
and two compositions
\be 
*_0: C_2\times^{s_0s_1, t_0t_1}_{C_0}  C_2\lra C_2, \quad *_1: C_2\times^{s_1, t_1}_{C_1}C_2\lra C_2,
\ee
  called the {\em horizontal} and
{\em vertical compositions} of 2-morphisms. These compositions are required to be associative
to satisfy the 2-dimensional associativity condition, as well as the unit conditions
for the $\1_i$. For more details, see \cite{kelly-street}.
See also  Definition \ref{def:n-cat} below for a more general notion of a strict globular $n$-category.

Any strict  monoidal category  $(\Cc, \otimes)$ gives rise to a 2-category
$2\Cat(\Cc)$ with one object, with 1-morphisms being objects of $\Cc$ and 2-morphisms
being morphisms of $\Cc$. 
The horizontal composition $*_0$ of 2-morphisms is given by the monoidal
structure $\otimes$ on $\Cc$. The vertical composition  $*_1$ 
is given by the composition of morphisms in $\Cc$.
Applying this construction to the monoidal category $\Cat_\otimes(G^\bullet)$,
we associate to a crossed module $G^\bullet$ a 2-category $2\Cat(G^\bullet)$ with one object. 
See \S \ref{subsec:strict-n-cat} for a more general construction, due to Brown-Higgins
\cite{brown-higgins}, which associates an $n$-category to a crossed complex. 

The classifying space of the 2-category $2\Cat(G^\bullet)$ will be denoted $B(G^\bullet)$.
It is a connected topological space with
$$\pi_1(B(G^\bullet))\,\,=\,\,H^0(G^\bullet), \quad \pi_2(B(G^\bullet))\,\,=\,\, H^{-1}(G^\bullet).$$

\paragraph{C. Free crossed modules.} 
 Let $G$ be any group and $(g_a)_{a\in A}$ be any family of elements of $G$.
Let $C$ be the group generated by symbols $\sigma(\gamma,a)$ for all $\gamma\in G$ and $a\in A$,
which are subject to the relations
\be\label{eq:free-CM-relations}
\sigma(\gamma,a) \sigma(\delta, b) \sigma(\gamma,a)^{-1} \,\,=\,\,\sigma(\gamma g_a \gamma^{-1}, b).
\ee
We have a homomorphism $\partial: C\to G$, sending $\sigma(\gamma,a)$ to $\gamma g_a \gamma^{-1}$.
We also have an action of $G$ on $C$ which on generators is given by
$$\beta(g)(\sigma(\gamma,a))\,\,=\,\, \sigma(g\gamma, a).$$
With these data, $\bigl\{ C\buildrel \partial\over\lra G\bigr\}$
is a crossed module, known as the {\em free crossed $G$-module} on the set
of generators $(g_a)$.  It can be characterized by a universal property in the category $\CcM$,
see \cite{sharko}, Def. 6.5.

At the group-theoretical level, the group $H^0=\Coker(\partial)$ for this crossed module is
the quotient of $G$ by the normal subgroup generated by the $g_a$, i.e., the
result of imposing additional relations $g_a=1$ in $G$.

At the 2-categorical level, 
$2\Cat\{C\buildrel \partial\over\lra G\}$ is obtained by adding to $G$
(considered as a usual category with one object) new 2-isomorphisms $\sigma(a): g_a\Rightarrow 1$. 
The generator $\sigma(\gamma,a)$ of $C$ corresponds to the conjugate 2-morphism $\gamma *_0\sigma(a)*_0 \gamma^{-1}$. 

At the topological level, this construction corresponds to adding 2-cells to a CW-complex.
More precisely, 
we have the following fact, due to J. H. C. Whitehead \cite{whitehead-note-on-adding}.

\begin{thm}\label{thm:whitehead1}
Let $Y$ be a connected CW-complex, $X$ a subcomplex such that $Y-X$ consists of
2-dimensional cells only, and $x\in X$. Then  $G^\bullet(Y,X,x)$
is isomorphic to the free crossed $\pi_1(X,x)$-module on the set of generators corresponding to the
boundaries of the attached 2-cells. \qed
\end{thm}

The following particular case, also due to Whitehead \cite{whitehead-note-on-adding}, will be important for us.
For a group $G$ we denote $G_{\ab}=G/[G,G]$ the abelianized group.
For a CW-complex $K$ we denote by $Z_i (K, \ZZ) \subset C_i(K, \ZZ)$
the groups of cellular $i$-cycles and $i$-chains of $K$ with integer coefficients. 

\begin{thm}\label{thm:whitehead2}
Assume, in addition, that $X$ is equal to the 1-skeleton of $Y$ (so $Y$ is 2-dimensional),
and let $G^\bullet=G^\bullet(Y,X,x)$. Denote $\widetilde Y$ the universal covering of $Y$ (with the
CW-decomposition lifting that of $Y$) and let $\widetilde x\in\widetilde Y$ be one of 
the preimages of $x$.
Then:

(a) We have $H^0(G^\bullet)=\pi_1(Y, x)$, while 
$$H^{-1}(G^\bullet) =\pi_2(Y, x) =\pi_2(\widetilde Y, \widetilde x)=H_2(\widetilde Y, \ZZ)
=Z_2(\widetilde Y, \ZZ). $$

(b) Moreover, we have a commutative diagram
 $$\xymatrix{G^{-1} \ar[r] \ar[d]_\partial&G^{-1}_{\ab}\ar[r]^\sim\ar[d]_{\partial_{\ab}}&C_2(\widetilde Y, \ZZ)
 \ar[d]^b\cr
 \Im(\partial)\ar[r]&\Im(\partial)_{\ab}\ar[r]^\sim&C_1(\widetilde Y, \ZZ)
  }
 $$
 in which the left square is Cartesian, the two right arrows are isomorphisms,
 and $b$ is the chain differential. 
 \end{thm} 
 
 For convenience of the reader we add a proof that uses modern terminology.
  The very first equality in (a) follows from the long exact sequence of
 relative homotopy groups and the fact that $\pi_2(X)=0$. Subsequent equalities reflect
 respectively, the invariance of $\pi_2$ in coverings, the Hurewitz theorem, and the
 2-dimensionality of $Y$. 
 
 To see (b), note that $G^0=\pi_1(X,x)$ is free and so $F_0=\Im(\partial)$ is free as a subgroup
 in a free group. Lifting a system of free generators of $F_0$ into $G^{-1}$ in an arbitrary way, we realize $G^{-1}$
 as a semidirect product of $F_0$ and $\Ker(\partial)$. Since $\Ker(\partial)$ is contained in
 the center of $G^{-1}$, we have that $G^{-1}\simeq F_0\times\Ker(\partial)$ is a direct product,
 and $\partial$ is identified with the projection to the first factor. Since $\Ker(\partial)$
 is abelian, this implies that the left square in (b) is Cartesian.  Further, $G^0_{\ab}=H_1(X,\ZZ)=
 C_1(Y,\ZZ)$ and $\widetilde Y_{\leq 1}$,  the 1-skeleton of $\widetilde Y$,  is the covering of $X$ corresponding to
 the subgroup $\Im(\partial)\subset G^0$. This implies that $H_1(\widetilde Y_{\leq 1}, \ZZ) =
 C_1(\widetilde Y, \ZZ)$ is identified with $\Im(\partial)_{\ab}$. 
 
 Finally, we define a homomorphism $f: G^{-1}\to C_2(\widetilde Y, \ZZ)$  using the
 identification of $G^\bullet$ as a free crossed $G^0$-module on the set of generators
 corresponding to the boundaries of 2-cells. If $(c_a)_{a\in A}$ are the 2-cells, and
 $g_a\in G^0$ is the boundary path of $c_a$, then the generator $\sigma(\gamma,a)$, $\gamma\in G^0$,
 is sent by $f$ to the cell $(\overline\gamma)^* c_a$, where $\overline\gamma\in\pi_1(Y,x)$
 is the image of $\gamma$, and $(\overline\gamma)^*$ means the image under the
 covering transformation corresponding to $\overline\gamma$. One verifies that $f$ is
 well defined. Since its target is abelian, $f$ gives rise to $f_{\ab}: G^{-1}_{\ab}\to 
  C_2(\widetilde Y, \ZZ)$ making the right square in (b) commutative. Moreover,
  $f_{\ab}$ identifies $\Ker (\partial_{\ab})$ with $\Ker(b)$ by the Hurewitz theorem,
  so we conclude that $f_{\ab}$ is an isomorphism as well. \qed
  
  \begin{ex}[The crossed module of cubical membranes] \label{ex:cubic-2-branes}
   Let $F(n)$ be the free group on $n$ generators $X_1, ..., X_n$. Consider the following
  family of elements of $F(n)$:
  $$g_{ij} \,=\,[X_i, X_j]\,=\, X_iX_jX_i^{-1}X_j^{-1}, \quad 1\leq i<j\leq n.$$
   The free crossed $F(n)$-module on the set of generators $(g_{ij})$
   will be denoted by $\Box^\bullet_\CM(n)$ and called the {\em crossed module of cubical membranes
   in $\RR^n$.}
   Thus $\Box^0_\CM(n)=F(n)$, while $\Box^{-1}_\CM(n)$ is generated by
   the  symbols $\sigma_{ij}(\gamma)=
   \sigma(\gamma, ij)$,
   $\gamma\in F(n)$ subject to the relations  as in \eqref{eq:free-CM-relations}.
   
   \vskip .2cm
   
   At the topological level, let 
   $\RR^n_\Box$ be the CW-decomposition of $\RR^n$ into cubes
   of the standard cubical lattice, and $\TT^n=\RR^n_\Box/\ZZ^n$ be the $n$-torus  with CW-decomposition
   obtained by identifying
   the opposite sides of the unit $n$-cube. Thus $\TT^n$ has $n\choose i$ cells of dimension $i$.
   For $i\geq 0$ denote by $R^n_{\Box, \leq i}$ and $\TT^n_{\leq i}$ the $i$-skeleta of these CW-complexes.
   In particular, $\TT^n_{\leq 0}$ consists of a single point which we denote $0$.  Then
   Theorem \ref{thm:whitehead1} implies that
   $$\Box^\bullet_\CM(n) \,\,=\,\, G^\bullet(\TT^n_{\leq 2}, \TT^n_{\leq 1}, 0).$$
   Since the universal cover of $\TT^n_{\leq 2}$ is $\RR^n_{\Box, \leq 2}$, 
   Theorem \ref{thm:whitehead2} implies  that
   $\Box^{-1}_\CM(n)$ is the semidirect product of $[F(n), F(n)]$ and the group $Z_2(\RR^n_\Box, \ZZ)$
   of cellular 2-cycles in $\RR^n_\Box$.  
   
   \vskip .2cm
   
   At the  level of monoidal categories, $\Cat_\otimes(\Box^\bullet_\CM(n))$ is generated by invertible objects
   $X_1, ..., X_n$ and invertible morphisms $\sigma_{ij}: X_i\otimes X_j\to X_j\otimes X_i$ which are subject
   to no relations other than those implied by the definition of a strict monoidal category. 
   
   \vskip .2cm
   
   Passing to the 2-category $2\Cat(\Box^\bullet_\CM(n))$ with one object, we notice that its 1-morphisms,
  i.e., elements  of $F(n)$, can be seen as a lattice paths in $\RR^n_\Box$ starting from 0, see, e.g., 
     \cite{kapranov-NCG}.  The end point of the path corresponding to $\gamma\in F(n)$
   is the image of $\gamma$ in $\ZZ^n=F(n)_\ab$. A  2-morphism from $\gamma$ to $\gamma'$
  (existing only if the endpoints of $\gamma$ and $\gamma'$ coincide)  can be visualized
  as a  membrane formed out of 2-dimensional squares of $\RR^n_\Box$ and connecting $\gamma$
  and $\gamma'$. This explains the name ``crossed module of cubical membranes". 
  For instance, the element $\sigma_{ij}(\gamma)\in \Box^{-1}_{\CM}(n)$ can be seen as a
  ``lasso" formed by the square in the direction $(i,j)$ attached at the end of the path $\gamma$.
  This lasso is a 2-morphism from $[X_i,X_j]*_0 \gamma$ to $\gamma$. 
  
  The term ``lasso"  as well as the corresponding ``lasso variables"  were
   introduced by L. Gross \cite{gross-lasso} in connection with quantization of the
  Yang-Mills theory.

  \end{ex}

\subsection{Crossed modules of Lie groups and Lie algebras.}
There are three main constructions relating groups and Lie algebras, and we discuss
their effect on crossed modules.

\paragraph{A. From Lie groups to Lie algebras.}
If  $G$ is a Lie group, we denote by $\Lie(G)$ its Lie algebra (over $\RR$).
The following is then immediate.

 \begin{prop} Let $G^\bullet$ be a crossed module of Lie groups,
and $\Gg^i = \Lie(G^i)$. Then $\Gg^\bullet$ is a crossed module of Lie $\RR$-algebras
(Example \ref{ex:Lie-crossed}) and thus, by Proposition \ref{prop:crossed-Lie-dg}, it
 has a structure of a dg-Lie $\RR$-algebra situated in degrees
$0$ and $(-1)$. \qed
\end{prop}

\paragraph{B. Lower central series and the Magnus Lie algebra.}
If $G$ is any group, we have the {\em lower central series} 
defined in terms of group commutators:
\be\label{eq:low-cen-ser-grp}
G= \gamma_1(G)\supset \gamma_2(G)\supset \cdots, \quad \gamma_{r+1}(G)=[G, \gamma_r(G)].
\ee

  Each $\gamma_{r+1}(G)$ is a normal subgroup
in $\gamma_r(G)$ with abelian quotient.
Let $k$ be a field of characteristic 0. The $k$-vector space
\be
L(G) \,\,=\,\,\bigoplus_{r=1}^\infty \,\,(\gamma_r(G)/\gamma_{r+1}(G))\otimes_\ZZ k
\ee
 with bracket induced by the group commutator in $G$,
is a Lie $k$-algebra.
known as the {\em Magnus Lie algebra} of $G$. 

More generally, let $G^\bullet$ be a crossed module of groups.
We equip $G^{-1}$ by the {\em lower $G^0$-central series} 
\be\label{eq:LCS-CM}
G^{-1}=\gamma_1(G^0, G^{-1})\supset \gamma_2(G^0, G^{-1}) \supset ...
\ee
 where
 $\gamma_{r+1}(G^0, G^{-1})$ is the normal subgroup in $G^{-1}$,
normally  generated by elements of the form
\be\label{eq:LCS-CM-details}
\beta(z)(x)\cdot x^{-1}, \quad     z\in G^0, \,\,x\in \gamma_r(G^0, G^{-1}).
\ee
See \cite{mikhailov-passi},  p. 93.
We will say that $G^\bullet$ is {\em nilpotent}, if both series $(\gamma_r(G^0))$ and
$(\gamma_r(G^0, G^{-1}))$ terminate. 

\begin{prop}  (a)   Both
$$\bigl\{
\gamma_r(G^0, G^{-1})\buildrel \partial \over\lra \gamma_r(G^0)\bigr\} \quad {\rm and}\quad
\bigl\{G^{-1}/\gamma_r(G^0, G^{-1}) \buildrel \partial \over\lra G^0/\gamma_r(G^0)\bigr\}$$
inherit the structures of  crossed modules of groups. 
 
(b) Successive
quotients $\gamma_r(G^0, G^{-1})/\gamma_{r+1}(G^0, G^{-1})$ are abelian, and 
the group commutator in $G^{-1}$ makes
$$L(G^0, G^{-1}) \,\,=\,\,\bigoplus_{r=1}^\infty \,\,(\gamma_r(G^0, G^{-1})/
\gamma_{r+1}(G^0, G^{-1}))\otimes_\ZZ k
$$
into a Lie $k$-algebra. The homomorphism $\partial: G^{-1}\to G^0$ gives rise to a homomorphism
of Lie algebras
$$L (\partial): L(G^0, G^{-1}) \lra L(G^0)
$$

(c) The formula
$$\bigl( \gamma\in \gamma_p(G^0), \, \delta\in \gamma_q(G^0, G^{-1})\bigr) \,\,\longmapsto\,\,
\beta(\gamma)(\delta)\cdot \delta^{-1}\,\,\in\,\, \gamma_{p+q}(G^0, G^{-1})$$
defines an action of $L(G^0)$ on $L(G^0, G^{-1})$ by derivations and makes
$$L (G^\bullet)\,\,=\,\,\bigl\{ L(G^0, G^{-1}) \buildrel L(\partial)\over\lra L(G^0)\bigr\}$$
into a crossed module of Lie algebras. 
 
\end{prop}

\noindent {\sl Proof:}  Straightforward, left to the reader. \qed

\paragraph{C. The Malcev theory.} Let $k$ be a field of characteristic 0. 
If $\gen^\bullet$ is any  dg-Lie $k$-algebra, its {\em  lower central series}   is defined
in terms of Lie algebra commutators:
\be\label{eq:low-cen-ser-Lie} 
\gen^\bullet = \gamma_1(\gen^\bullet)\supset\gamma_2(\gen^\bullet)\supset \cdots, \quad 
\gamma_{r+1}(\gen^\bullet) = [\gen^\bullet, \gamma_r(\gen^\bullet)],
\ee
similarly to
\eqref{eq:low-cen-ser-grp}. As usual, we say that $\gen^\bullet$ is {\em nilpotent}, if
$\gamma_r(\gen^\bullet)=0$ for some $r$. 

\vskip .2cm

Let $\gen$ be a finite-dimensional nilpotent Lie $k$-algebra (with trivial dg-structure).
In this case the Malcev theory  produces 
  a unipotent algebraic group $\exp(\gen)$ over $k$.  More explicitly,
  the augmentation ideal $I\subset U(\gen)$ is in this case topologically
  nilpotent: $\bigcap \, I^d=0$, and so $U(\gen)$ is embedded into the $I$-adic
  completion $\widehat U(\gen) = \varprojlim U(\gen)/I^d$,
  which is a topological Hopf algebra.
  The group $\exp(\gen)(k)$ of $k$-points can be identified with the
  group of group-like elements of $\widehat U(\gen)$, i.e., of $g$
  such that $\Delta(g)=g\otimes g$. 
As a set, it  consists precisely of  elements of the form $\exp(x)=\sum_{i=0}^\infty {x^i/i!}$
for $x\in\gen$. Similarly for points with values in an arbitrary commutative
$k$-algebra.  A moprhism $f: \gen\to\gen'$ of finite-dimensional nilpotent
Lie algebras gives rise to a morphism $\exp(f): \exp(\gen)\to\exp(\gen')$
of algebraic groups. 

If $\gen$ is any Lie $k$-algebra, we define its pro-nilpotent completion to be
$\widehat \gen = \varprojlim \ken$, where $\ken$ runs over 
finite-dimensional nilpotent quotients of $\gen$. We have then the pro-algebraic
group $\exp(\widehat \gen) = \varprojlim \exp(\ken)$. 

\vskip .2cm

Let now $\gen^\bullet$ be a crossed module of Lie algebras. We define the
{\em  lower $\gen^0$-central
series} of $\gen^{-1}$:

\be\label{eq:LCS-CM-Lie}
\begin{split}
\gen^{-1} = \gamma_1(\gen^0, \gen^{-1})\supset\gamma_2(\gen^0,\gen^{-1}\supset\cdots\cr
\gamma_{r+1}(\gen^0, \gen^{-1}) \,\,=\,\,\alpha(\gen^0)\bigl(\gamma_r(\gen^0, \gen^{-1})\bigr).
\end{split}
\ee
 
 \begin{prop} We have the equality
 $$
\gamma_r(\gen^\bullet) \,\,=\,\,
 \bigl\{ \gamma_r(\gen^0, \gen^{-1})\buildrel d\over\lra \gamma_r(\gen^0)\bigr\}.
 $$
 That is, the diagram on the right
 is a crossed submodule of Lie algebras in $\gen^\bullet$,
 and the corresponding dg-Lie algebra  coincides with the $r$-term
 of the lower central series of the dg-Lie algebra  corresponding to $\gen^\bullet$.
   \end{prop}
  
  \noindent {\sl Proof:} Let $x\in\gen^0, y\in \gen^{-1}$. The commutator $[x,y]$
  in the dg-Lie algebra corresponding to $\gen^\bullet$, is equal to
  the element $\alpha(x)(y)$. Therefore the second line in \eqref{eq:LCS-CM-Lie}
  can be written, in terms of the dg-Lie algebra structure,  as follows:
  \[
  \gamma_{r+1}(\gen^0, \gen^{-1})\,\,=\,\,[\gen^0, \gamma_r(\gen^0, \gen^{-1})].
  \]
  This implies the comparison. \qed
  
  \vskip .2cm

So we will use the term {\em nilpotent} for a crossed module $\gen^\bullet$
of Lie algebras to signify that $\gamma_r(\gen^\bullet)=0$, i.e., $\gen^\bullet$
is nilpotent as a dg-Lie algebra. 

\begin{prop}\label{prop:engel-CM}
Let $\gen^\bullet$ be a crossed module of Lie algebras. Then:

(a) If $\gen^\bullet$ is nilpotent, then both
 $\gen^0$ and $\gen^{-1}$ are nilpotent as Lie algebras.

(b) Assume that $\gen^\bullet$ is finite-dimensional. Then $\gen^\bullet$ is nilpotent
if and only if $\gen^0$ is a nilpotent Lie algebra. 

\end{prop}

\noindent {\sl Proof:} (a) Follows from the definition of the bracket $[x,y]_{-1}=[x,d(y)]$
on $\gen^{-1}$ in terms of the dg-Lie algebra structure on $\gen^\bullet$.

(b) The ``only if" part is clear. To see the ``if" part, note that nilpotence
of $\gen^0$ and finite-dimensionality of $\gen^{-1}$ imply that 
 $\gamma_r(\gen^0, \gen^{-1})=0$ for some $r$ in virtue of the Engel theorem. \qed

\begin{prop}\label{prop:2.2.11}
 Let $\gen^\bullet$ be a nilpotent crossed module of finite-dimensional Lie algebras.
Then
$$\exp(\gen^\bullet) \,\,=\,\,\bigl\{ \exp(\gen^{-1})\buildrel \exp(d)\over\lra \exp(\gen^0)\bigr\}$$
has a natural structure of a crossed module of algebraic groups over $k$. This crossed
module is nilpotent. The correspondence $\gen^\bullet\mapsto \exp(\gen^\bullet)$ establishes
an equivalence between the categories of nilpotent crossed modules of finite-dimensional
Lie $k$-algebras and  of nilpotent crossed modules of algebraic groups
over $k$. 
\end{prop}

\noindent {\sl Proof:} Follows from the fact that Malcev's construction $\gen\mapsto\exp(\gen)$
is an equivalence of categories between finite-dimensional nilpotent Lie $k$-algebras
and unipotent algebraic groups over $k$. 
\qed

\begin{rem}
Despite its rather straightforward 
construction,  it is not easy to express the crossed module $\exp(\gen^\bullet)$ in terms of
$\gen^\bullet$ as a dg-Lie algebra. This becomes noticeable for example, when $\gen^\bullet$ is itself obtained from
some other dg-Lie algebra as the maximal crossed module quotient.  In particular, we do not know how to
extract  all of $\exp(\gen^\bullet)$ out of the universal enveloping dg-algebra $U(\gen^\bullet)$
using some analog of the group-like property. Note that $\gen^{-1}$, being an abelian ideal
in $\gen^\bullet$, contributes an ideal in $U(\gen^\bullet)$
isomorphic to the exterior algebra
$\Lambda^\bullet(\gen^{-1})$. On the other hand,  Malcev's construction  
in Proposition \ref {prop:2.2.11} makes use of the (ungraded) enveloping algebra
$U(\gen^{-1}, [-,-]_{-1})$ whose size is that of the symmetric algebra $S^\bullet(\gen^{-1})$.
\end{rem}

\subsection{Unparametrized paths and branes.}\label{subsec:unpar}
In this section we  recall the concept of unparametrized paths and branes 
which based on  the technique  of thin homotopies introduced in 
  \cite{baez-schreiber, martins-picken}.  

\vskip .2cm

 \paragraph{ A.  Differentiable spaces. The space of paths.} 
 We denote by $\Mc an$ the category of $C^\infty$-manifolds
 with corners.  The idea of following definition goes back to K.-T. Chen
 \cite{chen1}. We present a slightly modified version,  closer to that given by R.  M. Hain
(\cite{hain}, Def. 4.1). 
 
 \begin{Defi} A {\em differentiable space} $Y$ is a contravariant functor
 $h=h_Y: \Mc an \to\Sc et$ satisfying the following {\em gluing condition}:
 \begin{itemize}
 \item[(G)]
 Let $(M_i)_{i\in I}$ is any covering of a manifold $M\in\Mc an$  by open subsets. 
 Call a system
 $
 (\phi_i)_{i\in I} \,\in\, \prod_{i\in I} h(M_i)
 $
{\em compatible}, if for each $i,j\in I$  the images of $\phi_i$ and $\phi_j$ in
 $h(M_i\cap M_j)$ coincide. 
 Then the natural map
 $
 h(M) \lra  \prod_{i\in I} h(M_i)
 $
 identifies $h(M)$ with the set  of compatible systems.
 \end{itemize}
 Elements of $h(M)$ are called {\em plots} of $Y$ of type $M$. 
  \end{Defi}
  
\vskip .2cm

\noindent For instance, any $C^\infty$-manifold $N\in\Mc an$ is considered as a differentiable space
via the representable functor $M\mapsto C^\infty(M,N)$. 
For a general differentiable space $Y$ 
it is convenient to  think of $h_Y(M)$ as the set of smooth maps $M\to Y$ and 
to denote this set $C^\infty(M, Y)$ as above. 
In particular, the underlying set of $Y$ is recovered as $h_Y(\pt)$. 

\vskip .2cm

For a differentiable space $Y$ we define the space of smooth differential $p$-forms
on $Y$ (in particular, of smooth functions, for $p=0$) by
\be
\Omega^p_Y \,\,\,=\,\,\,\varprojlim_{(\phi: M\to Y)} \Omega^p_M,
\ee
the limit over all plots. In other words, a $p$-form on $Y$ is a compatible system
of $p$-forms on all the plots.

 \begin{ex}
  Let $X$ be a $C^\infty$-manifold
 and  $PX$ be the space of
parametrized  smooth paths $\gamma: [0,1]\to X$.  
The set $PX$ is made into a differentiable space by putting
$C^\infty(M, PX)$ to be the space of  smooth maps 
$\phi: M\times [0,1]\to X$. 
 The tangent space  $T_\gamma PX$  
is understood as the space   of smooth sections of the vector bundle
$\gamma^* T_X$ on $[0,1]$. We will typically denote such sections by $\delta\gamma$. 
A differential $p$-form $\Phi$ on $PX$ gives a function 
$\Phi(\gamma; \delta_1\gamma, \cdots \delta_p\gamma)$ of a point
$\gamma\in PX$ and $p$  elements of $T_\gamma PX$, and $\Phi$
is uniquely determined by this function. 

We define $P_x^yX$ to be the subspace of paths $\gamma$ as above
such that $\gamma(t)\equiv x$ for $t$ sufficiently close to $0$,
and $\gamma(t)\equiv y$ for $t$ sufficiently close to $1$. 
   Thus $T_\gamma P_x^y X$ consists of 
$\delta\gamma$ vanishing  on some neighborhoods of $0$ and $1$. 
\end{ex}

\begin{rem}\label{rem:constancy}
 In requiring that paths from  $P_x^y X$ 
are {\em constant} on some neighborhoods of $0$ and $1$, we follow
\cite{caetano-picken}. This requirement ensures
that the composition of two such paths is again smooth, not
just piecewise smooth.  On the other hand, any geometric
path can be easily parametrized in this way, so this requirement
leads to no loss of generality. 

Sticking to actual smooth maps 
 becomes particularly important
when we pass from paths to membranes which are maps $\Sigma: [0,1]^p\to X$, 
see Definition \ref{def:param-p-branes} below.
The definition of ``piecewise" smoothness for such maps  
is not quite clear (what kind of ``pieces" are to be allowed so that the desired
constructions go through?).  Instead of attempting such
a definition, we impose on $\Sigma$, as in  \cite{baez-schreiber, martins-picken},
 additional conditions  of constancy (in certain directions)   near the boundary,
 so as to ensure the all the necessary compositions of such membranes are
 again smooth.

\end{rem}

\paragraph {B. Thin homotopies. Unparametrized paths.} 
Let $X$ be a $C^\infty$-manifold and $\gamma, \gamma'\in P_x^y X$ be two parametrized
paths with common source $x$ and target $y$. 

\begin{Defi}
A thin homotopy between $\gamma$ and $\gamma'$ is a $C^\infty$-map $\Xi=\Xi(a_1, a_2): [0,1]^2\to X$
such that:
\begin{enumerate}
\item[(1)] $\Xi$ is constant in some neighborhoods of te faces $a_1=0$, $a_1=1$.

\item[(2)] $\Xi$ depends only on $a_1$ in some neighborhoods of the faces $a_2=0$, $a_2=1$.

\item[(3)] $\Xi(a_1,0)=\gamma(a_1)$, $\Xi(a_1, 1)=\gamma'(a_1)$ for any $a_1\in[0,1]$.

\item[(4)] The rank of the differential $d_a\Xi: T_a[0,1]^2 \to T_{\Xi(a)} X$ at any $a\in[0,1]^2$ is $\leq 1$. 

\end{enumerate}
\end{Defi} 

The first condition means that all the intermediate paths $\gamma_{a_2}: a_1\mapsto \Sigma(a_1, a_2)$
lie in $P_x^y X$. The second condition means that $\gamma_{a_2}$ coincides with $\gamma$ for $a_2$
close to 0 and with $\gamma'$ for $a_2$ close to 1. It ensures that being thin homotopic is an equivalence
relation on $P_x^y X$ which we denote $\gamma\sim\gamma'$. The 4th condition (``thinness") 
means that the relation 
$\sim$ includes, in particular, 
reparametrization of paths as well as
 cancellation of a segment and the same segment run in the opposite direction
 immediately after.  
 
 We denote by $\Pi_x^y X = P_x^y X/\sim$ the set of this homotopy classes of paths from $P_x^yX$.
 By the above, we can consider elements of $\Pi_x^y X$ as {\em unmarametrized paths} from $x$ to $y$. 
We define the {\em groupoid of unparametrized paths} $\Pi_{\leq 1} X$ to be the category with
\[
\Ob (\Pi_{\leq 1} X) = X, \quad \Hom_{\Pi_{\leq 1} X}(x,y) = \Pi_x^y X.
\]
 Composition of paths is given by
 the standard concatenation (which gives a smooth path, see Remark \ref{rem:constancy}):
 \be\label{eq:comp-paths}
 (\gamma*\gamma')(a) = \begin{cases} \gamma'(2a), \quad {\rm if}\quad 0\leq a\leq {1\over 2} \cr 
 \gamma(2a-1), \quad {\rm if}\quad {1\over 2}\leq a\leq 1.
 \end{cases}
 \ee
 Because of the thin homotopy relation, this composition
 is associative, and each path is invertible.

 \paragraph{C. Globes and branes.} 
Let $I=[0,1]$ be the unit interval, so $I^p$ is the unit $p$-cube
with coordinates $a_1, ..., a_p \in[0,1]$. The following definition is a reformulation
of one from \cite{brown:new}

\begin{Defi}\label{def:param-globes} Let $X$ be a topological  space. 
 A {\em singular $p$-globe} in $X$ is a continuous
 map $\Sigma: I^p\to X$ which for each $i=1, ..., p$ satisfies
 the following condition:
 \begin{itemize}
 \item[($\operatorname{Glob}_i$)] The restrictions of $\Sigma$  to  the
 faces $\{t_i=0\}$ and $\{t_i=1\}$ of $I^p$  depend only on the coordinates $t_1, ..., t_{i-1}$
 (in particular, are constant, if $i=1$). 
  \end{itemize}
  We denote by $\Glob_p(X)$ the set of singular  $p$-globes in $X$.
  \end{Defi} 
 
 The conditions $(\Glob_i)$ mean that singular $p$-globes factor through the
 {\em universal singular $p$-globe} $\alpha_p: I^p \to \bigcirc^p$. More precisely,  $\bigcirc^p$ is the quotient
 of $I^p$ by the identifications coming from the $(\Glob_i)$. As well known $\bigcirc^p$ can be identified
 with the unit $p$-ball $D^p= \{ x\in \RR^p: \| x\|\leq 1\}$, see \cite{brown:new}, Def. 2.1 for an explicit
 map $\beta_p: I^p\to D^p$ establishing this identification. Further, the cell structure on $I^p$ given by the faces,
 is contracted by $\beta_p$ into the cell structure on $D^p$ given by the interior  open ball and the hemispheres:
 \[
 D^p \,\,=\,\, e^p \cup e_{\pm}^{p-1} \cup e_\pm^{p-2} \cup \cdots \cup e_\pm^0,
 \]
 where 
 \[
 e_\pm^i \,\,=\,\,\bigl\{ x\in D^p: \,\, \|x\|=1, \,\, x_j=0 \text{ for } j+i<p, \,\, \pm x_{p-i}>0.\bigr\}
 \]
The conditions $(\Glob_i)$ implies that for a singular $p$-globe $\Sigma: I^p\to X$ and $  q<p$
the restrictions of $\Sigma$ to the faces $a_q=0$ and $a_q=1$ define singular $p$-globes
which we denote $s_q\Sigma$ and $t_q\Sigma$. These globes are called the {\em $p$-dimensional
source and target} of $\Sigma$. 
Thus we have maps
\be
\label{eq:s_p}
s_q = s_q^p, \, T_q=t_q^p:\,\, \Glob_p(X)\lra \Glob_q(X).
\ee
 
 \begin{ex}
 For the universal globe $\beta_p: I^p\to D^p \simeq \bigcirc^p$ the globe $s_q\beta_p$ takes values in the
 closed hemisphere $e_-^q$ and $t_q\beta_p$ takes values in $e_+^p$. So we will denote these
 hemispheres by $s_q\bigcirc^p$, $t_q\bigcirc^p$. 
 \end{ex}
  
  We now consider the  smooth setting and introduce a slight modification of the above concepts, adapted along the lines of
  Remark \ref{rem:constancy}. 
  
\begin{Defi}\label{def:param-p-branes} Let $X$ be a differentiable space. 
 A {\em parametrized $p$-brane} in $X$ is a  smooth
 map $\Sigma: I^p\to X$ which for each $i=1, ..., p$ satisfies
 the following condition:
 \begin{itemize}
 \item[($\operatorname{Br}_i$)] There exist neighborhoods of the
 faces $\{t_i=0\}$ and $\{t_i=1\}$ of $I^p$ in which
 $\Sigma$ depends only on the coordinates $t_1, ..., t_{i-1}$
 (in particular, is constant, if $i=1$). 
  \end{itemize}
   \end{Defi} 
  We denote by $\widetilde \Pi_p(X)$ the set of parametrized $p$-branes in $X$.
Thus, for a smooth manifold $X$ we have $\widetilde\Pi_p(X) \subset\Glob_p(X)$.    
 It is clear that the maps $s_q, t_q$ from \eqref{eq:s_p} take $\widetilde\Pi_p(X)$
 to $\widetilde\Pi_q(X)$. 
 Note further that a parametrized $p$-brane in $X$ descends, via $\alpha_p$,  to a {\em smooth}
 map $\bigcirc^p=D^p \to X$. 
    
    \begin{ex}
     We have $\widetilde\Pi_0(X)=X$, while
    $\widetilde\Pi_1(X)\subset\ PX$ is the space of paths $\gamma: I\to X$
     which are constant in some neighborhoods of
$0$ and $1$.  Let $\Sigma\in\widetilde\Pi_2(X)$. Then by
($\operatorname{Br}_2$) we have that $\Sigma$ is constant on
some neighborhoods of the intervals $\{0\}\times I$ and $\{1\}\times I$. 
    Denoting the images of these
 neighborhoods  by $x$ and $y$, we can associate to  $\Sigma$ a
 smooth path   
  \[
 \sigma = \tau\Sigma: [0,1]\to P_x^yX, \quad \sigma(a)(b)=\Sigma(b,a), \,\, a,b\in [0,1],
 \] 
 called the {\em transgression} of $\Sigma$.

 \end{ex}

 \paragraph{D. Thin homotopy of branes. Unparametrized branes.} 
 Let $X$ be a $C^\infty$-manifold and $\Sigma, \Sigma': \bigcirc^p \to X$
 be two parametrized $p$-branes. 
 
 \begin{Defi}\label{def:hom-branes}
 A {\em homotopy} between $\Sigma$ and $\Sigma'$ is a smooth
 map $\Xi: I\times \bigcirc^p\to X$ such that:
 \begin{enumerate}
\item[(1)] For each $b\in I$ the map $\Xi_b= \Xi(b, -): \bigcirc^p\to X$ is a parametrized $p$-brane in $X$.

\item[(2)] $\Xi_b$ is independent of $b$ for $b$  in some neighborhoods of  $0$ and $1$ and equals
in these neighborhoods to $\Sigma$ and $\Sigma'$ respectively.   
 \end{enumerate}
 \end{Defi}
 
 If $q<p$, then by restricting to $I\times s_q\bigcirc^p$, the homotopy $\Xi$ induces  a homotopy $\Xi^{s_q}$ between
 $s_q\Sigma$ and $s_q\Sigma'$. We similarly obtain a homotopy $\Xi^{t_q}$ between $t_q\Sigma$
 and $t_q\Sigma'$.

 \begin{Defi}
 A homotopy $\Xi: I\times
 \bigcirc^p\to X$ between $\Sigma$ and $\Sigma'$ is called {\em thin}, if:
 \begin{enumerate}
 \item[$(T_p)$] The rank of the differential of $\Xi$ at any point of $I\times\bigcirc^p$ is $\leq p$.
  \item[$(T_{< p})$] For any $q<p$, the homotopies $\Xi^{s_q}: I\times s_q\bigcirc^p \to X$, 
  $\Xi^{t_q}: I\times t_q\bigcirc^p \to X$ are such that their differential at each point has rank $\leq q$. 
  \end{enumerate}
 \end{Defi}
 
 \begin{rems}\label{rems:thin}
 (a) In particular, for a thin homotopy $\Xi$ between $\Sigma$ and $\Sigma'$ we have $s_0\Sigma=s_0\Sigma'$
 amd $t_0\Sigma=t_0\Sigma'$, as $\Xi^{s_0}$ and $\Xi^{t_0}$ have differentials of rank $0$ and so are
 constant maps. 
 
 (b) A requirement equivalent to $(T_p)$ would be tto say that the image of $\Xi$ has Hausdorff
 dimension $\leq p$, similarly $(T_{<p})$ can be expressed by saying that the images
 of $\Xi^{s_q}$, $\Xi^{t_q}$ have Hausdorff dimension $\leq q$. 
 \end{rems}
 
 Note, in particular, that thin homotopies include reparametrizations of branes as well
 as ``$p$-dimensional cancellations" (preserving the boundary in a compatible way). 
 We also see that being think homotopic is an equivalence relation on $\Glb_p(X)$,
 which we denote $\sim$.
 Define the {\em set of unparametrized $p$-branes} in $X$ as the quotient
 \be
 \Pi_p(X) = \Glb_p(X)/\sim.
 \ee
 Elements of $\Pi_p(X)$
  can be thought of as ``geometric pieces of $p$-dimensional surfaces in $X$", free
of the choice of parameters.

\paragraph{E. The 2-groupoid of unparametrized 2-branes.} Let $X$ be a $C^\infty$-manifold.
The sets of parametrized $p$-branes in $X$ for $p\leq 2$, organize themselves into a diagram
\be\label{eq:tilde-pi-leq-2}
\xymatrix{\widetilde \Pi_2(X)\ar@<.5ex>[r]^{ s_1}\ar@<-.5ex>[r]_{t_1}&
\widetilde \Pi_1(X) \ar@<.5ex>[r]^{ s_0\hskip 3mm}\ar@<-.5ex>[r]_{t_0\hskip 3mm}&
\widetilde \Pi_0(X)=X,
}  
\ee
which descents to diagram of sets of unparametrized branes 
\be\label{eq:pi-leq-2}
\xymatrix{\Pi_2(X)\ar@<.5ex>[r]^{ s_1}\ar@<-.5ex>[r]_{t_1}&
\Pi_1(X) \ar@<.5ex>[r]^{ s_0\hskip 3mm}\ar@<-.5ex>[r]_{t_0\hskip 3mm}&\Pi_0(X)=X.
} 
\ee
Following  \cite{baez-schreiber, martins-picken}, 
 we make  \eqref{eq:pi-leq-2} into a 2-category, in fact a 2-groupoid $\Pi_{\leq 2}(X)$
with the set of $i$-morphisms being $\Pi_i(X)$, $i=0,1,2$. 
For this, we first introduce the composition \eqref{eq:comp-paths} on $\Glb_1(X)$ and 
 two compositions $*_0$ and $*_1$ on $\widetilde \Pi_2(X)$,
 called the {\em horizontal} and {\em vertical}
composition of 2-branes, defined similarly to \eqref{eq:comp-paths} by
\be
\begin{aligned}
 (\Sigma*_0\Sigma')(a_1, a_2) = \begin{cases} \Sigma'(2a_1,a_2), \quad {\rm if}\quad 0\leq 1\leq {1\over 2} \cr 
 \Sigma(2a_1-1,a_2), \quad \text { if }\quad {1\over 2}\leq a_1\leq 1,
 \end{cases} \text{ when }  s_0 \Sigma = t_0 \Sigma';
 \cr
(\Sigma*_1\Sigma')(a_1,a_2) = \begin{cases} \Sigma'(a_1,2a_2), \quad {\rm if}\quad 0\leq a_2\leq {1\over 2} \cr 
 \Sigma(a_1,2a_2-1), \quad \text { if }\quad {1\over 2}\leq a_2\leq 1,
 \end{cases}\quad \text{ when }  s_1\Sigma=t_1\Sigma'.
 \end{aligned}
\ee
These operations do not satisfy the axioms of a 2-category (in particular, they are not associative)
but they ``descend" to operations on $\Pi_1(X), \Pi_2(X)$ which do. Such descent is completely
clear for the composition of paths in $\Pi_1(X)$ and for the $*_0$-composition of thin homotopy
classes from $\Pi_2(X)$. 

For $*_1$ the descent needs more explanation.  Let us denote by $[\Sigma]$ the thin homotopy class
of a parametrized p-brane $\Sigma$, $p=1,2$. Then we need to define $[\Sigma]*_1 [\Sigma']$
for any parametrized 2-branes $\Sigma$ and $\Sigma'$ such that the parametrized path
$s_1\Sigma$ is thin homotopic (not necessarily equal) to $t_1\Sigma'$. In order to do this, we consider the corresponding
thin homotopy as a parametrized 2-brane $\Xi$ with $s_1\Xi=t_1\Sigma'$ and $t_1\Xi=s_1\Sigma$, and
define
\be\label{eq:filling}
[\Sigma] *_1 [\Sigma'] \,\,=\,\, [ (\Sigma*_1\Xi) *_1\Sigma']. 
\ee

It is shown in \cite{baez-schreiber, martins-picken}
that these composition define a 2-category, in fact a 2-groupoid $\Pi_{\leq 2}(X)$
with the set of $i$-morphisms being $\Pi_i(X)$, $i=0,1,2$. We will
call $\Pi_{\leq 2}(X)$ the {\em 2-groupoid of 2-branes} in $X$.

\subsection{2-dimensonal holonomy.}\label{subsec:2-dim-hol}
Here we give a summary of  the main points of the theory of connections
with values in gerbes (crossed modules) and their 2-dimensional holonomy
as developed by Breen and Messing  \cite{breen-messing}
and Baez-Schreiber \cite{baez-schreiber}. 
Our exposition of the holonomy follows the approach  
 \cite{baez-schreiber}   based on the 
``covariant" generalization of Chen's theory of iterated integrals,
as developed in 
\cite{hofman} \cite{baez-schreiber}. Additional details of this theory can be
found in the Appendix.

\paragraph{A. The Schlessinger formula.}
 Let $X$ be a $C^\infty$-manifold, let $G$ be a Lie group with Lie algebra $\gen$,  let
 $Q$ be a principal $G$-bundle on $X$, and $\nabla$ a connection in $Q$.
 We denote by $\Ad(Q)$ the bundle of Lie algebras on $X$ associated to $Q$
 via the adjoint representation of $G$. Each fiber $\Ad(Q)_x$ is a Lie algebra
 is isomorphic (non-canonically) to $\gen$. 
 In a trivialization of $Q$ the connection $\nabla$ is given as $d_{\DR}-A$ for
 a $\gen$-valued 1-form
  $A\in\Omega^1_X\otimes \gen$.
   Denote by 
$F_\nabla\in\Omega^2_X\otimes \Ad(Q)$ the curvature of $\nabla$.
Locally, if $\nabla$ corresponds to $A$ as above, then
$$
F_\nabla \,\,= \,\,F_A \,\,= \,\, dA - {1\over 2}[A,A]\,\,
\in \,\, \Omega^2_X\otimes \gen
$$
For any parametrized
 smooth path $\gamma: [a,b]\to X$ (defined on any interval $[a,b]\subset\RR)$
we have  the holonomy
$M_\nabla(\gamma) : Q_{\gamma(a)} \to Q_{\gamma(b)}$.
In particular, for any $x,y\in X$ the holonomy gives
 a smooth function
  $$
  M_\nabla: P_x^yX\lra \Hom_G(Q_x,Q_y).
  $$ 
  We denote the 
   logarithmic (group-theoretic) differential of $M_\nabla$ by
\[
  B= M_\nabla^{-1} d_{\DR}(M_\nabla) \,\,\in \,\, \Omega^1_{P_x^yX}\otimes \Ad(Q)_x.
  \]
Then, for any $\delta\gamma$ vanishing near the ends
we have  the {\em Schlessinger formula}:
\be\label{eq:form-B}
B(\gamma; \delta\gamma) \,\,=\,\,\int_0^1 M_\nabla(\gamma_{\leq t})^{-1} \cdot F_\nabla\bigl(
\gamma(t); \gdot (t), \delta\gamma(t)\bigr)
\cdot M_\nabla(\gamma_{\leq t}) dt,
\ee
where $\gamma_{\leq t}: [0,t]\to X$ is the restriction of $\gamma$ to $[0,t]$.
An equivalent (integrated) formulation of \eqref{eq:form-B} is that for any smooth path  $\sigma: [0,1]\to P_x^yX$ 
beginning at  $\sigma(0)=\gamma$ and ending at $\sigma(1)=\zeta$, the ratio 
$M_\nabla(\zeta)^{-1} M_\nabla(\gamma)$
can itself  be represented as the holonomy of the connection ${d\over ds}-\sigma^*(B)$ on $[0,1]$.
This integrated statement was proved by   L. Schlessinger in his 
1928 paper  \cite{schlessinger}, formula (22), the holonomy of ${d\over ds}-\sigma^*(B)$
being precisely his ``gemischtes Doppelintegral".  

The Schlessinger formula  can be used to give a clear proof of the following well known fact.

 \begin{prop} Let  $G$ be a Lie group, $Q$ a principal $G$-bundle on $X$
 and $\nabla$ a connection in $Q$. If $\gamma_0$ and $\gamma_1$ are rank-1
 homotopic, then the corresponding holonomies are equal:
 $M_\nabla(\gamma_0)=M_\nabla(\gamma_1)$. 
 The holonomies of   $\nabla$  
 define then a functor from $\Pi_{\leq 1}(X)$ to the category of $G$-torsors. 
  \end{prop}
  
  \noindent {\sl Proof:}  The pullback of the curvature 2-form along a rank-1
  homotopy $I^2 \to X$ is zero as a 2-form on $I^2$,
   so the Schlessinger formula implies 
  the first statement. The second statement is obvious. \qed

\paragraph{B. Covariant transgression.}
Let  $\beta: G\to \Aut(V)$ be a smooth representation of $G$ in a finite-dimensional
$\RR$-vector space $V$, and $V(Q)$ be the vector bundle on $X$ associated to $Q$.
We 
define the {\em $\nabla$-transgression} map
\be\label{eq:covariant-transgression}
\gathered
\oint_\nabla: \Omega^{m+1}_X\otimes V(Q) \lra\Omega^{m}_{P_x^y X}\otimes V(Q)_x,\quad 
\left(\oint_\nabla(\Phi)\right)(\gamma;  \delta_1\gamma, \cdots, \delta_{m}\gamma)\,\, :=\cr
=\, \,\int_0^1 \beta \bigl(M_\nabla(\gamma_{\leq t})\bigr)^{-1}
\bigl(\Phi(\gamma(t); \gdot(t), \delta_1\gamma(t), \cdots , \delta_{m}\gamma(t))\bigr) 
dt.
\endgathered
\ee
When $Q$ is the trivial bundle, and $\nabla = d-A$ is given by a 1-form $A\in\Omega^1_X\otimes\gen$,
we will write $\oint_A$ for $\oint_\nabla$. 

  \begin{ex}
  (a)  The usual (non-covariant)  transgression of (scalar) differential forms
\be\label{eq:usual-transgression}
\oint: \Omega^{m+1}_X \lra\Omega^m_{P_x^y X}
\ee
 is obtained as a particular case,  when  
 $\nabla$ is the trivial connection in the trivial bundle
  and $V$ is the trivial 1-dimensional representation.
In this case the action by the holonomy drops out.

(b) As another example, note that  Schlessinger's  form $B$ above corresponds to the case when
 $V=\gen$ is the adjoint
representation: $B=\oint_\nabla(F_\nabla)$. 
 \end{ex}
 
 When it is important to emphasize the dependence of the covariantly transgressed form on the points $x,y\in X$,
 we will use the notation
 \[
 \biggl(\oint_\nabla (\Phi)\biggr)_x^y \,\, \in \,\, \Omega^m_{P_x^y X}\otimes V(Q)_x. 
 \]
 We will not discuss the relations among such forms for different $x,y\in X$. More precisely,
 let $x,y,z\in X$ be three points. We then have the composition map
 \[
 \on{com}: P_y^z X \, \times \, P_x^yX \lra P_x^z X, 
 \quad \on{com}(\gamma, \gamma')(t) = 
 \begin{cases}
 \gamma'(2t),& \text{ if } t\in [0, {1\over 2}],
 \\
 \gamma(2t-1), & \text {if } t\in [{1\over 2}, 1].
 \end{cases}
 \]
 For $\gamma\in P_y^ X$ we have the {\em left translation map}
 \[
 l_\gamma: P_x^y X \lra p_x^z X, \quad r_\gamma(\gamma') = \on{com}(\gamma, \gamma').
 \]
For $\gamma'\in P_x^yX$ we have the {\em right translation map}
\[
r_{\gamma'}: P_y^x Z \lra P_x^z X, \quad l_{\gamma'}(\gamma) = \on{com}(\gamma, \gamma'). 
\] 

The following ``translation invariance" property is a direct consequence of the definitions.
\begin{prop}\label{prop:trans-trans}
(a) For each $\gamma\in P_y^z X$ we have
\[
r_\gamma^*  \biggl(\oint_\nabla (\Phi)\biggr)_x^z \,\,=\,\,  \biggl(\oint_\nabla (\Phi)\biggr)_x^y.
\]
(b)  For each $\gamma'\in P_x^y X$ we have
\[
l_\gamma^*  \biggl(\oint_\nabla (\Phi)\biggr)_x^z \,\,=\,\, M_\nabla(\gamma)^{-1} \biggl(\oint_\nabla (\Phi)\biggr)_y^z
\]
(the result of applying the isomorphism $M_\nabla(\gamma)^{-1}: V(Q)_y\to V(Q)_x$ to an
$V(Q)_y$-valued $m$-form). \qed

\end{prop}

\paragraph{C. Connections with values in crossed modules and on the space of paths.}
 Let now $G^\bullet$ be a crossed module of Lie groups and $\gen^\bullet$ the corresponding
crossed module of Lie algebras.   

\vskip .2cm

We first consider $\gen^\bullet$ as a dg-Lie algebra situated
in degrees $-1,0$. Let $X$ be a $C^\infty$-manifold and 
and let $A^\bullet\in(\Omega^\bullet_X\otimes \gen^\bullet)^1$ be
 a graded connection on $X$ with values in $\gen^\bullet$, i.e., 
a $\gen^\bullet$-valued differential form of total degree 1,
see \S \ref {subsec: connection}. 
Thus $A^\bullet$ has two
components:
$$A^1\in\Omega^1_X\otimes \gen^0, \quad A^2\in\Omega^2_X\otimes \gen^{-1}.$$
 We can consider $A^1$ as a usual connection in the trivial $G^0$-bundle on $X$. 
Let $F^\bullet $ be the curvature of  the graded connection $A^\bullet$,
see \eqref{eq:graded-curvature}. 
 Thus $F^\bullet$ has two components:
 \be F^2 = F_{A^1} - d_{\gen^\bullet}(A^2)\in\Omega^2_X\otimes\gen^0, \,\,\,
 F^3=d_{\DR}(A^2) - [A^1, A^2]\in \Omega^3_X\otimes \gen^{-1}.
  \ee
  
  Now  consider $\gen^\bullet$ as a crossed module of Lie algebras, as it
  was originally defined. From this
point of view, 
  $A^\bullet$ is  a 2-connection in the trivial $G^\bullet$-gerbe \cite{breen-messing} 
 on $X$.  The component $F^2$ is  traditionally called the {\em fake curvature}, while
  $F^3=\nabla_{A^1}(A^2)$  is called   the {\em 3-curvature}
  of  $A^\bullet$, see \cite{breen-messing}. 
  We say that the 2-connection $A^\bullet$
  is {\em semiflat}, if the fake curvature vanishes: $F^2=0$.

  \vskip .2cm
  
  Fixing points $x,y\in X$, we associate to  $A^\bullet$  a usual connection
 in the trivial $G^{-1}$-bundle on the space $P_x^yX$ given by the $\gen^{-1}$-valued
 1-form
 $$\oint_{A^1} (A^2) \,\,\in\,\,\Omega^1_{P_x^y X} \otimes \gen^{-1}.$$ 
   We have the following fundamental fact, see \cite{baez-schreiber}, formula (2.57).
 
 \begin{prop}\label{prop:curv-tramsgression}
 If the 2-connection $A^\bullet$ is semiflat, then the curvature of $\oint_{A^1}(A^2)$
 is equal to $\oint_{A^1}(F^3)$. 
  \end{prop}
  
  For convenience of the reader we give here a proof,  referring to the
  Appendix for necessary background. By defiiniton, the curvature in question is equal to
  \[ 
  d_\DR \oint_{A^1} (A^2) \,\,\, - \,\,\, {1\over 2} \biggl[  \oint_{A^1} (A^2), \,  \oint_{A^1} (A^2)\biggr].
  \]
  Note first that the commutator in the
  second term vanishes. Indeed, by definition, this  is the commutator of two
  integrals over $[0,1]$, which is expressible as an integral over $[0,1]^2$ of the
  pointwise commutators of the two  integrands. This new integrand is
  a function on $[0,1]^2$ which is antisymmetric with respect to the 
  interchange of the two variables, so the integral over $[0,1]^2$ vanishes. 
     
  Second, we recall Example \ref{ex:two-trans} which implies that we can apply 
  Proposition \ref{prop:nabla-covariant-iter-integral}(b) to calculate the exterior derivative of
  the covariant transgression. It gives that the 
    first term in the above formula for the curvature is equal to 
  \be\label{eq:d-of-int-A2}
  \begin{gathered}
   d_\DR \oint_{A^1} (A^2) \,\,=\,\,\oint_{A^1} (\nabla_{A^1}(A^2)) \,\,+\,\,\loint_{A^1}( F_{A^1}, A^2)
   \,\,-\,\,\loint_{A^1}(A^2, F_{A^1})\,\,=\cr
   =\,\, \oint_{A^1} (F^3) \,\,+\,\, \biggl[    \oint_{A^1} (F_{A^1}), \,\, \oint_{A^1} (A^2)
   \biggr],
   \end{gathered}
  \ee
  where the last term is the wedge product of two 1-forms on $P_x^yX$ followed by
  the commutator, i.e., action of $\alpha: \gen^0\otimes\gen^{-1}\to\gen^{-1}$. 
   Now, the condition of semiflatness $F_{A^1}=d_\gen(A^2)$ implies that the last
  wedge commutator above is equal to 
  \[
  \biggl[ d_\gen \oint_{A^1} (A^2), \,\, \oint_{A^1} (A^2) \biggr],
  \] 
  which is zero since for each $y\in\gen^{-1}$ we have $\alpha(d_\gen(y))(y) = d_\gen [y,y]=0$. 
  \qed

\paragraph{F.  2-dimensional holonomy.} Let $G^\bullet=\{G^{-1}\buildrel \partial\over\to G^0\}$
 be a crossed module of Lie
groups, $\gen^\bullet$ the corresponding crossed module of Lie algebras, and
$A^\bullet=(A^1, A^2)$ be a connection on $X$ with coefficients in $\gen^\bullet$,
as in \S D. Let $\Sigma$ be a parametrized 2-brane in $X$.  The {\em 2-dimensional holonomy}
of $A^\bullet$ along $\Sigma$ is defined as
\be\label{eq:def-2-hol}
M_{A^\bullet}(\Sigma) \,\,=\,\, M_{\oint_{A^1}(A^2)} (\sigma) \,\,\in \,\,G^{-1}.
\ee
Here $\sigma: [0,1]\to P_x^yX$ is 
the   path in the space of paths corresponding to $\Sigma$,  and
$\oint_{A^1}(A^2)$ is 
the connection on $P_x^y X$ obtained as
 the covariant transgression of $A^2$ with respect to $A^1$.
 The following is the main result of \cite{baez-schreiber} (specialized for the case of
 2-connections in trivial 2-bundles).
 
 \begin{thm}\label{thm:2-mon}
  Assume that $A^\bullet$ is semiflat. Then:
 
 (a) $M_{A^\bullet}(\Sigma))$ represents a 2-morphism in $2\Cat(G^\bullet)$ from
 $M_{A^1}(\partial_0  \Sigma)$ to $M_{A^1}(\partial_1\Sigma)$, i.e., 
 $$\partial(M_{A^\bullet}(\Sigma))\,\,=\,\, M_{A^1}(\partial_1 \Sigma) \cdot M_{A^1}(\partial_0 \Sigma)^{-1}
 \in G^0.$$
 
 (b) $M_{A^\bullet}(\Sigma)$ is unchanged if $\Sigma$ is replaced by a thin homotopic
2-brane.

 (c) The correspondence
 $$\begin{gathered}
 x\in X= \Pi_0(X)\,\,\longmapsto \,\, pt \in \Ob \, 2\Cat(G^\bullet),\cr
 \gamma\in\Pi_1(X)\,\,\longmapsto \,\, M_{A^1}(\gamma)\in G^0 =  1\operatorname{Mor} (2\Cat(G^\bullet))
 \cr
 \Sigma\in\Pi_2(X)\,\,\longmapsto\,\, M_{A^\bullet}(\Sigma) \in 2\Hom_{2\Cat(G^\bullet))}
 (M_{A^1}(\partial_0\Sigma), M_{A^1}(\partial_1 \Sigma)
 \end{gathered}
 $$
 defines a 2-functor $\MM_{A^\bullet}: \Pi_{\leq 2}(X)\to 2\Cat(G^\bullet)$. 
   \end{thm}
       
   \noindent {\sl Proof:} (a) The condition of semiflatness, i.e., of vanishing of the fake curvature,
   means that $d_{\gen^\bullet}(A^2) = F_{A^1}$. Therefore 
   $$\partial\bigl(M_{A^\bullet}(\Sigma)\bigr) \,\,=\,\,\partial\biggl( M_{\oint_{A^1}(A^2)}(\sigma)\biggr)\,\,=\,\,
   M_{\oint_{A^1}(F_{A^1})}(\sigma)\,\,=\,\,M_B(\sigma),
   $$
  the last equality being the Schlessinger formula. So the statement follows from the fact
  that $B$ is the logarithmic differential of $M_{A^1}$. 
  
  \vskip .2cm
  
  (b) For a thin homotopy $[0,1]^3\to S\to X$ between two unparametrized 2-branes
  the pullback of the 3-curvature $F^3$ is the zero 3-form on $[0,1]^3$.
  So the covariant transgression of this pullback with respect to the pullback of $A^1$
  is zero as well. Our statement now follows from Proposition 
  \ref{prop:curv-tramsgression}. 
  
  \vskip .2cm
  
  (c) We need to prove compatibility of $\MM_{A^\bullet}$ with the two compositions
  $*_0, *_1$. Compatibility with $*_1$ follows from compatibility of the holonomy of
  usual connections with composition of paths. Indeed, if $\sigma, \sigma'$ are the
  paths in $P_x^y X$ corresponding to parametrized membranes $\Sigma, \Sigma'$
  s.t. $\partial_0\Sigma=\partial_1\Sigma'$, then 
$\Sigma *_1\Sigma'$ corresponds to the composition of $\sigma$ and $\sigma'$. 

Let us now prove that 
$\MM_{A^\bullet} (\Sigma *_0\Sigma') = \MM_{A^\bullet}(\Sigma) *_0
\MM_{A^\bullet}(\Sigma')$
 whenever $\Sigma,\Sigma'$ are $*_0$-composable.
 For this it is enough to assume that either $\Sigma$ or $\Sigma'$ is a
 1-morphism, i.e., a path considered as a membrane. Indeed, 
 in any 2-category the $*_0$-composition
of  any two 2-morphisms
can be expressed using $*_1$-composition of 2-morphisms as well as 
$*_0$-composition involving a 2-morphism and a 1-morphism. 

We now recall the definition of the $*_0$-composition of a 1-morphism and a 2-morphism
in $2\Cat(G^\bullet)$. Since $2\Cat(G^\bullet)$ has one object, any 1-morphism and
2-morphism can be $*_0$-composed in either order. Now, 1-morphisms are identified
with elements of $G^0$, and 2-morphisms with triples $T= (g_0, g_1, h)$ with
$g_i\in G^0, h\in G^{-1}$ such that $\partial(h) = g_0^{-1}g_1$. Such a $T$
is a 2-morphism $T: g_0\Rightarrow g_1$. With these conventions, if $u\in G^0$,
then
\be\label{eq;comp-crossed-2cat}
\begin{split}
u*_0 T = (ug_0, ug_1, h): ug_0\Rightarrow ug_1, 
\cr
 T*_0 u = \bigl(g_0u, g_1u, \beta(u^{-1})(h)\bigr):
g_0u\Rightarrow g_1u.
\end{split}
\ee
Let now assume that $\Sigma=\gamma\in P_y^z X$ is a path, and $\Sigma'$ is a parametrized
2-brane corresponding to a path $\sigma': [0,1]\to P_x^y X$. 
So we denote by $\sigma'_s = \Sigma(-, s) \in P_x^y X$ the path corresponding to $s\in [0,1]$. 
Then $\gamma*_0\Sigma'$
corresponds to a path in $[0,1]\to P_x^z X$ which takes $s\in [0,1]$
into the composition $\gamma*\sigma'_s$, so $\sigma'_s$ is run first and then followed by $\gamma$.
This means that   $\bigl(\oint_{A^1}(A^2)\bigl)_{\gamma*\sigma'_s}$, the
covariant transgression evaluated at $\gamma*_0\sigma'_s$ is equal to
$\bigl(\oint_{A^1}(A^2)\bigl)_{\sigma'_s}$
the result of evaluation along $\sigma'_s$ alone. Indeed, the integral 
in \eqref{eq:covariant-transgression} splits into the sum of two integrals
corresponding to $\sigma'_s$ and $\gamma$, and the second integral is equal to 0
because the pullback of the 2-form $A^2$ to the interval along $\gamma$ vanishes.
So   $\MM_{A^\bullet}(\gamma*_0\Sigma')$ as an element of $G^{-1}$ is equal to
$M_{A^\bullet}(\Sigma')$, in agreement with the first formula in \eqref{eq;comp-crossed-2cat}.

Assume now that $\Sigma'=\gamma'\in P_x^y X$ is a path, and $\Sigma$
  is a parametrized
2-brane corresponding to a path $\sigma:  [0,1]\to P_y^z X$. Then $ \Sigma*_0\gamma'$
corresponds to a path in $[0,1]\to P_x^z X$ which takes $s\in [0,1]$
into the composition $ \sigma_s*\gamma'$, so $\gamma'$ is run first and then followed by $\sigma_s$.
As before, the integral for $\bigl(\oint_{A^1}(A^2)\bigl)_{\sigma_s*\gamma'}$
splits into two integrals, one corresponding to $\sigma_s$, the other to $\gamma$.
As before, the integral corresponding to $\gamma'$ is equal to 0. But since $\sigma_s$
now
follows $\gamma$ in the order of integration, the integral corresponding to $\sigma_s$ will 
be equal to the
conjugation of $\bigl(\oint_{A^1}(A^2)\bigl)_{\sigma_s}$ by 
$M_{A^0}(\gamma')$, the holonomy of the connection $A^0$ along $\gamma'$.
This agrees exactly with the second formula in \eqref{eq;comp-crossed-2cat}. \qed

\subsection{2-dimensional holonomy and the crossed module of formal 2-branes.}
Consider the dg-Lie algebra $\fen^\bullet(\RR^n)_{\CM}$
 defined as  in \eqref{eq:cm-quotient}. It is situated in degrees $[-1,0]$. 
 We will denote this dg-Lie algebra
 as well as the corresponding crossed module of Lie algebras,  by $\gen^{\geq -1}_n$.
 The notation is chosen to indicate that  $\gen^{\geq -1}_n$ is the truncation
 of a longer (and more fundamental) crossed complex of Lie algebras, see
 \S \ref{subsec:form-R-n}. 
 Here is a concise summary of the properties of $\gen_n^{\geq -1}$.

 \begin{prop}\label{prop:g-n}
  (a) As a dg-Lie algebra, $\gen^{\geq -1}_n$ is
  is generated by symbols $Z_i, i=1, ..., n$,
 in degree 0 and $Z_{ij}$, $1\leq i<j\leq n$ in degree $-1$, subject to the relations:
 $$d(Z_{ij}) =[Z_i, Z_j], \quad [Z_{ij}, Z_{pq}]=0,\quad [[Z_i, Z_j], Z_{pq}] = [Z_{ij}, [Z_p, Z_q]].$$
 
 (b) The degree 0 part $\gen_n^0$ is the free Lie algebra $\FL(Z_1, ..., Z_n)$.
 As a Lie algebra with the bracket $[-,-]_{-1}$, the degree $-1$ part $\gen^{-1}_n$
 is isomorphic to a direct product of a free Lie algebra (isomorphic to $[\gen_n^0,\gen_n^0]$)
 and an abelian Lie algebra isomorphic to $\Gamma_2^{\cl}(\RR^n)$. 
 
 (c) In particular, $H^0(\gen^{\geq -1}_n)=\RR^n$ (the abelian Lie algebra),
 while $H^{-1}(\gen^{\geq -1}_n)=\Gamma_2^{\cl}(\RR^n)$. 
 
 \end{prop}
 
 \noindent {\sl Proof:}   (a) By definition, see  \eqref{eq:cm-quotient}, $\gen_n^{\geq -1}$
 is generated by the $Z_i$, $Z_{ij}$ as stated, with the differential defined as stated, but
 then quotiented out by all elements of degrees $\leq -2$ as well as by l elements
 of the form $[dx,y]-[x,dy]$ for $x,y$  arbitrary elements of degree $-1$. 
 Quotienting by all elements of degrees $\leq -2$ is equivalent to imposing
 the second series of relations in (a): that  $[Z_{ij}, Z_{pq}]=0$. Further, the last series of
 relations in (a) amounts to imposing the conditions $[dx,y]-[x,dy]=0$ for $x=Z_{ij}$,
 $y=Z_{pq}$ being the  degree $-1$ generators. Now, the Jacobi identity implies that
 any other element of the form $[dx,y]-[x,dy]$ will belong to the dg-Lie ideal  generated
 by the above particular elements, so defining $\gen_n^{\geq -1}$ by the relations as in (a),
 has the same effect as \eqref{eq:cm-quotient}. 
 
 (b) This was established in the course of the proof of Theorem \ref{thm:cartesian},
 see Eq. \eqref{eq:Lie-direct product-whitehead}.

 (c) follows from (b). 
 \qed

 \vskip .2cm

 Firther, the crossed module of Lie algebras 
 $\gen_n^{\geq -1}$  carries the universal translation invariant
 semiflat connection $A^\bullet$ with 
 $$A^1=\sum_i Z_i dt_i, \quad A^2=\sum_{i<j} Z_{ij} dt_idt_j.$$
 
 We denote by $\gen^{\geq -1}_{n,d}$ the quotient of $\fen^\bullet(\RR^n)_{\CM}$
 by the $(d+1)$th term of the lower central series. Thus
$\gen^0_{n,d}= \FL(\RR^n)/\FL_{\geq d+1}(\RR^n)$ is the free degree $d$ 
nilpotent Lie algebra on $n$ generators. We regard $\gen^{\geq -1}_{n,d}$
as a nilpotent crossed module of finite-dimensional Lie algebras.

\begin{ex} Let $d=2$. Then $\gen_{n,2}^0$ is the ``universal Heisenberg Lie algebra"
generated by $Z_1, ..., Z_n$ such that $\hbar_{ij}:=[Z_i, Z_j]$ are linearly independent central
elements whose span is naturally identified with $\Lambda^2(\RR^n)$.
 Further, the Lie algebra $\gen_{n,2}^{-1}$ has a basis formed  by the $Z_{ij}, i<j$,
 spanning $\Lambda^2(\RR^n)$,
and the $[Z_i, [Z_{pq}]]$, spanning $\RR^n\otimes\Lambda^2(\RR^n)$. The bracket
$[-,-]_{-1}$ on $\gen_{n,2}^{-1}$ vanishes, while the differential sends $Z_{ij}$ to $\hbar_{ij}$.
Thus $\Ker\{d: \gen^{-1}_{n,2}\to\gen^0_{n,2}\}$ is identified with $\RR^n\otimes\Lambda^2(\RR^n)$.

Note the difference with the full  dg-Lie algebra $\gen^{\geq -1}_n$, in which
\[
d[Z_i, Z_{pq}] = [Z_i, [Z_p, Z_q]] \neq 0. 
\]
In this case we need to recall that the space  spanned by the $[Z_i, [Z_{pq}]]$,
splits, $GL_n(\RR)$-equivariantly as 
\[
\RR^n \otimes\Lambda ^2 (\RR^n) \,\,= \,\, \Lambda^3(\RR^n) \oplus \Sigma^{2,1}(\RR^n)
\]
of which the second summand maps by $d$ isomorphically to $\FL_3(\RR^n)\subset \gen_n^0$
and the first summand is a part of $\Ker(d) = \Gamma_2^{\cl}(\RR^n)$.

\end{ex}

Let $G^{\geq -1}_{n,d} = \exp(\gen^{\geq -1}_{n,d})$ be the crossed module
of nilpotent Lie group corresponding to $\gen^{\geq -1}_{n,d}$. We have then the
2-dimensional holonomy functor
\be\label{eq:m-n-d}
\MM_{n,d}^{\leq 2}: \,\, \Pi_{\leq 2}(\RR^n) \lra 2\Cat(G^{\geq -1}_{n,d})
\ee
associated to the tautological semiflat connection with values in $\gen^{\geq -1}$. 
Consider the projective limit crossed module $\widehat G^{\geq -1} _n = 
\varprojlim G^{\geq -1}_{n,d}$. We call
$\widehat G^{\geq -1}_n$ the {\em crossed module of formal 2-branes} in $\RR^n$. 
By construction, the
zero part  $\widehat G_n^0$ is the free prounipotent group over $\RR$ with $n$ generators.

\begin{prop} 
We have $\pi_0(\widehat G^{\geq -1}_n)=\RR^n$, while
$\pi_1(\widehat G^{\geq -1}_n) = (\Omega^{2, \cl} (\RR^n))^*$ is the algebraic
dual of the space of closed polynomial 2-forms in $\RR^n$. 
\end{prop}

\noindent {\sl Proof:} This follows from Proposition \ref{prop:g-n}(c),
as $(\Omega^{2, \cl} (\RR^n))^*$ is the pro-finite-dimensional completion of 
$\Gamma_2^{\cl}(\RR^n)$.

\qed

\vskip .2cm
   
 By passing to the limit we get the functor
\be\label{eq:2-chen}
\widehat \MM_n^{\leq 2}  = \varprojlim_d \MM_{n,d}^{\leq 2}: \Pi_{\leq 2}(\RR^n) \lra 2\Cat(\widehat G^{\geq -1}_{n}).
\ee
This functor is translation invariant: a parallel shift of a membrane   does not affect the value  of the
functor $\widehat \MM_n^{\leq 2}$. 
On the level of 1-morphisms,  $\widehat \MM_n^{\leq 2}$  associates to each unparametrized path $\gamma$ in $\RR^n$,   
the noncommutative power series $E_\gamma(Z_1, ..., Z_n)$ from (0.2), which is a group-like element of
$\RR\langle \langle Z_1, ..., Z_n\rangle\rangle$ and so is  an element of  $\widehat G_n^0$. 

\begin{ques}\label{ques:2d}
 K.-T. Chen proved in \cite{chen-faithful} that the correspondence $\gamma\mapsto E_\gamma$
is faithful modulo translations. That is, for two piecewise smooth paths $\gamma$ and $\gamma'$ the equality
 $E_\gamma=E_{\gamma'}$  in $\widehat G_n^0$ 
 implies that $\gamma'$ differs from $\gamma$ only by a rank 1 homotopy (reparametrization + 
cancellations)  and an overall translation in $\RR^n$. 
 Can one generalize this to the functor $\widehat \MM_n^{\leq 2}$? 

\end{ques}

 \vfill\eject

\section{The crossed complex of formal $n$-branes.}\label{sec:crossed-complex-n-branes}

\subsection{Crossed complexes of groups and Lie algebras}\label{subsec:crossed-complexes}

\paragraph{A. Groups.}  The following definition is taken from \cite{brown-higgins};
we restrict to the case of crossed complexes consisting of groups and not groupoids,
as in  \cite{brown-higgins}.
  
 \begin{Defi} A  crossed complex (of groups) is 
 a sequence of groups and
homomorphisms
$$G^\bullet \,\, =\,\,  \bigl\{ \cdots  \buildrel \partial_{-3} \over\longrightarrow   G^{-2}\buildrel
\partial_{-2}\over\longrightarrow  G^{-1}\
\buildrel \partial_{-1}\over\longrightarrow G^0\bigr\},$$
 equipped with an action of $G^0$ on each
$G^i$, via $\beta_i: G^0\to\Aut(G^i)$ satisfying the following properties:\hfill\break
(a) The part $\bigl\{G^{-1}\buildrel \partial_{-1}\over\longrightarrow G^0\bigr\}$
is a crossed module of groups; in particular, $\beta_0$ is the action
of $G^0$ on itself by conjugation.

(b) For $i\leq -2$ the group $G^i$ is abelian, and the action $\beta_i$ factors through an action
of $H^0(G^\bullet) = \Coker(\partial_{-1})$ on $G^i$. 

(c) The maps $\partial_i$ commute with the action of $G^0$.

(d) The composition $\partial_i \partial_{i-1}$ is the trivial homomorphism. 

\end{Defi}

Note that for any crossed complex $G^\bullet$  and any $i\leq 0$ the truncation
$G^{\geq i}$ is still a crossed complex (a crossed module, if $i=-1$).

\begin{ex} \label{ex:crossed-complex-filtered} (a)  Let 
$$X_\bullet \,\, =\,\,  \bigl\{ x_0\in X_1\subset X_2 \subset ...\bigr\}$$
be a filtered topological space, with $x_0$ being a point. Then setting
$G^{-i} = \pi_i(X_i, X_{i-1}, x_0)$ to be the relative homotopy groups of $X_i$ modulo $X_{i-1}$
and $\partial_{-i}$ to be the boundary map on the relative homotopy groups, we get a crossed
complex. We denote this crossed complex by $G^\bullet(X_\bullet)$.
See \cite{brown-higgins}. 

By considering crossed complexes of groupoids  \cite{brown-higgins} one can generalize this
example to the case of an arbitrary filtered space. 

\vskip .2cm

(b) In particular, let $X=\TT^n$ be the $n$-dimensional torus with
the cubical cell decomposition, as in Example \ref{ex:cubic-2-branes}. 
We denote by  $\TT^n_{\leq i}$
the $i$-skeleton of this decomposition. This gives a filtered CW-complex
which we denote $\TT^n_\bullet$. 
   We will call $G^\bullet(\TT^n_\bullet)$ 
  the {\em crossed complex of cubical branes in} $\RR^n$. Its truncation in
  degrees $[-1,0]$ is the crossed module discussed in 
  Example \ref{ex:cubic-2-branes}.
  
  \end{ex}

\paragraph{B. Lie algebras.}
One easily extends the concept of a crossed complex to the Lie algebra case.

 \begin{Defi}\label{def:2.1.5}
 A crossed complex of Lie $k$-algebras consists of a sequence
of Lie $k$-algebras and homomorphisms
$$\Gg^\bullet \quad = \quad \bigl\{ \cdots  \buildrel d_{-3} \over\longrightarrow   \Gg^{-2}\buildrel
d_{-2}\over\longrightarrow  \Gg^{-1}
\buildrel d_{-1}\over\longrightarrow \Gg^0\bigr\}$$
equipped with an action of $\Gg^0$ on each
$\Gg^i$, via $\alpha_i: \Gg^0\to\Der(\Gg^i)$ satisfying the following properties:\hfill\break
(a) The part $\bigl\{\Gg^{-1}\buildrel d_{-1}\over\longrightarrow \Gg^0\bigr\}$
is a crossed module of Lie algebras; in particular, $\alpha_0$ is the adjoint action
of $\Gg^0$ on itself.\hfill\break
(b) For $i\leq -2$ the Lie algebra $\Gg^i$ is abelian, and the action $\alpha_i$ factors through an action
of $H^0(\Gg^\bullet) = \Coker(d_{-1})$ on $\Gg^i$. \hfill\break
(c) The maps $d_i$ commute with the action of $\Gg^0$.\hfill\break
(d) The composition $d_i d_{i-1}$ is the zero homomorphism.
\end{Defi}

Propositions \ref{prop:semiabelian-conditions} and  \ref{prop:crossed-Lie-dg}  
 easily imply the following.

 \begin{prop}\label{prop:2.1.6}
 Let $\Gg^\bullet$ is a crossed complex of Lie algebras.
Then defining the bracket $[x,y]\in \Gg^{-i-j}$ for $x\in\Gg^{-i}, y\in\Gg^{-j}$ by
$$[x,y] = \alpha_j(x)(y), \,\, i=0, \quad [x,y] = -\alpha_i(y)(x), \,\, j=0, \quad
[x,y]=0, \,\, i,j<0,$$ 
we make $\Gg^\bullet$ into a semiabelian dg-Lie algebra. This correspondence
establishes an equivalence between crossed complexes of Lie algebras on the one hand
and semiabelian dg-Lie algebras on the other hand. \qed
\end{prop}

In particular,  Example \ref{ex:crossed-complex-filtered}  (a)
is analogous to  the Lie-algebraic Example \ref{ex:semiabelian-filtered}.

 \begin{cor}\label{cor:2.1.7} Let $G^\bullet$ be a crossed complex of Lie groups,
and $\Gg^i = \Lie(G^i)$. Then $\Gg^\bullet$ is a crossed complex of Lie $\RR$-algebras
and thus has a natural structure of a semiabelian dg-Lie $\RR$-algebra. \qed
\end{cor}

\vfill\eject

\paragraph{C. Lower central series.} Let $G^\bullet$ be a crossed complex of groups,
with operations written multiplicatively.
Its {\em lower central series}  consists of crossed subcomplexes
\be\label{eq:LCS-Cros-com-gr}
\begin{gathered}
G^\bullet = \gamma_1(G^\bullet) \supset\gamma_2(G^\bullet)\supset\cdots\cr
\gamma_r(G^\bullet) \,\,=\,\,\bigl\{ \cdots\lra \gamma_r(G^0, G^{-2})\lra \gamma_r(G^0, G^{-1})
\lra\gamma_r(G^0)\bigr\}.
\end{gathered}
\ee
Here $\gamma_r(G^0, G^{i})\subset G^i$ is defined, for $i=-1$, by \eqref{eq:LCS-CM} and
\eqref{eq:LCS-CM-details} and for $i\leq -2$ it is defined as the subgroup in
the abelian group $G^i$ generated by
  elements of the form
 $\beta(z)(x)\cdot x^{-1}$ for $z\in G^0$ and  $x\in \gamma_r(G^0, G^i)$.
We say that $G^\bullet$ is {\em nilpotent}, if $\gamma_r(G^\bullet)=1$
(consists of 1-element groups) for some $r$. 

\vskip .2cm

Similarly, let  $\gen^\bullet$ be a crossed complex of Lie $k$-algebras,
 Its {\em lower central series}  consists of crossed subcomplexes of Lie algebras
\be\label{eq:LCS-Cros-com-LA}
\begin{gathered}
\gen^\bullet = \gamma_1(\gen^\bullet) \supset\gamma_2(\gen^\bullet)\supset\cdots\cr
\gamma_r(\gen^\bullet) \,\,=\,\,\bigl\{ \cdots\lra \gamma_r(\gen^0, \gen^{-2})\lra \gamma_r(\gen^0, \gen^{-1})
\lra\gamma_r(\gen^0)\bigr\}.
\end{gathered}
\ee
Here the $\gamma_r(\gen^0, \gen^{i})\subset \gen^i$ are defined, generalizing  
\eqref{eq:LCS-CM-Lie} from the case $i=-1$, by
\be
\gamma_{r+1}(\gen^0, \gen^i) \,\,=\,\,\alpha(\gen^0)\bigl(\gamma_{r}(\gen^0, \gen^i)\bigr).
\ee

 \begin{prop} The dg-Lie algebra corresponding to the
 crossed complex $\gamma_r(\gen^\bullet)$ is equal to the $r$th term
 of the lower central series of the dg-Lie algebra corresponding to $\gen^\bullet$. 
 \qed
  \end{prop}
  
  We say that  crossed complex of Lie algebras $\gen^\bullet$ is {\em nilpotent}, if
  $\gamma_r(\gen^\bullet)=0$ for some $r$, i.e., if $\gen^\bullet$ is nilpotent
  as a dg-Lie algebra.  As in the case of crossed modules (Proposition
  \ref{prop:engel-CM}), the Engel theorem implies the following: 
  
  \begin{prop} Let $\gen^\bullet$ be a finite crossed complex of finite-dimensional
  Lie algebras. Then $\gen^\bullet$ is nilpotent, if and only if $\gen^0$ is a nilpotent
  Lie algebra.\qed
  \end{prop}
  
  \paragraph{D.  Malcev theory for crossed complexes.} Let $\gen^\bullet$ be a
  finite crossed complex of finite-dimensional
  Lie algebras, or, what is the same, a finite-dimensional semiabelian dg-Lie algebra
  over $k$. Assume that $\gen^\bullet$ is nilpotent. 
  Applying the Malcev theory to each nilpotent Lie algebra $\gen^i$, we get a sequence
  of groups and homomorphisms
  \be
  \exp(\gen^\bullet) \,\,=\,\,\bigl\{ \cdots\lra \exp(\gen^{-2})
  \buildrel\exp(d_{-2})\over \lra\exp(\gen^{-1})
  \buildrel \exp(d_{-1})\over \lra\exp(\gen^0)\bigr\}.
  \ee
  Exponentiating the action of $\gen^0$ on $\gen^i$ by derivations, we get an action
  of $\exp(\gen^0)$ on $\exp(\gen^i)$ by automorphisms, and we see easily:
  
  \begin{prop}\label{prop:malcev-CC}
  With the actions described, $\exp(\gen^\bullet)$ is a nilpotent crossed complex
  of algebraic groups over $k$. We have, therefore, three equivalent categories,
  with equivalence between (ii) and (iii) given by the functor $\exp$: 
  
  (i) Finite-dimensional, nilpotent,  semiabelian dg-Lie $k$-algebras;
  
  (ii) Finite, nilpotent crossed complexes of finite-dimensional Lie $k$-algebras.
  
  (iii) Finite, nilpotent crossed complexes of algebraic groups over $k$.
   \qed
   \end{prop}
   
   \begin{rem} In \cite{getzler}, E. Getzler associated a group-like object to
   any nilpotent dg-Lie algebra situated in degrees $\leq 0$. However, in 
   general, his construction is far from being reduced to some data of algebraic
   groups and their homomorphisms. It is the assumption of being
   semiabelian which allows for such a reduction.

   \end{rem}

\subsection{Strict $n$-categories and $n$-groupoids.}\label{subsec:strict-n-cat}
  We recall the concept of
(small, strict, globular) $n$-categories, see \cite{street} for more background.

\paragraph{A. Globular sets.} 

 \begin{Defi}\label{def:globular-set} A globular set $Y_\bullet$ consists of sets
$Y_i, i\geq 0$ and maps $s_i, t_i: Y_{i+1}\to Y_i$, $\1_i: Y_i\to Y_{i+1}$
satisfying the identities:
$$s_is_{i+1} = s_it_{i+1}, \quad t_i s_{i+1} = t_it_{i+1}, \quad s_i \1_i = t_i \1_i = \id. $$
\end{Defi}

Elements of $Y_i$ are called $i$-{\em cells} of $Y_\bullet$, the maps $s_i$ and $t_i$
are called the $i$-{\em dimensional source and target maps} for $(i+1)$-cells,
and $\1_i$ is called the $i$-{\em dimensional identity map}. For a globular set $Y$ one defines
the $i$-dimensional source and target maps for $j$-cells, $j\geq i$, by
\be\label{2.2.2}s_i = s^j_i = s_i s_{i+1} ... s_{j-1},  \,\,\,\, t_i = t_i^j = t_it_{i+1} ... t_{j-1}:
C_j\to C_i, 
\ee
and the iterated identity maps
\be\label{2.2.3}
\1_i^j = \1_{j-1} \1_{j-2} ... \1_i: Y_i\to Y_j, \quad i\leq j.
\ee
In particular, for $i=j$ we set $s_i^i=t_i^i=\1_i^i=\Id$. 

By an $n$-{\em globular set} we will mean a datum as in Definition \ref{def:globular-set} but with $Y_i$
defined only for $i\leq n$. An equivalent point of view is to say that an $n$-globular set
is a globular set $Y_\bullet$ with all the $Y_j, j\geq n$, identified with $Y_n$ via
$ \1^j_n$. 

\begin{exa}\label{ex:branes-transgression}
 (a)  Let $\Glb_n(X)$
be the set of parametrized $n$- branes in $X$, see  Definition \ref{def:param-p-branes}. 
   The sets $\Glb_n(X)$, $n\geq 0$,  give rise to a  globular set $\Glb_\bullet(X)$
via the  maps $s_i, t_i, \1_i$ defined before  \eqref{eq:s_p}. That is
\[
\begin{split}
(s_{n-1}\Sigma)(a_1, ..., a_{n-1}) = \Sigma(a_1,  ..., a_{n-1}, 0), \cr
(t_{n-1}\Sigma)(a_1, ..., a_{n-1}) = \Sigma(a_1, ..., a_{n-1}, 1),\cr
(\1_n\Sigma)(a_1, ..., a_{n+1})= \Sigma(a_1,  ..., a_{n}), \cr
\Sigma\in\Glb_n(X), \quad a_i\in I = [0,1].
\end{split}
\]
In particular, the points given by the images of $\{0\}\times I^{n-1}$ and $\{1\}\times I^{n-1}$ 
are, in this notation,  $s_0\Sigma$ and $t_0\Sigma$, the $0$-dimensional source and
target of $\Sigma$. 

Note that  a differential $n$-form on $X$ can be integrated over a parametrized $n$-brane:
\be
\int_\Sigma\omega \,\,:=\,\,\int_{I^n} \Sigma^*(\omega).
\ee

 \vskip .2cm
 
 (b) For $x,y\in X$ let $\Glb_\bullet(X)_x^y$ be the globular subset in $\Glb_\bullet(X)$ consisting
 of $\Sigma$ such that $s_0\Sigma=x$ and $t_0\Sigma=y$. Each  parametrized $n$-brane
 $\Sigma\in\Glb_n(X)_x^y$ gives rise to a parametrized  $(n-1)$-brane
 $$\tau(\Sigma) \in\Glb_{n-1}(P_x^yX)$$
 in the path space, called the {\em transgression} of $\Sigma$. Explicitly, $\tau(\Sigma)(a_1, ..., a_{n-1})$
 is the path $I\to X$ sending $a$ into $\Sigma(a, a_1, ..., a_{n-1})$. This operation is
 compatible with transgression \eqref{eq:usual-transgression} of differential forms:  if $\omega\in\Omega^n_X$
 is an $n$-form on $X$, then
 \be
 \int_{\tau\Sigma} \oint (\omega) \,\,=\,\,\int_\Sigma\omega.
 \ee

\end{exa}

\paragraph{B. $n$-categories, $n$-groupoids and crossed complexes.}   

\begin{Defi}\label{def:n-cat}
An $n$-category is an $n$-globular set $C$ equipped with
 composition maps
$$*_i: C_j\times_{s_i,t_i} C_j\to C_j.$$
 These compositions should
satisfy the following properties:\hfill\break
(a) Associativity.\hfill\break
(b) $j$-morphisms from the image of $\1_i^j$ are units for $*_i$.\hfill\break
(c) The compositions are compatible with the source and target maps:
$$s_i(u*_j v) = s_i(u) *_j s_i(v), \quad t_i(u*_j v) = t_i(u)*_j t_i(v).$$
(d) 2-dimensional associativity: 
$$(x*_i y)*_j (z*_i t) = (x*_j z) *_i (y*_j t).$$
\end{Defi}

Elements of $C_i$ for an $n$-category $C$ are commonly called $i$-{\em morphisms} of $C$.
For $i=0$ we will call $0$-morphisms simply {\em objects}. 

\vskip .2cm

One can also define $n$-categories inductively, as small categories enriched in the
monoidal category  of small $(n-1)$-categories, starting with 0-categories which
are the same as sets. See \cite{street}, Thm. 1.5.  

\vskip .2cm

By an $n$-{\em groupoid} we mean an $n$-category with all $j$-morphisms invertible for all
the compositions $*_j$. 
 An $n$-{\em group} is, by definition, an $n$-groupoid with one object. 
 We have the following result of Brown and Higgins \cite{brown-higgins}
 which generalizes the case $n=2$ discussed in \S \ref{sec:crossed-modules-groups}.

 \begin{thm}\label{2.2.6} Let $C$ be an $n$-group. Set 
$$G^0 = C_1, \quad G^{-i} = \Ker\{ s_i: C_{i+1}\to C_i\} \subset C_{i+1}, \,\, i=1, ..., n-1,$$
and define $\partial_{-i}: G^{-i}\to G^{-i+1}$ to be the restriction of $t_i$.
Define the action $\beta_i$ of $G^0$ on $G^{-i}$ by
$$\beta_i(g_0)(g_{-i}) = \1_1^i(g_0)*_0 g_{-i} *_0 \1_1^i(g_0^{-1}), 
\quad g_0\in C_1, \, g_i\in\Ker(s_i)\subset C_{i+1}.$$
Then
 $G^\bullet$ with these actions is  a crossed complex.
This construction establishes an equivalence between $n$-groups and crossed
complexes of groups situated in degrees $[-n+1, 0]$. \qed
\end{thm}

\vskip .2cm

For future reference we recall the inverse equivalence as well, cf. 
\cite{brown-higgins}, \S 4. It associates to a crossed complex $G^\bullet$
in degrees $[-n+1, 0]$, an $n$-group which we denote $n\Cat(G^\bullet)$. 
  Set $n\Cat(G^\bullet)_0= \{\pt\}$,
and
\be\label{2.2.7}
n\Cat(G^\bullet)_m = G^{-m+1}\times ... \times G^{-1}\times G^0, \quad m=1, ..., n. 
\ee
Define $s_{m-1}, t_{m-1}: n\Cat(G^\bullet)_m\to n\Cat(G^\bullet)_{m-1}$ by
\be\label{2.2.8}
\begin{split}
s_{m-1}(g_{-m+1}, ..., g_0) = (g_{-m+2}, ..., g_0), \cr
t_{m-1}(g_{-m+1}, ..., g_0) = (g_{-m+2}\cdot \partial(g_{-m+1}), g_{-m+3}, g_{-m+4}, ..., g_0),
\end{split}
 \ee
and
$\1_m: n\Cat(G^\bullet)_m\to n\Cat(G^\bullet)_{m+1}$ by
$$\1_m(g_{-m+1}, ..., g_0) = (1, g_{-m+1}, ..., g_0).$$
Let now $x= (g_{-m+1}, ..., g_0)$ and $y=(h_{-m+1}, ..., h_0)\in n\Cat(G^\bullet)_m$ be given.
Then, for a given $i<m$,  the condition $s_i(x) = t_i(y)$ means that
$$g_\nu=h_\nu, \,\, \nu>-i, \quad g_{-i} = h_{-i}\cdot \partial (h_{-i+1}).$$
In this case we set
 \be\label{2.2.9}
 \gathered
 x *_i y = 
   \bigl(g_{-m+1}\cdot  h_{-m+1}, ..., g_{-i}\cdot h_{-i}, h_{-i+1}, ..., h_0\bigr), \quad  i>0,
 \\
 x*_0y = \biggl(g_{-m+1}\cdot (\beta_{-m+1}(g_0) (h_{-m+1})), g_{-m+2} \cdot (\beta_{-m+2}(g_0) (h_{-m+2})), ...,
 \\
...,  g_{-1}\cdot (\beta_{-1}(g_0) (h_{-1})), g_0h_0\biggr). 
 \endgathered
 \ee
These structures make $n\Cat(G^\bullet) $ into an $n$-group.

 \begin{rem}\label{rem:c-groupoids}
  More generally,  one can extend the above to an equivalence between arbitrary $n$-groupoids 
 and crossed complexes of
 groupoids \cite{brown-higgins} of length $n+1$.
 \end{rem}

\paragraph{C. The fundamental $n$-groupoid of a filtered space.}
Let
$$X_\bullet \,\, =\,\,  \bigl\{ X_0\subset  X_1\subset X_2 \subset ...\bigr\}$$
be a filtered topological space. The crossed complex $G^{\geq -n+1}(X_\bullet)$
(described in Example \ref {ex:crossed-complex-filtered} (a) for the case $X_0$ is
a point) gives, by the equivalence of crossed complexes and $n$-groupoids (Remark
\ref{rem:c-groupoids}),
a certain $n$-groupoid $\varpi_{\leq n}(X_\bullet)$. We recall here a more direct
description of $\varpi_{\leq n}(X_\bullet)$ following \cite{brown:new}. 

Consider the $p$-globe $\bigcirc^p\simeq D^p$ with its filtration by skeleta:
\[
\on{sk}_q\bigcirc^p \,\,=\,\, s_q \bigcirc^p \,\,\cup\,\, t_q \bigcirc^p. 
\]

\begin{Defi}
(1) A filtered singular $p$-globe in $X_\bullet$ is a morphism of filtered topological spaces
\[
\Sigma: \bigl\{ \on{sk}_0 \bigcirc^p, \on{sk}_1\bigcirc^p, \cdots, \bigcirc^p\bigr\} \lra
X_\bullet  =  \bigl\{ X_0\subset  X_1\subset X_2 \subset ...\bigr\}. 
\]
(2) A filtered homotopy between two filtered singular $p$-globes $\Sigma$ and $\Sigma'$ in $X_\bullet$,
is a map $\Xi: I\times\bigcirc^p\to X$ such that for each $b\in I$ the map
$\Xi(b, -):$ is a filtered singular $p$-globe in $X$, and
\[
\Xi(0,-) = \Sigma, \,\,  \Xi(1,-) =\Sigma'. 
\]
\end{Defi}
The set of filtered homotopy classes of filtered simgular $p$-globes in $X_\bullet$ will be denoted
$\varpi_p(X_\bullet)$. The maps $s_i,  t_i, \1_i$ from Example \ref{ex:branes-transgression}
 (a)  descend on the filtered homotopy classes
 and make the collection $\varpi_{\leq n}(X_\bullet) = \{\varpi_p(X_\bullet)\}_{p\leq n}$
 into an $n$-globular set. 
 
 \begin{thm}\label{thm:brown-fund}\cite{brown:new} The operations $*_i$ on $\Glob_p(X)$, $p\leq n$, $i\leq p-1$
 given by 
 \[
 (\Sigma*_i\Sigma')(a_1, \cdots, a_p) =\begin{cases}
 \Sigma(a_1, \cdots, a_i, 2a_{i+1}, a_{i+2}, \cdots, a_p),&\text{ if } a_{i+1}\leq 1/2,
 \\
  \Sigma(a_1, \cdots, a_i, 2a_{i+1}-1, a_{i+2}, \cdots, a_p),&\text{ if } a_{i+1}\geq 1/2
 \end{cases}
 \]
 whenever $s_i\Sigma=t_i\Sigma'$, descend to well defined operations on $\varpi_{\leq n}(X_\bullet)$
 which make it into an $n$-groupoid with the set of objects $X_0$.
  If $X_0=\{x_0\}$ is a point, then the crossed complex associated to the $n$-group
 $\varpi_{\leq n}(X_\bullet)$ by Theorem \ref{2.2.6}, is identified with the crossed complex
 from Example \ref{ex:crossed-complex-filtered}(a). 
 \end{thm}
 
 We refer to \cite{brown:new} for an outline of the proof. This proof, based on
 \cite{BH-colimit}, uses a certain cubical Kan fibration property which can be seen as a
   higher-dimensional analog
 of the construction \eqref{eq:filling} connecting the thin homotopic parts of the   boundaries
 of $\Sigma$ and $\Sigma'$
 to make the composition possible. 
 
 An alternative approach would be to first define the crossed complex of groupoids $G^\bullet(X_\bullet)$ 
 associated to
 $X_\bullet$ (as done in Example    \ref{ex:crossed-complex-filtered}(a) for $X_0=\{x_0\}$)
 and then use the general  (groupoid) form of Theorem \ref{2.2.6} to produce am $n$-groupoid
 $\varpi_{\leq n}(X_\bullet)$.

 \paragraph{D. The $n$-groupoid of unparametrized branes in a manifold.}
 Let $X$ be a $C^\infty$-manifold of dimension $n$ with a filtration 
 $$X_\bullet \,\, =\,\,  \bigl\{ X_0\subset  X_1\subset X_2 \subset ...\subset X \bigr\}$$
 by closed subsets (not necessarily submanifolds). Then we can form a $C^\infty$-version of
 the $n$-groupoid $\varpi_{\leq n}(X_\bullet)$.
 
 \begin{Defi}
 (1) A filtered parametrized $p$-brane in $X_\bullet$ is a parametrized $p$-brane $\Sigma: \bigcirc^p\to X$
 which take each $\on{sk}_i\bigcirc^p$ to $X_i$. We denote by $\Glb_p(X_\bullet)$ the set of
 filtered parametrized $p$-branes in $X_\bullet$. 
 
 (2) A filtered homotopy of parametrized filtered $p$-branes $\Sigma, \Sigma'$ is a $C^\infty$-homotopy
 $\Xi: I\times\bigcirc^p\to X$ in the sense of Definition \ref{def:hom-branes} such that 
 $\Xi(b, \on{sk}_i\bigcirc^p) \subset X_i$ for each $b\in I$ and $i\leq p$. We denote by $\varpi^{C^\infty}_p(X_\bullet)$ the
 set of filtered homotopy classes of filtered parametrized $p$-branes in $X_\bullet$. 
 \end{Defi}
 
 The following is a natural $C^\infty$-modification of Theorem \ref {thm:brown-fund}. The proof
 may be obtained by the same steps as   in \cite{brown:new} for the topological case.
 
 \begin{prop}
 The operations $*_i$  defined as in Theorem \ref{thm:brown-fund}, descend to give an $n$-groupoid
 structure on $\varpi^{C^\infty}_{\leq n}(X_\bullet) = \{\varpi_p^{C^\infty}(X_\bullet)\}_{p\leq n}$. \qed
  \end{prop}
  
  We now note that a thin homotopy between parametrized $p$-branes in $X$ can be seen as a version
  of a filtered homotopy but with respect to an ``indetermined" filtration $X_\bullet$ of $X$ such that
  $X_i$ has Hausdorff dimension $\leq i$, see Remark \ref{rems:thin}(b). We can therefore write the set
  of thin homotopy classes of parametrized $p$-branes as
  \[
  \Pi_p(X) \,\,=\,\,\varinjlim_{(X_\bullet)} \varpi^{C^\infty}_p(X_\bullet), 
  \]
  where the inductive limit is taken over the filtering poset  formed by  filtrations $X_\bullet$ such that the Hausdorff
  dimension of $X_i$ is $\leq i$. 
  (One can say that considering thin homotopies amounts to a differential-geometric analog of the
  skeletal filtration. )
  In this way we obtain:
  
  \begin{prop}
   The operations  $*_i$ in the $\varpi_{\leq n}^{C^\infty}(X_\bullet)$
  give rise to a well defined $n$-groupoid structure on $\Pi_{\leq n}X$. 
  \qed
 \end{prop}
 
 We will call $\Pi_{\leq n}X$ the {\em $n$-groupod of unparametrized branes} in $X$.

\subsection{Higher holonomy.}

\paragraph {A. Connections with values in crossed complexes.} 
Let $G^\bullet$ be a crossed complex of Lie groups situated in degrees $[-n+1,0]$, and let 
$\gen^\bullet$ be the corresponding crossed complex of Lie algebras.
Let us consider $\gen^\bullet$ as a
  semiabelian dg-Lie algebra. 
  Let $X$ be a $C^\infty$-manifold and 
 $A^\bullet\in(\Omega^\bullet_X\otimes \gen^\bullet)^1$ be a $\gen^\bullet$-valued differential form 
of total degree 1. Thus $A^\bullet$ has components $A^i\in\Omega^i_X\otimes\gen^{1-i}$, $i=1, ..., n+1$. 
We want to consider
$A^\bullet$ as a connection in the trivial $n$-bundle with structure $n$-group $G^\bullet$. 

We denote by 
\be
F^\bullet= d_{\DR}A^\bullet - {1\over 2}[A^\bullet, A^\bullet] \,\,
\in\,\,
(\Omega^\bullet_X\otimes \gen^\bullet)^2
\ee
 the curvature of $A^\bullet$. Thus $F^\bullet$
has components $F^i\in\Omega^i_X\otimes\gen^{2-i}$, $i=2, n+2$. 
We say that $A^\bullet$ is {\em semiflat}, if all $F^i=0$ for $i=2, ..., n+1$. At the level of
components, semiflatness means:
\be\label{eq:semiflat-p}
\begin{split}
F_{A^1}\,=\, d_{\gen^\bullet}(A^2),\cr
d_{\DR}A^i - [A^1, A^i] \,=\, d_{\gen^\bullet}(A^{i+1}), \quad i=2, ..., n.
\end{split}
\ee
Note that $d_{\DR}A^i - [A^1, A^i]=\nabla_{A^1}(A^i)$ is just the covariant differential
of $A^i$ with respect to the connection $A^1$.

\paragraph {B. Forms on the space of paths corresponding to a connection.}
We fix points $x,y\in X$ and consider the $\gen^{1-i}$-valued differential forms
$\oint_{A^1} A^i$ on $P_x^y X$, for $i=2, ..., n+1$. 

\begin{prop}
\label{prop:d-oint-A1-Ai}
 If $n\geq 2$ and $A^\bullet$ is semiflat, then:
 
 (a) The curvature of $\oint_{A^1}(A^2)$ is equal to $\oint_{A^1}(\nabla_{A^1} A^3)$. 
  
 (b) For $i=2, ..., n-1$ we have 
 $$d_{\DR}\oint_{A^1} (A^i) \,\,= \,\, \oint_{A^1} (\nabla_{A^1}A^{i+1}).$$
\end{prop}

\noindent {\sl Proof: } Part (a)
follows from Proposition \ref{prop:curv-tramsgression} applied to the truncated connection with
values in the crossed module $\gen^{\geq -1}$. Indeed, $\nabla_{A^1} A^3$ is the 3-curvature
of this truncated connection. 

Part (b)   follows from Proposition \ref{prop:nabla-covariant-iter-integral}(b) in exactly the same
way as Proposition \ref{prop:curv-tramsgression}. More precisely, we see, similarly to 
Eq. \eqref{eq:d-of-int-A2},  that 
\[
d_\DR \oint_{A^1} (A^i) \,\,=\,\, \oint_{A^1}  (\nabla_{A^1} A^{i+1}) \,\, + \,\, 
\biggl[ d_\gen \oint_{A^1} (A^2), \,\, \oint_{A^1} (A^{i+1}) \biggr],
\]
and the last commutator is equal to 0. Indeed,     part (c) of Definition
\ref{def:2.1.5} of crossed complexes of Lie algebras says that $[d_\gen x,y]=0$
for any $x\in\gen^{-1}, y\in\gen^{-i}$. 
\qed

\paragraph{C. The $p$-dimensional holonomy.} 
Let $G^\bullet$ and $\gen^\bullet$ be as before, and $A^\bullet$ be a semiflat form with
values in $\gen^\bullet$ of total degree $1$.  Let $p\geq 1$ and let $\Sigma: I^p\to X$ be
a parametrized $p$-brane in $X$, with $s_0\Sigma=x$ and $t_0\Sigma=y$.
Recall that $\tau(\Sigma)$  denotes the $(p-1)$-brane in $P_x^yX$
obtained as the transgression of $\Sigma$, see Example \ref{ex:branes-transgression} (b). 
We now define the {\em $p$-dimensional holonomy}
of $A^\bullet$ along $\Sigma$ to be the element
\be
M_{A^\bullet}(\Sigma) \,\,=\,\,\begin{cases}
M_{A^0}(\Sigma) \,\,=\,\, P\exp\int_\Sigma A^1\,\,\in\,\, G^0, \quad {\rm if} \quad p=1;\cr
M_{A^{\leq -1}}(\Sigma)\,\,=\,\, P\exp \int_{\tau(\Sigma)} \oint_{A^1}(A^2) \,\,\in\,\, G^{-1},
\quad {\rm if} \quad p=2;\cr
\exp\int_{\tau(\Sigma)}\oint_{A^1} (A^{p-1})\,\,\in\,\, G^{-p+1}, 
\quad {\rm if} \quad p\geq 3. 
\end{cases}
\ee
Here the first line is the usual holonomy of the connection $A^1$ along the path $\Sigma$,
and the second line is the 2-dimensional holonomy of the truncated connection $A^{\leq 2}$
with values in the crossed module $\gen^{\geq -1}$, as defined by Baez-Schreiber
\cite{baez-schreiber}, 
see \eqref{eq:def-2-hol}. In the third line, which extends the definitions of 
\cite{baez-schreiber} to arbitrary crossed complexes, $\exp: \gen^{-p+1}\to G^{-p+1}$
is the exponential map of the abelian Lie group $G^{-p+1}$. We then define
\be\label{eq:MM-coordinates}
\gathered
\MM_{A^\bullet}(\Sigma) \,\,=\,\,\bigl( M_{A^\bullet}(\Sigma), M_{A^\bullet}(s_{p-1}\Sigma), ... , 
M_{A^\bullet}(s_1\Sigma)\bigr) \,\, \in
\\
\in \,\,\,  G^{-p+1} \times G^{-p+2} \times ... \times G^0 \,\,=\,\, 
n\Cat(G^\bullet)_p. 
\endgathered
\ee

\begin{thm}\label{thm:p-holonomy} 
(a) The element $\MM_{A^\bullet}(\Sigma)$ depends only on the thin
homotopy class of $\Sigma$. In particular, it is invariant under reparametrizations of
$\Sigma$ identical near the boundary.

(b) The correspondence $\Sigma\mapsto \MM_{A^\bullet}(\Sigma)$
defines  a homomorphism (strict $n$-functor)  $\MM_{A^\bullet}: \Pi_{\leq n}(X)\to n\Cat(G^\bullet)$. 
\end{thm}

\noindent {\sl Proof:} (a) Follows from Proposition \ref{prop:d-oint-A1-Ai}. Let us prove path (b).
The first thing to prove is compatibility of $\MM_{A^\bullet}$ with the source and target maps.
It amounts to the equality, for each $p$ and for each $p$-brane $\Sigma$:
\[
\partial (M_{A^\bullet}(\Sigma)) \,\,=\,\, M_{A^\bullet}(t_{p-1}\Sigma) \cdot M_A(s_{p-1}\Sigma)^{-1} \,\,
\in \,\, G^{-p+2}. 
\] 
For $p=1,2$ this equality has already been proved in Theorem \ref{thm:2-mon}(a). 
 For $p>2$ this follows from the equality, in $\gen^{2-p}$, of the logarithms of the both sides which is established
 using the condition of  semiflatness and the Stokes formula:
 \[
 \gathered
d_\gen \int_{\tau(\Sigma}\oint_{A^1} (A^{p-1}) \,\, = \,\,  \int_{\tau(\Sigma}\oint_{A^1}(d_\gen A^{p-1})
\,\, \buildrel \eqref{eq:semiflat-p}\over = \,\, \int_{\tau(\Sigma} \oint_{A^1} \bigl( d_\DR A^{p-2} - [A^1, A^{p-2}]\bigr) = 
\\ = \int_{\tau(\Sigma)} \oint_{A^1} \bigl (\nabla_{A^1}(A^{p-2})\bigr) \,\,\buildrel { \eqref{prop:d-oint-A1-Ai}(b)  }
\over =\,\,\int_{\tau(\Sigma} d_\DR \oint_{A^1} (A^{p-2}) \,\,\buildrel {\text{Stokes}}\over = \,\, \int_{\partial \tau(\Sigma)}
\oint_{A^1} (A^{p-2}) \,\,= 
\\
\hskip -8cm = \,\, \log M_A(t_{p-1}\Sigma) -  \log M_A(s_{p-1}\Sigma). 
 \endgathered
 \]
The second thing to prove is  that $\MM_{A^\bullet}$ commutes with the compositions:
\be\label{eq:MM-compatible}
\MM_{A^\bullet}(\Sigma*_i \Sigma') \,\,=\,\, \MM_{A^\bullet}(\Sigma)  *_i \MM_{A^\bullet}(\Sigma'),
\ee
whenever $\Sigma$ and $\Sigma'$ are $*_i$-composable. Denote
\[
p=\dim(\Sigma),\,\, p'=\dim(\Sigma'), \,\, \mu=\min(p, p'), \,\, m=\max(p, p'). 
\]
We use induction on $m$, the case $m\leq 2$ being already established in Theorem 
\ref{thm:2-mon}(c). 
Further, it is enough to prove \eqref{eq:MM-compatible} under the assumption that $i=\mu-1$. 
 Indeed, 2-dimensional associativity implies that 
such ``irreducible"  types of compositions generate all other compositions of all higher 
morphisms in any strict $n$-category. Keeping this assumption, we
distinguish two cases and analyze, in each case, the nature of the $(m-1)$-membrane $\tau(\Sigma*_i\Sigma')$
in the path space.

\vskip .2cm

\noindent {\bf Case 1: $\mu = 1$.} That is, either $p=1$ or $p'=1$, and so $i=0$. We can suppose $m\geq 3$. 
Suppose $p=1$. This means that 
$\Sigma=\gamma\in P_y^z X$ is a path and $s_0\Sigma'=x, t_0\Sigma'=y$ (so that  $\gamma*_0\Sigma'$ is
defined). In this case we have that
\[
\tau(\gamma*_0\Sigma') \,\,=\,\, l_\gamma(\tau(\Sigma'))
\]
is the left translation of $\tau(\Sigma')$ by $\gamma$, and so Proposition \ref{prop:trans-trans}(a)
implies that 
\[
M_{A^\bullet}(\gamma*_0\Sigma') \,\,=\,\, M_{A^\bullet}(\Sigma') \,\,\in G^{-p'+1},
\]
which is the highest (degree $p'$) component of the desired equality \eqref{eq:MM-compatible}
with respect to the decompostion \eqref{eq:MM-coordinates}. 
 The  validity of the lower components follows by induction since
\[
s_{p'-1}(\gamma*_0\Sigma') = \gamma*_0 s_{p'-1}(\Sigma').
\]
The case when $p'=1$ is treated similarly using  Proposition \ref{prop:trans-trans}(b).

\vskip .2cm

\noindent {\bf Case 2: $\mu\geq 2$.} We can also assume $m\geq 3$.
 In this case all three branes $\Sigma, \Sigma'$ and $\Sigma*_i\Sigma'$,
$i=\mu-1$, have the same $s_0=x$ and $t_0=y$.
 Therefore $\tau(\Sigma)$, $\tau(\Sigma')$ and
$\tau(\Sigma*_i\Sigma')$ are branes in the same path space $P_x^yX$. and
\be\label{eq:add-int-sigma}
\tau(\Sigma*_i\Sigma')\,\, =  \,\, \tau(\Sigma) *_{i-1}\tau(\Sigma'),
\ee
 is identified with  the $*_{i-1}$-composition (which  can be seen, 
 geometrically, as the union) of $\tau(\Sigma)$ and $\tau(\Sigma')$
inside $P_x^yX$.

 We first concentrate on the highest
(degree $m$) component of  \eqref{eq:MM-compatible}. If $p=p'$, then the validity of this compoment
follows from additivity of the integral of the $(p-1)$-form $\oint_{A^1} (A^{p-1})$ with respect to the
decomposition \eqref{eq:add-int-sigma} of the domain of integration. The validity of lower components follows
by induction since, under our assumptions of $p=p'$ and $i=p-1$ we have 
\[
s_{p-1}(\Sigma *_i \Sigma') =s_{p-1}(\Sigma)*_{i-1} s_{p-1}(\Sigma'). 
\]
so for the $s_{p-1}$-parts we have the same situation of equal dimension and ``irreducible" composition. 

If $p\neq p'$, say $p<p'$, then $m=p'$ and the integral of $\oint_{A^1} (A^{m-1})$  over $\tau(\Sigma*_i\Sigma')$
is equal to the integral over $\tau(\Sigma')$ since the other part has smaller dimension. This equality is
precisely the degree $m$ component of \eqref{eq:MM-compatible}. The equality of components of smaller
degree follows by induction similarly to the above since
\[
s_{p'-1}(\Sigma *_i\Sigma') = \Sigma  *_i s_{p'-1}(\Sigma'), \]
 so for the $s_{p-1}$-parts we again  have the   situation of   ``irreducible" composition.
 
The case $p>p$ is treated similarly. \qed

\subsection{The crossed complex of formal branes in $\RR^n$ and higher holonomy.}
\label{subsec:form-R-n}

Fix $n\geq 1$ and denote by $\gen^\bullet_n$ the semiabelian dg-Lie algebra $\fen^\bullet(\RR^n)_\sab$,
see \S \ref{subsec:semiab-f}. It is situated in degrees $[-n+1, 0]$. 
We
will also consider $\gen^\bullet_n$  as a crossed complex of Lie algebras, alternating the two points of view. 

Let $\gen^\bullet_{n,d} = \gen^\bullet_n/\gamma_{d+1}(\gen^\bullet_n)$
denote  the semiabelian dg-Lie algebra (as well as the corresponding crossed complex
of Lie algebras) obtained by quotienting $\gen^\bullet_n$ by the $(d+1)$st
layer of the lower central series \eqref{eq:LCS-Cros-com-LA}.  As $\gen^\bullet_{n,d}$ is
nilpotent and finite-dimensional, the Malcev theory for crossed complexes
(Proposition \ref{prop:malcev-CC}) produces a crossed complex of unipotent Lie groups
$G^\bullet_{n,d}= \exp(\gen^\bullet_{n,d})$.
The universal translation invariant connection \eqref{eq:the-connection}  descends to
a translation invariant connection $A^\bullet_{n,d}$  on $\RR^n$ with values in $\gen^\bullet_{n,d}$. 
  This connection gives rise to the higher holonomy $n$-functor
  \be\label{eq:M-n-d}
  \MM_{n,d} = \MM_{A^\bullet_{n,d}}:\,\,  \Pi^{\leq n}(\RR^n) \lra n\Cat(G^\bullet_{n,d}).
  \ee
  We further consider the pro-unipotent crossed complex $\widehat G^\bullet_n = \varprojlim_d G^\bullet_{n,d}$
  which we call 
    the
{\em crossed complex of formal branes} in $\RR^n$. By passing to the projective limit of the $M_{n,d}$,
we get the $n$-functor
\be\label{eq:M-n-formal}
\widehat \MM_n \,\,=\,\,\varprojlim_d \MM_{n,d}: \,\, \Pi^{\leq n}(\RR^n) \lra n\Cat(\widehat G^\bullet_{n}).
\ee

\begin{ques}
Similarly to Question \ref{ques:2d}, it is interesting to understand to what extend is the functor
$\widehat \MM_n$ injective on higher morphisms. Consider in particular the sets
\[
\pi_i^{\text{geom}}(\RR^n, x) \,\, := \,\, \bigl\{ \Sigma\in \Pi^{\leq n} (\RR^n)_i \,  \bigl| s_{i-1}(\Sigma) = t_{i-1}(\Sigma) = x\bigr\} ,
\quad i\geq 1. 
\]
As for any strict $n$-groupoid, these sets are groups, abelian for $i\geq 2$. Translation of membranes
 in $\RR^n$ identifies them
for different $x\in\RR^n$, so we can assume $x=0$.

For any $i$, elements of $\pi_i^{\text{geom}}(\RR^n, 0)$ can be understood as geometric unparametrized
spheres (maps $\Sigma: (I^i, \partial I^i)\to (\RR^n, 0)$ modulo thin homotopies). 

The set  $\pi_1^{\text{geom}}(\RR^n, 0)$ is the group of  ``paths in $\RR^n$  modulo reparemetrization and
cancellation" beginning  at 0. Chen's theorem \cite{chen-faithful} implies that $\widehat \MM_n$
embeds it into the pro-unipotent completion of the free group on $n$ generators.

For $i\geq 2$, the functor $\widehat \MM_n$  sends any  $\Sigma\in \pi_i^{\text{geom}}(\RR^n, 0)$ to its associated current (functional on polynomial $i$-forms)
\[
I_\Sigma: \,\, \omega \,\,\mapsto \int_\Sigma \omega \,\, := \,\, \int_{I^i} \Sigma^*(\omega). 
\]
and the question is to what extent the thin homotopy class of $\Sigma$ can be recovered from $I_\Sigma$.

\end{ques}

\vfill\eject

\appendix

\section{Appendix: forms and connections on the space of paths.}
Here we collect the necessary facts about ``covariant"  generalizations of Chen's iterated integrals,
as developed in \cite{hofman, baez-schreiber}. We give a direct treatment,   
reducing these facts to classical results of Chen \cite{chen1}.

 \subsection {Scalar iterated integrals.}\label{subsec:scalar-integrated}
    Let $X$ be a $C^\infty$-manifold
 and $PX=C^\infty([0,1], X)$ be the space of smooth paths in $X$. 
 Suppose we are given differential forms 
 $\omega_\nu$, $\nu=1, ..., r$ 
on $X$ of degree $m_\nu+1$.  Chen's {\em iterated integral} of $\omega_1, ..., \omega_r$ is the form
$\oint (\omega_1, ..., \omega_r)$ on $PX$ of degree $m_1+...+m_r$, sending
tangent vectors $\delta_1\gamma$, ...,  $\delta_{m_1+...+m_r}\gamma$ at a point $\gamma$
to
\be\label{eq:iter-int-scal}
\begin{split} 
\operatorname{ALT} \int_{0\leq t_1\leq ...\leq t_r\leq 1} \prod_{\nu=1}^r
\omega_\nu \biggl(\gamma(t_\nu); \gdot(t_\nu),
 \delta_{m_1+...+m_{\nu-1}+1}\gamma(t_\nu),
...,  \delta_{m_1+...+m_{\nu}}\gamma(t_\nu) \biggr)
\cr
 \times dt_1 ... dt_r.
  \end{split}
\ee
Here $\operatorname{ALT}$ stands for alternation in the arguments $\delta_i\gamma$. 
Put differenty, let 
 \be
\Delta_r\,\,=\,\,
\{ 0\leq t_1\leq ...\leq t_r\leq 1\}
\ee
be the integration domain (simplex) in
\eqref{eq:iter-int-scal}, and consider the diagram
\be
P X\buildrel \pi_r\over\longleftarrow \Delta_r\times P X \buildrel e_r\over\lra X^r, \quad
e_r( t_1, ..., t_r, \gamma) = (\gamma(t_1), ..., \gamma(t_r)),
\ee
with $\pi_r$ being the projection. Then
\be
\oint(\omega_1, ..., \omega_r) \,\,= \,\,
(\pi_r)_* \, e_r^* (\omega_1\boxtimes  ...\boxtimes \omega_r),
\ee
where $(\pi_r)_*$ in the integration along the fibers of $\pi$, and 
$\omega_1\boxtimes  ...\boxtimes \omega_r$ is the wedge product of the pullbacks of the
forms $\omega_\nu$ from the factors. Cf.  \cite{getzler-jones-petrack}, \S 2.

One can replace $PX$ in the above definition  by the space of
Moore paths $C^\infty([a,b], X)$ for any interval $[a,b]\subset\RR$.
We denote the corresponding forms on this space by the same
symbol $\oint (\omega_1, ..., \omega_r)$.

\begin{ex}\label{ex:iter-integrals-chen}
 (a) Let $r=1$. Then  the correspondence which takes $\omega$ to $\oint (\omega)$ 
 restricted to the subspace $P_x^y X\subset PX$,  is the transgression
\eqref{eq:usual-transgression}
of differential forms. 

\vskip .2cm 

(b) Suppose that all
 $m_\nu=0$. Then $\oint (\omega_1, ..., \omega_r)$ is a function on $PX$. On the other hand,
 consider the associative algebra $\RR\llan Z\rran = \RR\llan Z_1, ..., Z_r\rran$
 of noncommutative formal power series in variables $Z_1, ..., Z_r$.
 Consider the differential form
 $
 A\,=\,\sum \omega_i Z_i\,\in\,\Omega^1_X\otimes \RR\llan Z\rran
 $
  and the corresponding
 connection $\nabla = d_{\DR}-A$ on $X$ with values in $\RR\llan Z\rran$. 
 Since $\RR\llan Z\rran$ is a projective limit of finite-dimensional associative algebras,
 this connection has a well defined holonomy which is a function
 $M_A=M_\nabla: PX\to \RR\llan Z\rran^\times$. The following fundamental formula of
 Chen can be used as an alternative definition of iterated integrals of 1-forms:
 \be\label{eq:iterated-integrals-holonomy}
 M_A \,\,=\,\,\sum_{d=0}^\infty \sum_{1\leq i_1, ..., i_d \leq r} \,\, Z_{i_1} ... Z_{i_p} \, \oint (\omega_{i_1}, ..., \omega_{i_d}).
 \ee

 \vskip .2cm
 
 (c) More generally, let $m_\nu$ be arbitrary. Following \cite{chen1}, \S 2, we consider first the
 algebra $\RR\la Z\ra = \RR\la Z_1, ..., Z_r\ra$ of noncommutative polynomials in the $Z_\nu$,
 equipped with grading $\deg(Z_\nu) =-m_\nu$. Then denote by $\RR\llan Z\rran$ the
 completion of $\RR\la Z\ra$ allowing power series only in those generators which have
 degree $0$. Then $\RR\llan Z\rran$ inherits the grading, and 
  \[
 A\,=\,\sum \omega_i Z_i\,\in\,(\Omega^\bullet_X\otimes \RR\la Z\ra)^1
 \]
 is an element of total degree $1$ called the {\em formal power series connection}
associated to the forms $\omega_1, ..., \omega_r$ of arbitrary degrees.
Now, defining $M_A$ by
 the infinite sum  as in \eqref{eq:iterated-integrals-holonomy}, we find that
 $M_A$ is an invertible element of $\Omega^\bullet_{PX}\otimes \RR\llan Z\rran$
 of total degree 0. 
 
 \vskip .2cm
 
 (d) Similarly, for any interval $[a,b]\subset \RR$, the sum as in  
  \eqref{eq:iterated-integrals-holonomy} defines an  invertible element  
  \[
  M_{A, [a,b]} \,\,\in\,\, (\Omega^\bullet_{C^\infty([a,b], X)} \otimes\RR\llan Z\rran)^0. 
  \]
 
\end{ex}

\vskip .2cm

We now recall the well known result of Chen, see \cite{chen1}, \S 4, 
  describing the products and differentials of iterated integrals.  
As usual, by an $(r,s)$-{\em shuffle} we mean a permutation  $w\in S_{r+s}$ such that
$w(i)<w(j)$ for all $i<j$ such that ether both $i,j\leq r$ or both $i,j \ge r$. 
The set of $(r,s)$-shuffles will be denote by $\Sh(r,s)$. 

\begin{prop}\label{prop:iter-int-properties} Let $x,y\in X$ and consider the
iterated integrals as differential forms on the subspace $P_x^y X\subset PX$. Then:

(a) The wedge product of  two iterated integrals 
is given by 
\[
\begin{gathered}
\oint (\omega_1, ..., \omega_r)  \,\,\,\wedge\,\,\, \oint (\omega_{r+1}, ...\omega_{r+s})\,\,=\,\,
\sum_{w\in \Sh(r,s)} (-1)^{\delta(w)}
\oint(\omega_{w(1)}, ..., \omega_{w(r+s)}),\cr
{\rm where}\quad \delta(w) \,\,=\,\,\sum_{i<j\atop w(i)>w(j)} m_i m_j.
\end{gathered}
\]

(b) The exterior differential of an iterated integral is given by
\[
\begin{split}
d \oint (\omega_1, ..., \omega_r)
\,\,=\,\,
-\sum_{\nu=1}^r 
(-1)^ {\sum_{\nu'<\nu}m_\nu'}
 \oint (\omega_1, ..., d\omega_\nu, ..., \omega_r) -
 \cr
 -\sum_{\nu=1}^{r-1}
 (-1)^{\sum_{\nu'\leq \nu}(m_\nu'+1)}
 \oint (\omega_1, ..., \omega_\nu \wedge\omega_{\nu+1}, ..., \omega_r).
 \end{split}
\]
\end{prop}

\noindent    Note that part (a) is obvious: it follows by decomposing the product of
two simplices $\Delta_r\times\Delta_s$ into $r+s\choose r$ simplices labelled by shuffles. 
The integral over each of these simplices is the iterated integral as stated. 

On the contrary, the standard proofs of part (b) are usually much less transparent. 
For future purposes, let us indicate a simple proof in the case when all $m_\nu=0$. 
  In this case we use Chen's formula \eqref{eq:iterated-integrals-holonomy}
for the holonomy of the connection $A$. The curvature of  $A$
is then  
\be\label{eq:F-A-iterated}
F\,=\, dA - {1\over 2} [A,A]\,\,=\,\,
\sum d\omega_\nu\cdot  Z_\nu \,- \,\sum_{\nu'<\nu} (\omega_{\nu'}\wedge\omega_\nu)\cdot
[Z_{\nu'}, Z_\nu].
\ee
After that part (b) is obtained by
applying the Schlessinger formula  \eqref{eq:form-B} to $A$, writing $dM_A = M_A\cdot B$, using part (a) to
evaluate the product and comparing coefficients
at the monomial $Z_1 ... Z_r$.  

\vskip .3cm

We now formulate a generalization
 of the Schlessinger formula  to higher degree forms, found by Chen.
 This is achieved by expressing  the original formula in a more diagrammatic fashion
 and then replacing functions by differential forms. 
  We use the conventions  of 
 Example \ref{ex:iter-integrals-chen}(c) and consider $A$ as a graded connection
 with values in the associative algebra $\RR\la Z\ra$ (considered as a Lie algebra).
 The curvature 
\[
F\in (\Omega^\bullet_{P_x^yX} \otimes \RR\la Z\ra )^2
\]
 of $A$ is then given by the same formula \eqref{eq:F-A-iterated},
 with the commutators understood in the graded sense. 
  For $t_0\in I$ consider the diagram
 \[
 P_x^yX \buildrel j_{t_0}\over\lra I\times P_x^y X \buildrel e_1\over\lra X,
 \quad j_{t_0}(\gamma) = (t_0, \gamma), \, e_1(t,\gamma)=\gamma(t),
  \] 
 and define
 \be\label{eq:contraction-F}
 {\bf i}(F, t_0) \,\, =\,\, j_{t_0}^*\bigl( i_{\partial\over \partial t} e_1^* F \bigr)\,\,\in\,\,
 (\Omega^\bullet_{P_x^yX}\otimes \RR\la Z\ra)^1. 
 \ee
Here $i_{\partial\over\partial t}$ is the contraction of differential forms with the vector
field $\partial\over\partial t$ on $I\times P_x^yX$.
 Further, for $t\in I$ consider the map  
  \[
\pi_{t}:  P_x^y X\lra C^\infty([0,t], X)
 \]
 given by restricting  paths
 $\gamma: I\to X$ to the subinterval $[0,t]$, and let
 \be\label{eq:M-A-leq-t}
  M_{A, \leq t} \,\,=\,\,\pi_t^* M_{A, [0,t]} \,\,\in\,\, (\Omega^\bullet_{P_x^yX} \otimes \RR\llan Z\rran )^0,
\ee
 see Example \ref{ex:iter-integrals-chen}(d).

\begin{thm}[\cite{chen1}, Thm. 2.3.2]\label{thm:chen-schlessinger}
 In this situation we have
\[
M_A\wedge d_\DR M_A\,\,=\,\,\int_{0}^{1} M_{A, \leq t}^{-1} \wedge {\bf i}(F, t)\wedge M_{A, \leq t} 
\cdot dt. \qed
\]
\end{thm}

One can deduce from this the general case of Proposition 
\ref{prop:iter-int-properties}(b) in the same way as indicated above for the case when
all $m_\nu=0$.

\begin{rems} (a) The notation $\oint$, different from Chen's but compatible with \cite{hofman, baez-schreiber},
is chosen to emphasize that we work with forms on $P_x^yX$, the subspace of paths 
having source and target fixed. In particular, if $x=y$, then the paths are closed. 

\vskip .2cm

(b) Let $(\Ac^\bullet, d)$ be a commutative dg-algebra. Its {\em bar-construction} is
a new commutative dg-algebra $\Bar^\bullet(\Ac^\bullet, d)$ which, as a graded vector space, is
 $\bigoplus_r \Ac[-1]^{\otimes r}$, the tensor algebra of the graded space $\Ac^\bullet$ 
 with degree shifted by $(-1)$.
 This space is
  equipped with the 
{\em shuffle multiplication}
\be
\begin{gathered}
(\omega_1\otimes \cdots\otimes\omega_r) \cdot (\omega_{r+1}\otimes \cdots\otimes \omega_{r+s}) 
\,\,=\,\, 
\sum_{w\in \Sh(r,s)} (-1)^{\delta(w)}\,\,\omega_{w(1)}\otimes \cdots \otimes \omega_{w(r+s)}\cr
 \delta(w) \,\,=\,\,\sum_{i<j\atop w(i)>w(j)} m_i m_j, \quad \omega_i\in\Ac^{m_i+1},
\end{gathered}
\ee
and with the differential
\be
\begin{split}
d(\omega_1\otimes\cdots\otimes \omega_r) \,\,= \,\,
-\sum_{\nu=1}^r 
(-1)^ {\sum_{\nu'<\nu}m_\nu'}
 \,\, \omega_1\otimes  ...\otimes  d\omega_\nu \otimes ...\otimes  \omega_r -
 \cr
 -\sum_{\nu=1}^{r-1}
 (-1)^{\sum_{\nu'\leq \nu}(m_\nu'+1)}
 \,\, \omega_1\otimes ...\otimes  (\omega_\nu \cdot \omega_{\nu+1})\otimes  ...\otimes  \omega_r.
 \end{split}
\ee
Proposition \ref{prop:iter-int-properties} can be formulated by saying that iterated integrals define a morphism
of commutative dg-algebras
\be
\oint: \Bar^\bullet(\Omega^\bullet_X, d_{DR}) \lra (\Omega^\bullet_{P_x^y X}, d_{DR}).
\ee
The main result of Chen \cite{chen1} is that $\oint$ is a quasi-isomorphism, and so 
$\Bar^\bullet(\Omega^\bullet_X, d_{DR})$ calculates the
coholomogy of $P_x^yX$. 
\end{rems}

\begin{rem}
Another approach to Theorem \ref{thm:chen-schlessinger} would be to deduce it from
the classical  Schlesinger formula   \eqref{eq:form-B} but applied to {\em supermanifolds},
see \cite{manin} for background. 
 More precisely, for any $C^\infty$-manifold $X$ we have the supermanifold
\[
\Sc X \,\,=\,\, \operatorname{Spec} (\Omega^\bullet_X) \,\,=\,\, \underline\Hom(\RR^{0|1}, X),
\]
see \cite{kontsevich, kapranov-vasserot}. Thus differential forms on $X$ can be seen as functions on $\Sc X$.
We can then extend the concept of a differentiable space  and the functor $\Sc$ to the  super-situation
 and, in particular, have 
the identification
\[
\Sc P_x^y X \,\,=\,\, P_x^y \Sc X \quad \subset \quad \Sc PX = \underline\Hom(\RR^{0|1} \times [0,1], X) = 
P\Sc X. 
\]
Further, by analyzing differential forms on $\Sc X$, i.e., functions on $\Sc\Sc X$,
as in \cite{kapranov-vasserot},  one can interpret $p$-forms on $X$, $p\geq 2$,  as certain 1-forms
on $\Sc X$. In this way iterated integrals of forms of higher degree on $X$ can be expressed
through iterated integrals of 1-forms on $\Sc X$ so  the Schlesinger formula  applied to
$\Sc X$, implies Theorem \ref{thm:chen-schlessinger}. These remarks  also apply to more sophisticated
iterated integrals considered further in this Appendix. 
\end{rem}

 \subsection{ Iterated integrals with values in an associative algebra.}
 Let $V_0, ..., V_r$
 be finite-dimensional $\RR$-vector spaces, and let 
 $\mu: V_1\otimes \cdots \otimes V_r\to V_0$ be a multilinear map. Given
 differential forms $\omega_\nu\in\Omega^{m_\nu+1}_X\otimes V_{\nu}$,
$i=1, ..., r$, we have the iterated integral
\be
\oint (\omega_1, ..., \omega_r)\,\,=\,\,
\oint^\mu  (\omega_1, ..., \omega_r)\,\,\,\in\, \,\,\Omega^{\sum m_\nu}_{PX}\otimes V_{0}.
\ee
It is obtained by tensoring the usual (scalar) iterated integral map
\[
\oint: \,\bigotimes_{n=1}^r{} \Omega^{m_\nu+1}_X \lra \Omega^{\sum m_\nu}_{PX}
\] 
(tensor product over $\RR$), with $\mu$. 

\vskip .2cm

Next, let $R^\bullet = \bigoplus_{d\leq 0} R^d$ be a $\ZZ_-$-graded associative $\RR$-algebra,
with each $R^d$ finite-dimensional. Let $A^\bullet\in (\Omega^\bullet_X\otimes R^\bullet)^1$
be a graded connection on $X$ with values in $R^\bullet$ (considered as a graded Lie algebra in
a standard way). As usual, we write $A^p\in\Omega^p_X\otimes R^{1-p}$ for the
component of $A^\bullet$ which is a $p$-form. 
Let $F^\bullet $ be the curvature of $A^\bullet$.
 Let $\mu_r: (R^\bullet)^{\otimes r}\to R^\bullet$ be the $(r-1)$-fold
product map. Consider the {\em Picard series}
\be\label{eq:picard-series}
M_A \,\,=\,\,\sum_{r=0}^\infty \oint^{\mu_r}
\underbrace{ (A^\bullet, \cdots, A^\bullet)}_r
 \,\,\in\,\, (\Omega^\bullet_{PX}
\otimes R^\bullet)^0.
\ee
Unlike the series \eqref{eq:iterated-integrals-holonomy} which is purely formal,
the Picard series involves actual infinite summation of real numbers because
the component $A^1$ can enter the iterated
integrals arbitrarily many times without raising the degree of the  resulting form on $PX$. 

\begin{prop} The Picard series converges to an invertible element of 
$(\Omega^\bullet_{PX}
\otimes R^\bullet)^0$.
\end{prop}

\noindent {\sl Proof:} The degree $0$ component $\sum_r \oint^{\mu_r}(A^1, ..., A^1)$
represents the standard Picard (ordered exponential) series for the holonomy of the
usual connection $A^1$, so it converges to an invertible
element of $\Omega^0_{PX}\otimes R^0$.  As for higher degree components, the question
reduces to the convergence of series like
\[
\sum_{N_1, ..., N_s=0}^\infty \oint^{\mu_r} 
\bigl(A^{p_1}, (A^1)^{N_1}, A^{p_2},  (A^1)^{N_2}, ..., A^{p_s},  (A^1)^{N_s}\bigr), \quad
r=s+\sum N_i,
\]
where $p_1, ..., p_s\geq 2$ are fixed and the notation $(A^1)^N$ stands for the fragment $A^1, ..., A^1$ ($N$ times). 
This convergence follows from that of the series for the holonomies of the connections in
the trivial bundles with fibers $R^{1-p_i}$ induced by $A^1$ via the left and right actions
of $R^0$ on $R^{1-p_i}$. \qed 

\vskip .2cm

As in \S \ref{subsec:scalar-integrated}, we denote by $M_{A, [a,b]}$ the  form
on the Moore path space $C^\infty([a,b], X)$ defined in the same way but with iterated integrals
taken over $[a,b]$. We define
\[
{\bf i}(F^\bullet, t)\in (\Omega^\bullet_{PX}\otimes R^\bullet)^1, \quad M_{A, \leq t}\in 
 (\Omega^\bullet_{PX}\otimes R^\bullet)^0, \quad t\in I,
\]
using the same formulas as in \eqref{eq:contraction-F}, \eqref{eq:M-A-leq-t}.

\begin{prop}\label{prop:chen-schlessinger-assoc}
 Fix $x,y\in X$. Then 
we have the equality of forms on $P_x^y X$:
\[
M_A\wedge d_\DR M_A\,\,=\,\,\int_{0}^{1} M_{A, \leq t}^{-1} \wedge {\bf i}(F, t)\wedge M_{A, \leq t} 
\cdot dt\,\,\,\in\,\,\,(\Omega^\bullet_{P_x^yX}\otimes R^\bullet)^1.
\]
\end{prop}

\noindent{\sl Proof:} Taking bases in $R^{1-p}$, $1\leq p\leq\dim(X)$, we can reduce the statement
to Theorem \ref{thm:chen-schlessinger}. That is, we can find forms $\omega_\nu\in\Omega^{m_\nu+1}_X$,
$\nu=1, ..., r$ with the corresponding grading on the free algebra $\RR\la Z\ra$ and
the corresponding graded connection $\widetilde A^\bullet = \sum Z_\nu \cdot\omega_\nu$, 
together with a homomorphism
$\phi: \RR\la Z\ra\to R^\bullet$ such that $\phi_*(\widetilde A^\bullet)=A^\bullet$
and therefore
  $\phi_*(F_{\widetilde A}) = F^\bullet$.  Further, applying $\phi$ to convergent infinite sums,
  we see that $\phi_*(M_{\widetilde A})=M_A$
etc. so the statement follows from the cited  theorem. Alternatively,
one can repeat, in the present context,
 the proof of  Theorem \ref{thm:chen-schlessinger}  given by Chen,  by using
the differential equation of  $M_{A, \leq t}$ with respect to $t$. \qed

\subsection{Covariant iterated integrals.}
Let $E_0, ..., E_r$ be $C^\infty$ real vector bundles on $X$, with connections
$\nabla_0, ..., \nabla_r$, and let $\mu: E_1\otimes ... \otimes E_r\to E_0$
be a morphism of vector bundles preserving the connections. 
We write $\nabla$  for the system $(\nabla_\nu)_{\nu=0}^r$. Let 
 $x,y\in X$ be two points. 
 Given
 differential forms $\omega_\nu\in\Omega^{m_\nu+1}_X\otimes E_{\nu}$,
$\nu=1, ..., r$,  we define, following 
 \cite{hofman, baez-schreiber},  the {\em (left)  covariant iterated integral},
 which is a form 
\be
\loint_\nabla (\omega_1, ..., \omega_r)\,\,=\,\,
\loint_\nabla^\mu  (\omega_1, ..., \omega_r)\,\,\,\in\, \,\,\Omega^{m_1+...+m_r}_{P_x^y X}\otimes E_{0,x}.
\ee
The value of this form at a path $\gamma$ and tangent vectors $\delta_i\gamma$,
$i=1, ..., m_1+...+m_r$, is defined to be
\be
\begin{gathered}
 \operatorname{ALT} \int_{0\leq t_1\leq ...\leq t_r\leq 1}  dt_1 ... dt_r \cdot 
 \mu\biggl\{
 \bigotimes_{\nu=1}^r\,\, \biggl[ 
M_{\nabla_\nu}(\gamma_{\leq t_\nu})^{-1}\biggl(
 \omega_\nu \bigl(\gamma(t_\nu); \cr
\gdot(t_\nu),
 \delta_{m_1+...+m_{\nu-1}+1}\gamma(t_\nu),
\cdots,  \delta_{m_1+...+m_{\nu+1}}\gamma(t_\nu) \bigr)\biggr)
\biggr]\biggr\}.
 \end{gathered}
\ee
The idea is simple: in order to be able to form the tensor product, we transport the value of each $\omega_\nu$
from the fiber of $E_\nu$ over  $\gamma(t_\nu)$, to the fiber over $x=\gamma(0)$. 
In the case when all the bundles are trivial and $\nabla_\nu=d-A_\nu$ we write $A$ for the system $(A_0, ..., A_r)$
are $\loint_A$ for $\loint_\nabla$. 

We can also transport the values to the fibers over $y=\gamma(1)$, to define
the {\em right covariant iterated integral} to be the form
\be
 \varointclockwise_\nabla^\mu (\omega_1, ... , \omega_r) \,\,\,:=\,\,\, M_{\nabla_0} \loint_\nabla^\mu 
  (\omega_1, ... , \omega_r)\, \,\,\in\,\,\,
  \Omega^{m_1+...+m_r}_{P_x^y X}\otimes E_{0,y}.
\ee
(Note the  difference between the symbols $\loint$ and $\roint$ for the left and right covariant iterated
integrals.)
We will use the notation
\be\label{eq:(())}
\loint_\nabla ((\omega_1, ..., \omega_r))\,\,:=\,\,\loint_\nabla^{\Id}(\omega_1, ..., \omega_r), \quad E_0=E_1\otimes ...
\otimes E_r
\ee
for the covariant iterated integral corresponding to $\mu=\Id$, and similarly for $\roint_\nabla 
((\omega_1, ..., \omega_r))$.

\begin{ex}\label{ex:two-trans}
The covariant transgression introduced in \S \ref{subsec:2-dim-hol} C, is a particular case of the construction
\eqref{eq:(())}
for $r=1$.
More precisely, let $G$ be a Lie group, $Q$ a principal $G$-bundle on $X$ and  $\alpha: G\to GL(V)$   a 
representation of $G$. A connection $\nabla$ in $Q$ induces a connection $\nabla_V$ in the vector bundle
$V(Q)$. In this notation, for each $\Phi\in\Omega^{m+1}_X\otimes V(Q)$ we have 
\[
\oint_\nabla (\Phi) \,\,=\,\, \loint_{\nabla_Q} ((\Phi)). 
\]
  \end{ex}

We now formulate the main properties of covariant iterated integrals, 
 cf.  \cite{hofman}, formula (15) and \cite{baez-schreiber}, Prop.2.7.

\begin{prop}\label{prop:nabla-covariant-iter-integral}
(a)  The tensor/wedge product of two covariant iterated integrals with $\mu=\Id$ is found by:
\[
\loint_\nabla ((\omega_1, ..., \omega_r)) \,\,\,\wedge\,\,\,
 \loint_\nabla ((\omega_{r+1}, ..., \omega_{r+s})) \,\,=\,\,
 \sum_{w\in \Sh(r,s)} (-1)^{\delta(w)} w^*
\loint_\nabla ((\omega_{w(1)}, ..., \omega_{w(r+s)})),
 \]
 where $\delta(w)$ is as in Proposition \ref{prop:iter-int-properties}(a), and 
 \[
 w^*:  E_{w(1),x} \otimes ... \otimes E_{w(r+s),x} \lra
\bigl (E_{1,x} \otimes ... \otimes E_{r,x}\bigr) \otimes
\bigl (E_{r+1,x}\otimes ... \otimes E_{r+s, x}\bigr)
 \]
 is the permutation map corresponding to $w$. 

\vskip .2cm

(b) For  any $\mu: E_1\otimes...\otimes E_r\to E_0$ as before,
the de Rham differential of a covariant iterated integral 
(as a form on $P_x^yX$ with values in $E_{0,x}$) is found by 
\[
\begin{split}
d_{\DR}  \loint_\nabla^\mu (\omega_1, ..., \omega_n)\,\,=\,\,
-\sum_{\nu=1}^r 
(-1)^ {\sum_{\nu'<\nu}m_\nu'}
 \loint_\nabla^\mu (\omega_1, ..., \nabla_\nu \omega_\nu, ..., \omega_r) -
 \cr
 -\sum_{\nu=1}^{r-1}
 (-1)^{\sum_{\nu'\leq \nu}(m_{\nu'}+1)}
 \loint _\nabla^{\mu'} (\omega_1, ..., \omega_\nu \wedge\omega_{\nu+1}, ..., \omega_r)\cr 
+\sum_{\nu=1}^r  (-1)^{\sum_{\nu' < \nu} (m_{\nu'}+1)}
 \loint_\nabla^{\mu''_\nu} (\omega_1, ..., \omega_{\nu-1}, F_{\nabla_\nu}, \omega_{\nu}, ..., \omega_n). 
\end{split} 
\]
Here $F_{\nabla_\nu}\in\Omega^2_X\otimes \End(E_\nu)$ is the curvature of $\nabla$, while
\[
\mu'_\nu: E_1\otimes ...\otimes (E_\nu\otimes E_{\nu+1})\otimes ...\otimes E_r \lra E_0
\]
is obtained from $\mu$ using the associativity of tensor product, and
\[
\mu''_\nu: E_1\otimes ...\otimes E_{\nu-1}\otimes \End(E_\nu) \otimes E_\nu\otimes ...\otimes E_r\lra E_0
\]
is obtained by composing $\mu$ with the action map $\End(E_\nu)\otimes E_\nu\to E_\nu$. 

\end{prop}

\noindent {\sl Proof:}   (a)  This follows, just like the corresponding part of Proposition
\ref{prop:iter-int-properties}, by decomposing the product of
two simplices $\Delta_r\times\Delta_s$ into $r+s\choose r$ simplices labelled by shuffles. 

\vskip .2cm

(b) First, it is enough to prove the statement under the assumption that
 $E_0=E_1\otimes...\otimes E_r$ and $\mu=\Id$.
Indeed, the general statement follows from that one by applying 
a given $\mu: E_1\otimes ...\otimes E_r\to E_0$ to the above ``universal" statement. 
So the rest of the proof will be under  the assumption that $\mu=\Id$. 

Second, it is enough to assume that each bundle $E_\nu$ is trivial. 
 Indeed, we can embed each $E_\nu$  a direct sum of a trivial
bundle: $E_\nu \oplus E'_\nu= X\times V_\nu$, where $V_\nu$
is a finite-dimensional $\RR$-vector space. We can 
 choose a connection
on $E'_\nu$ in an arbitrary way thus getting a connection
on $X\times V_\nu$. After this, each form $\omega_\nu$
can be considered as a form with values in $V_\nu$, and the statement
follows if we know it for the trivial bundles $X\times V_\nu$, by applying
the projection $V_{0}=E_{0,x}\oplus E'_{0,x}\to E_{0,x}$. 
So 
the rest of the proof will be under the additional assumption that $E_\nu= X\times V_\nu$
for each $\nu$.

Consider the tensor algebra
\be\label{eq:tensor-algebra}
T \,\,=\,\, \bigoplus_{d=0}^\infty (V_1\oplus ...\oplus V_r)^{\otimes d} \,\,=\,\,
\bigoplus_{d=0}^\infty \bigoplus_{\nu_1, ..., \nu_d=1}^r
  Z_{\nu_1} ... Z_{\nu_d}\cdot
 (V_{\nu_1}\otimes ... \otimes V_{\nu_d}).
\ee
Here the monomials in noncommuting variables $Z_1, ..., Z_r$ are added as a bookkeeping device,
and to emphasize the analogy with Example \ref{ex:iter-integrals-chen}(b-c). Thus
an element of $T$ can be seen as a noncommutative polynomial in the $Z_i$, but  with the coefficient
at each monomial lying in the corresponding 
tensor product of vector spaces as specified. We introduce a grading on $T$
by putting $\deg(Z_\nu)=1-m_\nu$. 
 We also consider the completion $\widehat T$
of $T$ with respect to the powers of the ideal generated by the $Z_\nu$ of degree 0. 
We view $\widehat T$ as the algebra of formal series in the $Z_\nu$, with coefficient at each monomial being
as above. For instance, if each $V_\nu= \RR$, then  
$T= \RR\la Z\ra$ is 
 the algebra of noncommutative polynomials, and $\widehat T =\RR\llan Z\rran$, 
 as in Example \ref{ex:iter-integrals-chen}(c). Let $\End(T)$ be the
 algebra of endomorphisms of $T$ as an $\RR$-vector space, and $\End(\widehat T)$
 be the algebra of continuous endomorphisms of $\widehat T$ as a topological
 vector space. 

 We have a natural connection in the trivial bundle $X\times T$ on $X$,
 induced by the $\nabla_\nu$ on $E_\nu=X\times V_\nu$, and
 denoted by $\nabla'$. It is compatible with the direct sum decomposition 
 in the RHS of \eqref{eq:tensor-algebra}  as well as with the tensor product. 
Let write the connection $\nabla_\nu$ as $\nabla_\nu = d_\DR-A'_\nu$, with
$A'_\nu\in\Omega^1_X\otimes \End(V_\nu)$  and write accordingly
$\nabla' = d_\DR-A'$. Using the terminology of \S \ref{subsec: connection},
we view $A'$ as a graded connection with values in the graded associative algebra
$\End(T)$ (considered as a graded Lie algebra in the usual way). 

In addition, we have a 1-form
\be
A'' \,\,=\,\,\sum_{\nu=1}^r  \omega_\nu  \cdot Z_\nu\,\,\in\,\, (\Omega^\bullet_X\otimes T)^1
\,\,\subset\,\,(\Omega^\bullet_X\otimes\End(T))^1,
\ee
where the embedding on the right is via action of $T$ on itself by left multiplication. 
We use it to construct a new graded connection 
 $A^\bullet = A' + A''{}^\bullet$ on $X$ with values
 in $T$,
which is no longer compatible with the product. 
 Let 
  $\End(\widehat T)^{\leq 0}\subset \End(\widehat T)$
  be the subalgebra of endomorphisms which take any homogeneous element  
  into a (possibly infinite) sum of elements of equal or lesser degrees.
  Then  $\End(\widehat T)^{\leq 0}$
   can be represented as a projective
limit of $\ZZ_-$-graded associative algebras $R^\bullet$ with all graded 
components finite-dimensional. So applying the Picard series 
\eqref{eq:picard-series} to such algebras $R^\bullet$ and taking the projective 
limit, we associate to $A^\bullet$ an
invertible  element $M_A\in (\Omega^\bullet_{P_x^yX}\otimes \End(\widehat T))^0$. 

 \begin{lem}  Denoting by  $1_T\in T$  the unit element,  we have:
 $$
M_{A}(1_T)  \,\,=\,\,\sum_{d=0}^\infty \sum_{\nu_1, ..., \nu_d=1}^r 
Z_{\nu_1}...Z_{\nu_d}\cdot \roint_\nabla ((\omega_{\nu_1}, ..., \omega_{\nu_d}))\,\,\,\in \,\,\, 
\Omega^\bullet_{P_x^yT}\otimes \widehat T. 
$$
Here in the RHS we use right covariant iterated integrals.
\end{lem}

\noindent {\sl Proof:} This follows directly from comparison with the Picard series. \qed

\vskip .2cm

Notice  that the curvature of the graded connection $A$ has the form
\be
F_{A}\,\,\,=\, \,\, F_{\nabla'}\,\, -\,\,\sum_{\nu' < \nu} [Z_{\nu'}, Z_\nu]\cdot \omega_{\nu'}\wedge \omega_\nu,
\ee
Notice further that the  Schlessinger-Chen formula 
 (Proposition \ref{prop:chen-schlessinger-assoc}) is applicable to our 
situation since it is a projective limit of situations considered
 in  that  proposition. So we write
  \be\label{eq:formula-d-wedge}
  d_\DR M_A= M_A\wedge B, 
  \ee
  where $B$ is given by the integral in
 the Schlessinger-Chen formula.  We then  compare coefficients in 
 \eqref{eq:formula-d-wedge} at
 $Z_1...Z_r$ and use
  use the (already proved)
 part (a) of Proposition \ref{prop:nabla-covariant-iter-integral}
 to express the products of iterated integrals appearing on the right. Note that 
 the multiplication by $M_A$ in the RHS of \eqref{eq:formula-d-wedge} accounts for the difference
 between $\loint$ and $\roint$. 
  This established part (b) of the Proposition.

 \vfill\eject

\vskip .3cm

   Kavli IPMU, 5-1-5 Kashiwanoha, Kashiwa, Chiba, 277-8583 Japan,

    {{\tt mikhail.kapranov@ipmu.jp}}

\ed

\end{document}